\documentclass{amsart}
\usepackage[utf8]{inputenc}
\usepackage[T1]{fontenc}
\usepackage{lmodern}
\usepackage{amsmath}
\usepackage{adjustbox}
\usepackage{amssymb}
\usepackage{stmaryrd}
\usepackage{mathtools}
\usepackage{latexsym}
\usepackage{graphicx}
\usepackage[shortlabels]{enumitem}
\usepackage{nicefrac}
\usepackage{microtype}
\usepackage{braket}
\usepackage{color} 
\usepackage{tikz-cd}
\usepackage{quiver}
\usepackage{MnSymbol}
\usepackage{enumitem} 
\setlist[enumerate,1]{label=\textup{(\arabic*)}}

\usepackage[all]{xy}
\newdir{ >}{{}*!/-5pt/@{>}}

\usepackage[pdftitle={Bornological K-theory},
pdfauthor={Devarshi Mukherjee},
pdfsubject={Mathematics}
]{hyperref}
\usepackage{cleveref}

\newcommand{\comment}[1]{}  
\newcommand{\edit}[1]{\marginpar{\footnotesize{#1}}}

\theoremstyle{plain}
\newtheorem{theorem}{Theorem}[section]
\newtheorem*{theorem*}{Theorem}
\newtheorem*{definition*}{Definition}
\newtheorem{lemma}[theorem]{Lemma}

\newtheorem{corollary}[theorem]{Corollary}
\newtheorem{proposition}[theorem]{Proposition}
\theoremstyle{remark}
\newtheorem{remark}
[theorem]{Remark}
\newtheorem{notation}[theorem]{Notation}

\theoremstyle{definition}
\newtheorem{definition}[theorem]{Definition}

\newtheorem{example}[theorem]{Example}
\numberwithin{theorem}{section}

\newcommand\C{\mathbb C}

\newcommand\N{\mathbb N}

\newcommand\R{\mathbb R}
\newcommand{\s}{\mathcal{S}}
\newcommand{\cF}{\mathcal{F}}
\newcommand{\KU}{\mathbf{KU}}

\newcommand{\cW}{\mathcal{W}}
\newcommand\Z{\mathbb Z}

\newcommand{\Fun}{\mathbf{Fun}}

\newcommand{\Ind}{\mathbf{Ind}}

\newcommand{\bD}{\mathbf{D}}
\newcommand{\bC}{\mathbf{C}}

\newcommand{\CMR}{\mathrm{CMR}}
\newcommand{\BAlg}{\mathbf{BAlg}}
\newcommand{\NBA}{\mathbf{NBAlg}}

\newcommand{\LBAlg}{\mathbf{LBAlg}}

\newcommand{\kk}{\mathbf{kk}^{\mathrm{CMR}}}
\newcommand{\bKK}{\mathbf{KK}^{\mathrm{born}}}
\newcommand{\nkk}{\mathbf{kk}^{\mathrm{nuc}}}

\newcommand{\KKzero}{\mathrm{KK}^\CMR}

\newcommand{\bkk}{\mathbf{kk}^{\mathrm{born}}}
\newcommand{\kkzero}{\mathrm{kk}^\CMR}
\newcommand{\BAlgk}{\BAlg_{\kappa}}
\newcommand{\LBAlgk}{\LBAlg_{\kappa}}
\newcommand{\cInd}{\mathbf{Ind}_{\aleph_0}}

\newcommand{\bdd}{\mathcal B}
\newcommand{\vN}{\mathsf{vN}}
\newcommand{\Cpt}{\mathsf{Cpt}}
\newcommand{\coma}{\widehat}

\newcommand{\cS}{\mathcal{S}}

\newcommand{\Bor}{\mathsf{Born}}
\newcommand{\CBorn}{\mathsf{CBorn}}

\newcommand{\defeq}{\mathrel{:=}} 

\newcommand*{\into}{\rightarrowtail}
\newcommand*{\onto}{\twoheadrightarrow}

\newcommand{\colim}{\mathrm{colim}}

\DeclarePairedDelimiter{\abs}{\lvert}{\rvert}
\DeclarePairedDelimiter{\norm}{\lVert}{\rVert}
\DeclarePairedDelimiter{\gen}{\langle}{\rangle}
\DeclarePairedDelimiterX{\setgiven}[2]{\{}{\}}{#1\,{:}\,\mathopen{}#2}

\newcommand{\rig}{\mathrm{rig}}
\newcommand\haotimes{\mathbin{\coma{\otimes}}}

\newcommand*{\an}{\mathrm{diff}}
\newcommand{\diff}{\mathrm{an}}

\newcommand{\Sp}{\mathbf{Sp}}

\newcommand{\K}{\mathbf{K}}

\newcommand{\nil}{\mathrm{nil}}

\newcommand{\cK}{\mathcal{K}}
\newcommand{\KH}{\mathbf{KH}}
\newcommand{\KD}{\mathbf{KD}}
\newcommand{\cJ}{\mathcal{J}}

\newcommand{\op}{\mathrm{op}}
\newcommand{\ev}{\mathrm{ev}}
\newcommand{\dvr}{V}
\newcommand{\Alg}{\mathbf{Alg}}

\newcommand{\nKK}{\mathbf{KK}^{\mathrm{nuc}}}

\DeclareMathAlphabet{\mathpzc}{OT1}{pzc}{m}{it}

\newcommand{\csK}{\mathcal{K}_{C^*}}
\newcommand{\sK}{\mathcal{K}_{\mathcal{S}}}

\newcommand{\Hom}{\mathbf{Hom}}

\newcommand{\KK}{\mathbf{KK}^{\mathrm{CMR}}}

\usepackage[backend=biber]{biblatex}
\begin{document}

\title{Bivariant \(K\)-theory for bornological algebras}

\author{Anupam Datta}
\email{anupam.datta@mathematik.uni-goettingen.de}

\author{Devarshi Mukherjee}
\email{devarshi.mukherjee@abdn.ac.uk}

\begin{abstract}
 We construct a presentable, stable \(\infty\)-category \(\bKK\), and exhibit it as the target of the universal functor from the category of bornological algebras into cocomplete, stable \(\infty\)-categories, preserving filtered colimits, and satisfying smooth homotopy invariance, excision, and stabilisation with respect to Schatten-class or smooth compact operators. We also construct a bivariant \(K\)-theory \(\nKK\) for suitably nuclear bornological algebras, and show that it is a right Bousfield localisation of \(\bKK\). The category \(\nKK\) is symmetric monoidal and admits a bivariant Chern character to local cyclic homology. The categories \(\nKK\) and \(\bKK\) corepresent Corti\~nas and Thom's diffeotopy \(K\)-theory. We use the Chern character to show that up to a suitable Schatten-class stabilisation, the complexified Chern character from the diffeotopy \(K\)-theory (resp. \(K\)-homology) of a smooth manifold to the compactly supported cohomology (resp. locally finite homology) of the manifold is an isomorphism. %
\end{abstract}

\maketitle

\tableofcontents
\section{Introduction}\label{sec:intro}

This article studies the stable homotopy theory of \emph{bornological algebras}, building on the work of Cuntz, Meyer and Rosenberg in \cite{CMR}. Roughly, a bornological vector space is the filtered cocompletion of the category of Banach spaces, and a bornological algebra is a bornological vector space equipped with a bounded multiplication. The category \(\BAlg\) of bornological algebras has all limits and colimits, and a well-behaved tensor product that commutes with filtered colimits. It contains several categories of topological algebras such as Fr\'echet algebras on the one hand, and \(\C\)-algebras without any topology on the other.  In \cite{CMR}, the authors construct a triangulated category \(\KKzero\) with properties analogous to Kasparov's \(KK\)-theory for separable \(C^*\)-algebras. As in the $C^*$-algebraic context, they prove that there is a natural functor \(\BAlg \to \KKzero\) which is initial with respect to functors into \emph{triangulated categories} satisfying smooth homotopy invariance, excision, and stability with respect to suitable operator ideals. But unlike Kasparov's \(KK\)-theory, the category \(\KKzero\) is not categorically rich enough for several applications. For instance, a powerful organisational principle in the Dirac-dual Dirac method used several proofs of the Baum-Connes conjecture with coefficients is the interpretation of the assembly map as the ``resolution" of a \(C^*\)-algebra by a proper \(C^*\)-algebra \cite{MN}, for which the conjecture is known \cite{chabert2000deux}. Unfortunately, this method cannot be applied to study related conjectures such as the Farrell-Jones conjecture and the Bost conjecture using bivariant \(K\)-theory, as it relies on the category having countable coproducts. A similar problem appears in the purely algebraic bivariant \(K\)-theory of Corti\~nas and Thom \cite{CortinasThom}. The lack of colimits also prevents a Universal Coefficient Theorem in \(\KK\) for many interesting classes of algebras, such as smooth functions on a manifold. Our first main theorem repairs these issues:

\begin{theorem}\label{thm:main:1}[Theorem \ref{prop:presentable-stable}]
    There is a canonical functor \(\BAlg \to \bKK\) into a presentable and stable \(\infty\)-category, that is initial for functors into cocomplete, stable \(\infty\)-categories satisfying smooth homotopy invariance, excision, stabilisation by Schatten class (or smooth compact) operators, and preserving countable filtered colimits of bornological algebras with injective structure maps. 
\end{theorem}

Our construction of \(\bKK\) is closely related to that of \(E\)-theory \cite{EisCASS} and localising motives \cite{blumberg2013universal}. We briefly sketch this as follows. Starting with \(\BAlg\), we first construct a stable \(\infty\)-categorical refinement \(\KK\) of the triangulated category \(\KKzero\) by taking the \(\infty\)-categorical Dwyer-Kan localisation of \(\BAlg\) at the \(\KKzero\)-equivalences. 
The homotopy category of \(\KK\) is indeed \(\KKzero\) by Theorem \ref{thm: ho KK is triangulated KKCMR}. We emphasize that our result, even in the 1-categorical context, is an improvement over the universal property of \(\KKzero\) as mentioned in \cite{CMR}, as the result in loc. cit. was just phrased with regards to triangulated categories. Having said that, the issues concerning colimits mentioned above do not really get resolved in the stable refinement. To remedy this, we take the \emph{Bousfield localisation} of the filtered cocompletion \(\Ind(\KK)\) at the set \(\cW_{\aleph_0}\) of maps of the form \(\colim_{n \in \N}yj(A_n) \to yj(\colim A_n)\), where \((A_n)\) is an inductive system of bornological algebras. The resulting category \(\bKK\) remains stable, but is now additionally presentable, and the map \[\bkk \colon \BAlg \overset{j}\to \KK \overset{y}\to \Ind(\KK) \to \bKK\] is essentially constructed to satisfy the universal property of Theorem \ref{thm:main:1}.  

The category \(\bKK\) has most of the features one desires from a bivariant \(K\)-theory, save one - a symmetric monoidal structure inherited from bornological algebras. At a technical level, this is because tensoring by an arbitrary bornological algebra does not preserve filtered colimits with injective structure maps. This is only true for nuclear bornological algebras. Actually, for our category to remain tractable in a sense that will be made precise shortly, we need to pass to a smaller category \(\NBA\) of \emph{strongly nuclear} bornological algebras. Roughly, a bornological vector space is strongly nuclear if it can be written as a filtered colimit of Banach spaces with structure maps that are \(p\)-Schatten class for every \(p>0\). These algebras appear in \cite{kelly2025localising} and also the setting of liquid vector spaces \cite{clausen2026condensed}. Notably, any nuclear Fr\'echet space is strongly nuclear. In \cite{kelly2025localising}, it is shown that the category of strongly nuclear algebras is a coreflective subcategory of bornological algebras, with right adjoint denoted by \((-)^\rig \colon \BAlg \to \NBA\). We then show the following:

\begin{theorem}\label{thm:main2}[Theorem \ref{thm:bousfield-nuc}]
The right Bousfield localisation \((-)^\rig :\BAlg \leftrightarrows \NBA: \iota\) descends to a right Bousfield localisation 
    \[
    (-)^\rig \colon \bKK \leftrightarrows \nKK: \iota.
    \]
\end{theorem}

The category \(\nKK\) is again presentable, stable, but now has a symmetric monoidal structure (Proposition \ref{prop:symmetric-monoidal}), and the canonical functor \[\nkk \colon \NBA \to \nKK\]
is constructed to be the universal symmetric monoidal functor into a presentable, stable symmetric monoidal category that is reduced, smoothly homotopy invariant, stable with respect to nuclear operators, and preserves filtered colimits of strongly nuclear bornological algebras. The notation \((-)^\rig\) and the construction of \(\nKK\) is inspired by the construction and proof of rigidity of noncommutative motives, which shows that the latter is generated by the image of nuclear categories (see \cites{nkp,efimov2025rigidity}). 

We now come to \(K\)-theory, the relevant version of which is \emph{diffeotopy} \(K\)-theory \(\K^\an\) due to Corti\~nas and Thom \cite{cortinas2008comparison}, adapted to the setting of bornological algebras. The functor \(\K^\an \colon \BAlg \to \Sp\) is smoothly homotopy invariant, excisive, and preserves filtered colimits with injective structure maps. For Fr\'echet algebras, the functor \(\K^\an\) agrees with the diffeotopy \(K\)-theory appearing in \cite{cortinas2008comparison}. Furthermore, by a version of Higson's automatic homotopy invariance theorem \cite{higson1988algebraic} for Fr\'echet operator ideals, the canonical map \[\KH(A \haotimes \mathcal{L}^p) \to \K^\an(A \haotimes \mathcal{L}^p)\] is an equivalence, where the left hand side is Weibel's homotopy algebraic \(K\)-theory. Let \(\cK\) denote the choice of stabilisation used to define \(\bKK\). We show that the functor \(\K^\an(\cK \otimes -)\) is corepresented by \(\C\) in our bivariant \(K\)-theory category:

\begin{theorem}\label{thm:main3}[Theorem \ref{thm:comparision-kd}]
    The canonical map \[\nKK(\C,A^\rig) \simeq \bKK(\C,A) \to \K^\an(\cK(A))\] is an equivalence of spectra, where \(\cK\) denotes the choice of stabilisation used to define \(\bKK\). When \(\cK = \sK \haotimes \mathcal{L}^p\) for \(p>1\), we have equivalences \[\nKK(\C,\C) \simeq \bKK(\C,\C) \simeq \KU\] of ring spectra.
\end{theorem}

Finally, we come to the Universal Coeffcient Theorem alluded to at the start. Let \(\bdd^{\nKK}\) denote the localising subcategory generated by \(\C\) inside \(\nKK\).

\begin{theorem}\label{thm:main4}[Proposition \ref{prop:bootstrap}]
    For \(A \in \bdd^{\nKK}\), we have an equivalence 
    \[
    \nKK(A,B) \simeq \mathbf{Mod}_{\K^\an(\sK)}(\K^\an(\sK(A)), \K^\an(\sK(B)))
    \] for all \(B \in \nKK\). 
\end{theorem}

We show that several interesting classes of algebras, such as smooth compactly supported functions on a manifold \(C_c^\infty(M)\) and Leavitt path algebras belong to the localising subcategory generated by \(\C\). An important drawback of the smooth stabilisation is that the spectrum \(\K^\an(\sK)\) is inexplicit, as observed in \cite{cuntz2004bivariant}. This stems from the fact that it is not clear whether the functor \(K_0(\sK \haotimes -)\) is smoothly homotopy invariant (without any further stabilisation), that is, we do not yet have a version of Higson's theorem for \(\sK\). On the other hand, the smooth stablisation \(\sK\) is required for the Bott periodicity of our theory. We circumvent this problem by working in the localising subcategory generated by \((\mathcal{L}^p)^\rig\) inside \(\nKK\). To justify this, consider the local cyclic homology functor \(\mathbf{HA} \colon \BAlg \to \bD^{\mathrm{loc}}(\C)\) with values in a certain derived category of chain complexes of bornological vector spaces, constructed in \cite{Meyer:HLHA}. The definition of \(\mathbf{HA}\) is somewhat technical, but conceptually, it is the functor obtained when Hochschild homology is forced to be topologically nilinvariant, which remarkably makes it a nondegenerate invariant for \(C^*\)-algebras (unlike periodic cyclic homology).  We show that the restriction of \(\mathbf{HA}\) to \(\NBA\) factors through \(\nKK\), so that we get a map \(\nKK(A,B) \to \mathbf{HL}(A,B)\) of spectra. We then have the following:

\begin{theorem}\label{thm:main5}[Theorem \ref{thm:UCT}]
    The canonical Chern character \[\mathsf{Ch}_{\C} \colon \nKK(\C,A) \otimes \C \to \mathbf{HA}(\C,A)\] is an equivalence for \(A\) in the localising subcategory generated by \((\mathcal{L}^p)^\rig\). Dually, the map \[\mathsf{Ch}_\C^{\vee} \colon \nKK(A,(\mathcal{L}^p)^\rig) \otimes \C \simeq \bKK(A, \mathcal{L}^p) \otimes \C \to \mathbf{HA}(A,\C)\] is an equivalence on the localising subcategory generated by \(\C\).
\end{theorem}

The cases covered under Theorem \ref{thm:UCT} are appropriately stabilised versions of \(C_c^\infty(M)\) and graph algebras. In particular, in Corollary \ref{cor:smooth-chern} we deduce that the complexified Chern character induces isomorphisms
\[
K_n^\an(C_c^\infty(M,\sK \haotimes \mathcal{L}^p)) \otimes \C \to \mathsf{HA}_n(C_c^\infty(M)) \cong \bigoplus_{j \in \Z}H_c^{n+2j}(M,\C), \quad n \in \Z 
\] where the right hand side denotes periodified compactly supported sheaf cohomology. Dually, we get isomorphisms 
\[
K_n^\an(M) = \pi_n(\nKK(C_c^\infty(M),\C)) \otimes \C \to \mathsf{HA}^n(C_c^\infty(M)) \cong \bigoplus_{j \in \Z}H_{n+2j}^{\mathrm{lf}}(M,\C), \quad n \in \Z
\] where the right hand side denotes periodified locally finite homology.

\textit{Acknowledgments: The authors would like to thank  Ulrich Bunke, Guillermo Corti\~nas, Johannes Ebert, Ralf Meyer, Maxime Ramzi, Melvin Weiss and Christoph Winges for several insightful discussions. 
Throughout the course of writing this paper, AD was supported by the Deutsche Forschungsgemeinschaft (DFG, German
Research Foundation) under Germany's Excellence Strategy – GZ 2047/1 -
390685813, Hausdorff Center for Mathematics, and by the Deutsche Forschungsgemeinschaft (DFG, German Research Foundation) - Research Grant 574320459, and GRK 2491 Project-ID: 398436923. DM was funded
by a DFG Eigenestelle (project number 534946574) and the Deutsche Forschungsgemeinschaft
(DFG, German Research Foundation) under Germany’s Excellence Strategy EXC 2044390685587,
Mathematics Münster: Dynamics Geometry Structure and a UK Research and Innovation Horizon
Europe Guarantee MSCA Postdoctoral Fellowship.}

\comment{\textcolor{red}{I suggest the following notation:}
\textcolor{blue}{
\(\BAlg\), \(\kkzero_p: \BAlg \to \KKzero_{\cK_p}\) for \(\cK_p = \sK \haotimes \cK_{\mathcal{L}^p}\), and the one-categorical CMR-theory, \(\bkk_{p} \colon \BAlg \to \KK_p\) for the \(\infty\)-categorical localisation relative to \(\cK_p\)-stablisation,. Similarly, for the smooth compact stablisations: \(\kkzero_{\s} \colon \BAlg \to \KKzero_{\s}\), \(\bkk_\s \colon \BAlg \to \bKK_\s\). For the nuclear thing, there's a unique functor \(\NBA \to \nKK_\s\). We know \[\KK_{p}(\C,\C) = \KU = \bKK_p(\C,\C)\], \[\KK_\s(\C,\C) = \K^\an(\sK) = \bKK_\s(\C,\C)\]. In general, \[\KK_{p}(\C,A) = \K^\an(\cK_p(A)) = \bKK_p(\C,A)\] and \[\nKK(\C,A) = \K^\an(\sK(A))\]. }}
\section{Background on bornological spaces and algebras}

In this section, we recall some basic terminology from bornological analysis following \cite{Meyer:HLHA} and \cite{HENRI}. A reader already familiar with bornologies may skip straight to Definition \ref{def:infinitely-nuclear}. 

\begin{definition}\label{def:bornologies}
 A \emph{bornology} on a set \(X\) is a collection \(\bdd_X\) of its subsets, called \emph{bounded sets}, such that singletons are bounded, finite unions of bounded subsets are bounded, and subsets of bounded sets are bounded. The pair \((X, \bdd_X)\) is called a \emph{bornological set}. A map \(f \colon X \to Y\) is said to be \emph{bounded} if it maps bounded subsets of \(X\) to bounded subsets of \(Y\). 
\end{definition}

Denote by $\mathfrak{B}\!\operatorname{Set}$ the category whose objects are bornological sets and morphisms are bounded maps. The category $\mathfrak{B}\!\operatorname{Set}$ is bicomplete, and the forgetful functor $\mathfrak{B}\!\operatorname{Set} \rightarrow \operatorname{Set}$ lifts all limits and colimits uniquely.

\subsection{Categories of bornological vector spaces}
We now define several categories, whose objects are vector spaces equipped with a bornology satisfying certain conditions. But before that, we need a preliminary definition.
\begin{definition}{(cf. \cite[Section 2.1]{CMR}, \cite[Section 1.1.2]{Meyer:HLHA})}
    Let $V$ be a vector space, and $B \subseteq V$ a subset. Consider the following properties on $B$.
    \begin{enumerate}
        \item \textbf{convex}: $tx+(1-t)y \in B$ for all $x,y \in B, t \in [0,1]$.
        \item \textbf{absolutely convex}: $\lambda.B \subseteq B$ for all $\lambda \in \C$ with $\lvert \lambda \rvert \leq 1$.
        \item \textbf{internally closed}: $B=\bigcap_{\varepsilon>0}(1+\varepsilon).B$.
        \item \textbf{absoring}: $\R_+.B=V$.
    \end{enumerate}
    A subset $B \subset V$ is called a \emph{disk} in $V$ if $B$ satisfies the first three properties above. It can then be shown that $B$ determines a unique seminorm on 
    \[
    V_B:=\R_+.B \subseteq V
    \] 
    such that $B$ is the closed unit ball in $V_B$ with respect to that semi-norm. Indeed, this semi-norm is the so-called \emph{gauge semi-norm} of $B$, given by\[
    \nu_B: V_B \rightarrow \R_{\geq 0}, ~~~~~\nu_B(v):=\inf\{r \in \R_{\geq 0} \mid v \in rB\}.
    \]
    Furthermore, $B$ becomes an absorbing disk in $V_B$. $B$ is said to be \emph{complete} if $V_B$ is a Banach space.
\end{definition}
\begin{definition}
    \hfill
    \begin{enumerate}
        \item Denote by $\mathfrak{B}\!\operatorname{Vect}$ the category of $\C$-vector space objects in $\mathfrak{B}\!\operatorname{Set}$, i.e., the objects are $\C$-vector spaces $V$ such that the addition $V \times V \rightarrow V$ and scalar multiplication $\C \times V \rightarrow V$ maps are bounded\footnote{the base field $\C$ is always equipped with the Euclidean bornology.}, and the morphisms are the bounded linear maps.
        \item Denote by $\Bor$ the full subcategory of $\mathfrak{B}\!\operatorname{Vect}$ spanned by \textit{convex} bornological vector spaces: where the bornology is stable under the formation of convex hulls.
        \item Denote by $\Bor^{\mathrm{sep}}$ the full subcategory of $\Bor$ spanned by \textit{separated} convex bornological vector spaces: where the only bounded vector subspace is $\{0\}$.
        \item Denote by $\CBorn$ the full subcategory of $\Bor^{\mathrm{sep}}$ spanned by \textit{complete} separated convex bornological vector spaces: where every bounded subset is contained in a bounded complete disk.
    \end{enumerate}
\end{definition}
\begin{proposition}\label{prop: the cats are bicomplete}
    Each of the categories described above are bicomplete.
\end{proposition}
\begin{proof}
    The case of $\mathfrak{B}\!\operatorname{Vect}$ is subsumed in the content of pages 32-35 of \cite{HENRI}. The assertion for the rest is the content of \cite[Proposition 1.126]{Meyer:HLHA}.
\end{proof}
In the remainder of this subsection, we want to understand certain structural properties of the categories $\Bor,\Bor^{\mathrm{sep}},\CBorn$ better. To this end, we recall the \emph{dissection} functor.

\begin{definition}\label{defi: diss}{(cf. \cite[Definition 1.138]{Meyer:HLHA})}
    Let $V \in \Bor$. Denote by $\mathfrak{S}_d(V)$ the directed set of bounded disks in $V$ with the partial order being the relation of \emph{absorption}, i.e., \(B \leq D\) if and only if \(B \subseteq cD\) for some positive real number $c$. Note that it is indeed a directed set requires that $V$ is convex, see \cite[Lemma 1.8]{Meyer:HLHA}. In what follows, $\mathsf{Norm}_{\frac{1}{2}},\mathsf{Norm},\mathsf{Ban}$ denote the categories of (small) semi-normed, normed and Banach spaces respectively (with the morphisms being bounded linear maps), and $\Ind(-)$ denotes taking the ind-completion\footnote{this is the notation used in \cite{Meyer:HLHA}. Later in this article, we also consider other ($\infty$-categorical) ind-completions, and denote that by $\mathbf{Ind}_{\aleph_0}(-)$.}.
    \begin{enumerate}
        \item The semi-normed spaces $V_B$ for $B \in \mathfrak{S}_d(V)$ form an inductive system in $\mathsf{Norm}_{\frac{1}{2}}$ because we have injective bounded linear maps $V_B \rightarrow V_D$ for $B, D \in \mathfrak{S}_d(V)$ with $B \leq D$. This construction is furthermore functorial. The resulting functor 
        \[
        \operatorname{diss}: \Bor \rightarrow \Ind(\mathsf{Norm}_{\frac{1}{2}})
        \]
        is called disssection functor.
        \item If $V \in \Bor^{\mathrm{sep}}$, then $V_B$ is a normed space for all $B \in \mathfrak{S}_d(V)$, so that dissection restricts to a functor 
        \[
        \operatorname{diss}: \Bor^{\mathrm{sep}} \rightarrow \Ind(\mathsf{Norm}).
        \]
        \item If $V \in \CBorn$, then we get an isomorphic object of $\Ind(\mathsf{Norm})$ if we replace $\mathfrak{S}_d(V)$ with its cofinal subset of complete bounded disks. This way, we can further restrict to a functor
        \[
        \operatorname{diss}: \CBorn \rightarrow \Ind(\mathsf{Ban})
        \]
    \end{enumerate}
\end{definition}
\begin{lemma}\label{lem: ind of}
    We have a sequence of adjunctions, where the left adjoints go from left to right, and the right adjoints are fully faithful
\begin{equation}\label{eqn: adjunction on ind semi norm etc}
    \begin{tikzcd}
	\Ind(\mathsf{Norm}_{\frac{1}{2}}) && \Ind(\mathsf{Norm}) && \Ind(\mathsf{Ban})
	\arrow["\Ind(\operatorname{Haus.~quotient})", shift left=3, curve={height=-18pt}, from=1-1, to=1-3]
	\arrow[hook', from=1-3, to=1-1]
	\arrow["\Ind(\operatorname{completion})", shift left=3, curve={height=-18pt}, from=1-3, to=1-5]
	\arrow[hook', from=1-5, to=1-3]
\end{tikzcd}.
\end{equation}
Furthermore, all the categories in the above diagram are bicomplete. 
\end{lemma}
\begin{proof}
    In light of \cite[Exercise 9.3.3.11]{kerodon}, the first part of the lemma follows from the familiar adjunctions
\begin{equation}\label{eqn: adjunction between semi norm norm ban}
    \begin{tikzcd}
	{\mathfrak{Norm_{\frac{1}{2}}}} && {\mathsf{Norm}} && {\mathsf{Ban}}
	\arrow["{\operatorname{Haus.~quotient}}", shift left=3, curve={height=-18pt}, from=1-1, to=1-3]
	\arrow[hook', from=1-3, to=1-1]
	\arrow["{\operatorname{completion}}", shift left=3, curve={height=-18pt}, from=1-3, to=1-5]
	\arrow[hook', from=1-5, to=1-3]
\end{tikzcd}
\end{equation}
where the first functor is taking the Hausdorff quotient.
The second part of the lemma follows from \cite[Proposition 1.135]{Meyer:HLHA}. 
\end{proof}

\begin{remark}\label{rem: enriched adjunction}
    Note that the categories $\mathsf{Norm}_{\frac{1}{2}}$, $\mathsf{Norm}$ and $\mathsf{Ban}$ are enriched over themselves. Moreover, they are closed symmetric monoidal categories, with the tensor product being the appropriate variants of the classical projective tensor product of Banach spaces. Since all the three categories have finite colimits, we formally get that their ind-completions are also closed symmetric monoidal.
    
    It follows directly from the definition of operator (semi-)norms that in  \Cref{eqn: adjunction between semi norm norm ban}, $\operatorname{Haus.~quotient}$ is a $\mathsf{Norm}_{\frac{1}{2}}$-enriched left-adjoint, and $\operatorname{completion}$ is a $\mathsf{Norm}$-enriched left-adjoint to the respective (enriched) inclusion functors. As a consequence, the adjunctions in \Cref{eqn: adjunction on ind semi norm etc} are also \emph{appropriately enriched}. We will frequently be using this phrase in what follows, and the meaning of it should be clear from the context. 
\end{remark}
 
\begin{proposition}\label{prop: diss and classical together adjunctions}\cite[Theorem 1.139]{Meyer:HLHA}
    Each of the dissection functors in \Cref{defi: diss} are fully faithful and have left-adjoints given by the inductive limit functor\footnote{Recall the usual notation of an adjunction, where the left-end of $\dashv$ points to the left-adjoint.}
\begin{equation}\label{eqn: the adjunction diagram}
\begin{tikzcd}
	{\Bor} &&& {\Bor^{\mathrm{sep}}} &&& {\CBorn} \\
	\\
	\\
	\Ind(\mathsf{Norm}_{\frac{1}{2}}) &&& \Ind(\mathsf{Norm}) &&& \Ind(\mathsf{Ban})
	\arrow["{\operatorname{diss}}"{pos=0.3}, hook, from=1-1, to=4-1]
	\arrow[hook', from=1-4, to=1-1]
	\arrow["{\operatorname{diss}}"{pos=0.3}, hook, from=1-4, to=4-4]
	\arrow[hook', from=1-7, to=1-4]
	\arrow["{\operatorname{diss}}"{pos=0.3}, hook, from=1-7, to=4-7]
	\arrow["\dashv"', shift left=4, from=4-1, to=1-1]
	\arrow["\Ind(\operatorname{Haus.~quotient})", shift left, curve={height=-18pt}, from=4-1, to=4-4]
	\arrow["\dashv"', shift left=4, from=4-4, to=1-4]
	\arrow[hook', from=4-4, to=4-1]
	\arrow["\Ind(\operatorname{completion})", shift left, curve={height=-18pt}, from=4-4, to=4-7]
	\arrow["\dashv"', shift left=4, from=4-7, to=1-7]
	\arrow[hook', from=4-7, to=4-4]
\end{tikzcd}
\end{equation}

\noindent in which every square that should commute does commute. All the adjunctions in sight are appropriately enriched.
\end{proposition}
\comment{\begin{proof}
    The assertions about the dissection functors are contained in \cite[Theorem 1.139]{Meyer:HLHA}. The commutativity assertion follows from a routine yoga with adjunctions. The assertion about the enrichment of the adjunctions follows from \Cref{rem: enriched adjunction} and \cite[Theorem 1.139.5]{Meyer:HLHA}.
\end{proof}}

\begin{proposition}\label{prop: adjunctions between categories of bvs}
    We have a sequence of appropriately enriched adjunctions, where the left adjoints go from left to right, and the right adjoints are fully faithful.
\begin{equation}\label{eqn: adjunctions of categories of bvects}
    \begin{tikzcd}
	{\mathfrak{B}\!\operatorname{Vect}} && {\Bor} && {\Bor^{\mathrm{sep}}} && {\CBorn}
	\arrow["{\operatorname{convexification}}", shift left=3, curve={height=-18pt}, from=1-1, to=1-3]
	\arrow[hook, from=1-3, to=1-1]
	\arrow["{\operatorname{separated~quotient}}", shift left=3, curve={height=-18pt}, from=1-3, to=1-5]
	\arrow[hook, from=1-5, to=1-3]
	\arrow["{\operatorname{completion}}", shift left=3, curve={height=-18pt}, from=1-5, to=1-7]
	\arrow[hook, from=1-7, to=1-5]
\end{tikzcd}
\end{equation}
\end{proposition}
\begin{proof}
    This follows directly from the results of this section, from \Cref{lem: ind of} to \Cref{prop: abstract symm mon loc}. Note that the left-adjoints labeled above are the relevant bottom-right-up composites in \Cref{eqn: the adjunction diagram}. For example, the completion functor is the composite
    \begin{equation}\label{eqn: completion itinerary}
        \Bor^{\mathrm{sep}} \xrightarrow{\operatorname{diss}} \Ind(\mathsf{Norm}) \xrightarrow{\Ind(\operatorname{completion})} \Ind(\mathsf{Ban}) \xrightarrow{\colim} \CBorn .
    \end{equation}
\end{proof}
\begin{remark}
    The following remarks are in order.
    \begin{enumerate}

        \item The enrichment of $\mathfrak{B}\!\operatorname{Vect}$ over itself is given by the \emph{uniform bornology}: given $V_1,V_2 \in \mathfrak{B}\!\operatorname{Vect}$, we equip $\operatorname{Hom}_{\mathfrak{B}\!\operatorname{Vect}}(V_1,V_2)$ with the bornology where equibounded subsets are declared to be bounded. See \cite[Proposition 1.111]{Meyer:HLHA} for more details. In light of \Cref{prop: abstract symm mon loc}, this also yields the obvious self-enrichments of each of the categories mentioned in \Cref{eqn: adjunctions of categories of bvects}.
        \item The symmetric monoidal structure on $\mathfrak{B}\!\operatorname{Vect}$ can be given explicitly: for $(V_1,B_1), (V_2,B_2) \in \mathfrak{B}\!\operatorname{Vect}$, denote by $V_1 \otimes_{\mathfrak{B}} V_2$ to be the vector space tensor product $V_1 \otimes V_2$ equipped with the vector bornology generated by the subsets
        \[
        B_1 \otimes B_2:= \{b_1 \otimes b_2 \mid b_1 \in B_1, b_2 \in B_2 \}.
        \]
    See \cite[Section 1.3.6]{Meyer:HLHA} for details. Again in light of \Cref{prop: abstract symm mon loc}, this yields an obvious symmetric monoidal structures on each of the categories mentioned in loc.\! cit.\! The symmetric monoidal structure one obtains this way on $\CBorn$ is the \emph{completed projective bornological tensor product}, denoted by $-\haotimes -$ henceforth. See \cite[Definition 1.83]{Meyer:HLHA} for more details and properties. 
    \end{enumerate}
\end{remark}
We now give examples of (complete) bornological vector spaces. See \cite[Section 1.1.4]{Meyer:HLHA} for more examples and elaborations. We mention at this juncture that there is an obvious notion of \emph{a convex (vector) bornology generated by a collection of subsets} in any vector space $V$, and we will be using it in some of the examples below.
\begin{example}\label{ex: BVect}
\begin{enumerate}
\hfill
    \item  For any vector space \(V\), there is a minimal choice of a vector bornology on \(V\): a subset $B \subseteq V$ is bounded if and only if there is a finite dimensional subspace $W \subseteq V$ such that $B \subseteq W$ and $B$ is bounded in $W$ in the usual sense. This is called the \emph{fine bornology}. It is always complete, and is the canonical choice of bornology on a finite dimensional vector space. Denote by \(\mathsf{Fin}(V)\) the vector space \(V\) with the fine bornology. 
    \item Let \(V\) be a locally convex topological vector space. We call a subset \(S\)  of $V$ \emph{von Neumann bounded} if \(p \colon V \to \R_{\geq 0}\) is bounded in the usual sense for every continuous seminorm \(p\) on $V$. The vector bornology on $V$ generated by the von-Neumann bounded subsets is called the \emph{von Neumann bornology} on $V$, and denoted by $\vN(V)$. Note that $\vN(V)$ is separated if $V$ is Hausdorff, and complete if $V$ is complete (see \cite{HENRI}). This yields a lax monoidal functor
    \[
    \operatorname{vN:(CpltLCTVS},\hat\otimes_{\pi}) \rightarrow (\CBorn,\haotimes).
    \]
    Restricted to the full subcategory of LF-spaces, the functor $\vN(-)$ is fully faithful (see \cite[Theorem 1.33]{Meyer:HLHA}). Furthermore, if $V=\bigcup_{n \in \N}V_n$ is an LF-space which is the strict inductive limit of the Fréchet spaces $V_n$, then $\vN(V) \cong \operatorname{colim}_{n \in \N}\vN(V_n)$, see \cite[Theorem 1.32]{Meyer:HLHA}. Furthermore, $\vN(-)$ restricts to a strongly monoidal functor on the full subcategory of \emph{nuclear} Fréchet spaces, see \cite[Proposition 4.2.29]{PhDBorn}.
    \item Let \(V\) be a locally convex topological vector space. We call a subset \(S\) of $V$ \emph{precompact} if for every neighbourhood \(U\) of the origin, there is a finite subset \(F\) such that \(S \subseteq F + U\). Equivalently, $S$ is precompact if and only if the closure of $S$ in the (topological) completion of $V$ is compact. The vector bornology generated by the precompact subsets is called the \emph{precompact bornology}. Note that $\Cpt(V)$ is separated if $V$ is Hausdorff, and complete if $V$ is complete. This yields a lax monoidal functor 
    \[
    \operatorname{Cpt:(CpltLCTVS,\hat\otimes_\pi)} \rightarrow (\CBorn,\haotimes).
    \]
    \emph{Unless otherwise mentioned, we equip a locally convex topological vector space $V$ with the precompact bornology.} As with the von Neumann bornology, the functor $\Cpt(-)$ is fully faithful (see \cite[Theorem 1.33]{Meyer:HLHA}) on the full subcategory of LF-spaces, and if $V=\bigcup_{n \in \N}V_n$ is an LF-space, then $\Cpt(V) \cong \operatorname{colim}_{n \in \N}\Cpt(V_n)$, see \cite[Theorem 1.32]{Meyer:HLHA}. Finally, the functor $\Cpt(-)$ restricts to a strongly monoidal functor on the full subcategory of Fréchet spaces, see \cite[Theorem 1.87]{Meyer:HLHA}.
\end{enumerate}
\end{example}

\subsubsection{Strongly nuclear bornological vector spaces}

A class of bornological vector spaces that will play an important role later on is that of strongly nuclear bornological vector spaces introduced in \cite{kelly2025localising}. 

\begin{definition}\label{def:infinitely-nuclear}
\begin{enumerate}
    \item A map \(f: V \to W\) between separable Banach spaces is said to be \emph{trace-class of order \(p\)} if it may be written as \begin{equation}\label{eqref:standard-rep}
f(x) = \sum_{n = 0}^\infty \lambda_n f_n(x) g_n
\end{equation} for bounded sequences \((f_n) \in V'\) and \((g_n) \in W\), and \(\lambda_n \in l^p(\N)\);
\item We say that \(f\) has order \(0\) if it is of \emph{order \(p\)} for every \(p \in (0,1]\).  By \cite[Chapitre 2, §1.2 Corollaire pp9]{grothendieck1955produits}, any such map has a standard representation as in Equation \eqref{eqref:standard-rep} with \((\lambda_n)\) a \emph{rapidly decreasing sequence}, that is, for each \(k\), the sequence \((n^k \lambda_n)\) is bounded;
\item We call a complete bornological vector space \emph{strongly nuclear} \(V\) if for any bounded subset \(S\), there is a complete bounded disc \(T\) such that the induced map \(V_S \to V_T\) is trace-class of order \(0\).
\end{enumerate}
\end{definition}

Denote by \(\mathsf{Nuc}^\infty \subseteq \CBorn\) the full subcategory of strongly nuclear complete bornological vector space. By \cite[Theorem 6.30]{kelly2025localising}, the inclusion has a right adjoint \((-)^\rig \colon \CBorn \to \mathsf{Nuc}^\infty\) that equips a complete bornological vector space \(V\) with the bornology where a subset \(S\) is bounded if there is an increasing sequence \((T_n)\) of complete bounded disks containing \(S\), such that the induced inclusions \(V_{T_i} \to V_{T_{i+1}}\) are trace-class of order zero for each \(i\). We accumulate the following facts about strongly nuclear spaces that will be needed later:

\begin{proposition}\label{prop:infinitely-nuc}
The following hold:
    \begin{enumerate}
        \item The category of strongly nuclear complete bornological vector spaces is closed under countable colimits, closed subspaces, tensor products;
        \item If \(V\) is a nuclear Fr\'echet space with the von Neumann (or equivalently, the precompact bornology), then \(V\) is strongly nuclear;
        \item The rigidification functor preserves countable filtered colimits, i.e., if \(V \cong \colim_n V_n\), then \(V^\rig \cong \colim_n V_n^\rig\);
        \item If \(V\) is strongly nuclear and \(W\) is an arbitrary complete bornological vector space, then \((V \haotimes W)^\rig \cong V \haotimes W^\rig\).
        \end{enumerate}
\end{proposition}

\begin{proof}
    The first claim is \cite[Proposition 6.32]{kelly2025localising}. For the second claim, let \(B\) be a bounded subset in a nuclear Fr\'echet space \(V\). Then there is a complete, bounded disk \(D\) such that \(B \subseteq D\) and \(V_B \to V_D\) is trace-class of order \(1\). Iteratively, we may find a sequence of bounded disks \(B \subseteq D_0 \subseteq D_1 \subseteq \cdots \), such that \(V_{D_i} \to V_{D_{i+1}}\) is trace-class of order \(1\) for each \(i\). By \cite[Theorem 1.154]{Meyer:HLHA}, we may find a bounded disk \(T\) that absorbs each \(D_n\). Therefore, for any \(n\), the inclusion \(B \subseteq T\) admits a factorisation \(V_{B} \to V_{D_0} \to \cdots V_{D_n} \to V_T\), where each intermediate map is trace-class of order \(1\). By \cite[Chapitre 2, §1.3 Corollaire pp
10]{grothendieck1955produits}, the map \(V_B \to V_T\) then has order at most \(\frac{2}{n-1}\). As the \(n\) was arbitrary, the map \(V_B \to V_T\) is trace-class of order zero. 

     For the third claim, in one direction, since \((-)^\rig\) is a right adjoint functor, it preserves kernels, so that we have inclusions \(V_i^\rig \to V_j^\rig\). Furthermore, \(V_j^\rig \subseteq V_j \to V\), inducing a canonical map \(\colim_n V_i^\rig \to V^\rig\). In the other direction, let \(S \subseteq V^\rig\) be a bounded subset. Then there is a sequence \((T_n)\) of bounded discs containing \(S\) such that the maps \(V_{T_n} \to V_{T_{n+1}}\) are trace-class of order \(0\). Since each \(T_n\) is bounded in \(V\), by the definiton of the filtered colimit bornology, there is an \(m(n)\) and a bounded disc \(S_{m(n)}\) such that \(T_n\) is contained in the image of \(S_{m(n)}\) in \(V\). As a consequence, the inclusions \(V_{T_n} \subseteq V^\rig \subseteq V\) factor through \(V_{T_n} \subseteq V_{S_{m(n)}}\) for each \(n\). We then have an inclusions \(V_{T_{n-1}} \subseteq V_{T_n} \subseteq V_{S_{m(n)}}\) which are trace-class of order \(0\), so that we have factorisations \(V_{T_{n-1}} \to \mathcal{S}(\N) \to V_{S_{m(n)}}\) through the Schwartz space for each \(n\). Since \(\s(\N)\) is strongly nuclear by part (ii), we have a factorisation \(\s(\N) \subseteq V_{S_{m(n)}}^\rig \subseteq V_{S_{m(n)}}\). Since \(S\) is contained in the image of \(V_{T_{n-1}}\), it is therefore contained in the image of \(V_{S_{m(n)}}^\rig\) inside \(V^\rig\) as required. 
    
    For the last claim, the counit of the adjunction yields a bounded linear map \(W^\rig \to W\), and therefore a bounded linear map \(V \haotimes W^\rig \to V \haotimes W\). Now since by the first part, the tensor product of strongly nuclear spaces is strongly nuclear, \(V \haotimes W^\rig\) is strongly nuclear, so that by the universal property of the functor \((-)^\rig\), we get a bounded linear map \(V \haotimes W^\rig \to (V \haotimes W)^\rig\). In the other direction, let \(X\) be strongly nuclear, and \(f: X \to V \haotimes W\) a bounded linear map. Since \(X\) and \(V\) are strongly nuclear, we may write them as \(X \cong \colim_n X_n\) and \(V \cong \colim_n V_n\) with transition maps that are trace-class of order zero. Since the map \(f\) is bounded, for each \(n\), there are \(m(n)\) and \(k(n)\), such that \(f\) restricts to a bounded map \(X_n \to V_{m(n)} \haotimes W_{k(n)}\). Since the maps \(X_n \to X_{n'}\) are trace-class of order zero for each \(n' \geq n\), and since the composition of a trace-class map (of any order) with an arbitrary bounded map between Banach spaces is trace-class of the same order, by increasing the indices \(m(n)\) and \(k(n)\) if necessary,  we may arrange that the maps \(X_n \to V_{m(n)} \haotimes W_{k(n)}\) are themselves of order zero. We now claim that \(X_n \to V_{m(n)} \haotimes W_{k(n)}\) factors through \(X_n \to V_{m(n)} \haotimes W_{k(n)}^\rig\), yielding maps \(X_n \to V \haotimes W^\rig\) for each \(n\), and therefore a unique bounded linear map \(X \to V \haotimes W^\rig\). Since \(V\) is in particular nuclear, it is flat, so that tensoring the inclusion \(W^\rig \to W\) by \(V\) remains an inclusion \(V \haotimes W^\rig \to V \haotimes W\), factoring the original map \(X \to V \haotimes W\). To prove the missing claim, we need to show that if \(h:V \to V'\), \(k: X \to V \haotimes W\) are order zero maps of Banach spaces, then the induced map \(X \overset{k}\to V \haotimes W \overset{h \otimes 1_W}\to V' \haotimes W\) factors through \(X \to V' \haotimes W^\rig\). For this, we write the order zero maps \(h\) and \(k\) in standard form \[h(v) = \sum_{m \in \N}\alpha_m \psi_m(v) v_m, \quad k(x) = \sum_{n \in \N} \beta_n \varphi_n(x) u_n,\] where \((\alpha_m)\) and \((\beta_n)\) are rapidly decreasing sequences, \((\psi_m)\), \((\varphi_n)\) are uniformly bounded sequences in the bornological duals \(\Hom(V,\C)\) and \(\Hom(X,\C)\), respectively, and \((u_n) \in V \haotimes W\) and \((v_m) \in V'\) are bounded sequences.  Consider the set \(\xi_{m,n} \defeq (\psi_m \haotimes 1_W)(u_n)\in W\); this is bounded as \((\psi_m)\) and \((u_n)\) are both bounded. Then \[(h \otimes 1_W)(k(x)) = \sum_{m,n} \alpha_m \beta_n \varphi_n(x)v_m \otimes \xi_{m,n}.\] 
    
    Consider the map \(L \colon \s(\N^2) \to W\), \(z \mapsto \sum_{m,n}z_{m,n}\xi_{m,n}\). By part (ii), since \(\s(\N^2)\) is strongly nuclear, the map \(L\) factors through \(W^\rig\).  Now let \(\gamma_{m,n} \defeq |\alpha_m\beta_n|^{1/2}\). Since \((\alpha_m)\) and \((\beta_n)\) are rapidly decreasing, the subset \(\{\gamma_{mn}e_{mn}:m,n \in \N\} \subseteq \s(\N^2)\) is bounded in the von Neumann bornology, where \(e_{mn}\) are the standard basis vectors in \(\s(\N^2)\).  As a consequence, there is a bounded disk \(D \subseteq W^\rig\) containing the bounded set \(a_{m,n} = L(\gamma_{m,n})\). Setting \(\theta_{m,n} = \alpha_m \beta_n/|\alpha_m \beta_n|\), we get a bounded linear map 
    \[
    X \to V' \haotimes W_D, \quad x \mapsto \sum_{m,n} \gamma_{m,n}\theta_{m,n} \psi_n(x) v_m \otimes a_{m,n}
    \]  that factors \((h \otimes 1_W)\circ k\).
\end{proof}

\subsection{Categories of bornological algebras}
From discussions in the previous subsection, we conclude that each functor in the below diagram is (at least) lax monoidal.
\begin{equation}\label{eqn: adjunctions omnibus vs}
\begin{tikzcd}
	{\Bor} &&& {\Bor^{\mathrm{sep}}} &&& {\CBorn} \\
	\\
	\\
	\Ind(\mathsf{Norm}_{\frac{1}{2}}) &&& \Ind(\mathsf{Norm}) &&& \Ind(\mathsf{Ban})
	\arrow["{\operatorname{separated ~quotient}}", shift left=3, curve={height=-18pt}, from=1-1, to=1-4]
	\arrow["{\operatorname{diss}}"{pos=0.4}, hook, from=1-1, to=4-1]
	\arrow[hook', from=1-4, to=1-1]
	\arrow["{\operatorname{completion}}", shift left=3, curve={height=-18pt}, from=1-4, to=1-7]
	\arrow["{\operatorname{diss}}"{pos=0.4}, hook, from=1-4, to=4-4]
	\arrow[hook', from=1-7, to=1-4]
	\arrow["{\operatorname{diss}}"{pos=0.4}, hook, from=1-7, to=4-7]
	\arrow["\dashv"', shift left=4, from=4-1, to=1-1]
	\arrow["\Ind(\operatorname{Haus.~quotient})", shift left, curve={height=-18pt}, from=4-1, to=4-4]
	\arrow["\dashv"', shift left=4, from=4-4, to=1-4]
	\arrow[hook', from=4-4, to=4-1]
	\arrow["\Ind(\operatorname{completion})", shift left, curve={height=-18pt}, from=4-4, to=4-7]
	\arrow["\dashv"', shift left=4, from=4-7, to=1-7]
	\arrow[hook', from=4-7, to=4-4]
\end{tikzcd}\end{equation}
Consequently, it induces a diagram of categories of algebras over each of the closed symmetric monoidal categories above. Since the monoidal structures were symmetric, the category of algebras also naturally inherit a symmetric monoidal structure; these monoidal structures are not closed, of course.  

An adjunction between monoidal categories where the left adjoint is strongly monoidal and the right adjoint is lax monoidal induces an adjunction on the category of monoids. However, the adjunctions are no longer enriched over the category of (bornological) vector spaces. Summarising, we get the following diagram of categories of algebras, where all the right adjoints are fully faithful. We omit labeling the functors to avoid notational clutter.
\[\begin{tikzcd}
	{\mathbf{Alg}(\Bor)} & {} & {\mathbf{Alg}(\Bor^{\mathrm{sep}})} && {\mathbf{Alg}(\CBorn)} \\
	\\
	\\
	{\mathbf{Alg}(\Ind(\mathsf{Norm}_{\frac{1}{2}}})) && {\mathbf{Alg}(\Ind(\mathsf{Norm}))} && {\mathbf{Alg}(\Ind(\mathsf{Ban}))}
	\arrow[shift left=3, curve={height=-18pt}, from=1-1, to=1-3]
	\arrow[hook, from=1-1, to=4-1]
	\arrow[hook', from=1-3, to=1-1]
	\arrow[shift left=3, curve={height=-18pt}, from=1-3, to=1-5]
	\arrow[hook, from=1-3, to=4-3]
	\arrow[hook', from=1-5, to=1-3]
	\arrow[hook, from=1-5, to=4-5]
	\arrow["\dashv"', shift left=4, from=4-1, to=1-1]
	\arrow[shift left, curve={height=-18pt}, from=4-1, to=4-3]
	\arrow["\dashv"', shift left=4, from=4-3, to=1-3]
	\arrow[hook', from=4-3, to=4-1]
	\arrow[shift left, curve={height=-18pt}, from=4-3, to=4-5]
	\arrow["\dashv"', shift left=4, from=4-5, to=1-5]
	\arrow[hook', from=4-5, to=4-3]
\end{tikzcd}\]


We now want to discuss limits and colimits in the categories of algebras. This is contained in the next result.
\begin{proposition}\label{prop: alg prop}
    Let $\mathcal{C}$ denote any of the categories in \Cref{eqn: adjunctions omnibus vs}. Then, the category $\mathbf{Alg}(\mathcal{C})$ is bicomplete. Furthermore, the forgetful functor $U: \mathbf{Alg}(\mathcal{C}) \rightarrow \mathcal{C}$ creates all limits and filtered colimits, and admits a left adjoint. 
\end{proposition}
\begin{proof}
    Note that the \emph{free monoid monad} $T : \mathcal{C} \rightarrow \mathcal{C}$ exists, and is given by 
    \[
    T(A) = \coprod_{n \in \N} A^{\otimes n},
    \]
    see \cite[Example 5.1.4.ii, Exercise 5.1.i]{CATINCON}. Furthermore, note that $T$ preserves filtered colimits since the symmetric monoidal structure on $\mathcal{C}$ is closed, as can be seen by performing a standard yoga on tensor products and colimits, and finally using that the diagonal in the product of a directed set with itself is cofinal\footnote{it is crucial that the colimit is filtered for this argument.}. 
    
    Since the category of algebras over the free monoid monad in a symmetric monoidal category is equivalent to the category of algebras in the symmetric monoidal category (see \cite[Example 5.2.6.iii]{CATINCON} for a proof when the base category is $\operatorname{Set}$, but the proof readily generalizes to any symmetric monoidal category), the bicompleteness of $\mathbf{Alg}(\mathcal{C})$ follows from \cite[Theorem 5.6.12]{CATINCON}. The assertion about limits and filtered colimits being created follows from \cite[Theorem 5.6.5]{CATINCON}. The assertion about the left adjoint of $U$ is \cite[Theorem VI.3.2]{MACLANE}.
\end{proof}
\begin{definition}
    We denote by $\BAlg$ the symmetric monoidal category $\mathbf{Alg}(\CBorn)$, and call its objects \textit{complete bornological algebras} over $\C$. Unless otherwise mentioned, we henceforth work with complete bornological algebras, and therefore simply refer to such algebras as bornological algebras. 
\end{definition}

\begin{remark}
    This recovers the definition of bornological algebra as given in \cite[Definition 2.5]{CMR}. See also \cite[Example 1.109]{Meyer:HLHA}.
\end{remark}
We end this section by a brief discussion on the category of \emph{strongly nuclear} bornological algebras, denoted $\NBA$. Recall the full subcategory \(\mathsf{Nuc}^\infty \subseteq \CBorn\) generated by strongly nuclear spaces in the sense of Definition \ref{def:infinitely-nuclear}. Since the tensor product of strongly nuclear bornological vector spaces is strongly nuclear, the inclusion \(\mathsf{Nuc}^\infty \subset \mathsf{CBorn}\) is strongly monoidal, so that its right adjoint is lax monoidal. We then get an adjunction 
\begin{equation}\label{eqn: NBA adj}
    \NBA \defeq \Alg(\mathsf{Nuc}^\infty) \leftrightarrows \Alg(\mathsf{CBorn}) = \BAlg.
\end{equation}
Note that the category $\NBA$ is bicomplete, as follows from the bicompleteness of $\BAlg$ and \Cref{lem: bunke seven}.
\section{Smooth homotopies, compact operators and extensions}
In this section, we recall the definitions and basic properties of the terms mentioned in the title. These will turn out to be the \emph{defining} properties of $\KKzero$, as we will see in Section 4.
\subsection{Smooth homotopies}

The notion of homotopy that works most conveniently in our context is \emph{smooth homotopies} or \emph{diffeotopies}, that we now recall. Consider the Fr\`echet algebra \(C^\infty([0,1])\), which we view as a complete bornological algebra with respect to the precompact bornology.

\begin{definition}\label{def:smooth-htpy}
    We say that two bounded algebra homomorphisms \(f_0,f_1 \colon  A \rightrightarrows B\) are \emph{smoothly homotopic} (or \emph{diffeotopic}) if there is a bounded algebra homomorphism \[
    h: A \rightarrow C^{\infty}([0,1],B)\]
    such that $ev_i \circ h=f_i$, with \(i = 0,1\).
\end{definition}
Although it is not immediately clear from the definition, it is shown in \cite[Page 92]{CMR} that this is an equivalence relation. The key to this demonstration is that the notion of smooth homotopy can be equivalently formulated using the subalgebra of $C^{\infty}([0,1],B)$ comprising of those functions whose $n^{\operatorname{th}}$ derivatives at $0$ and $1$ vanish for all $n \geq 1$.
\begin{definition}\label{def:htpy-functor}
    A functor $
    F: \BAlg \rightarrow \mathcal{C}
    $
    into an \(\infty\)-category \(\mathcal{C}\) is said to be \emph{smoothly homotopy invariant} if it maps smoothly homotopic maps to equivalences. Denote by $\mathbf{Ho}(\BAlg)$ the category whose objects are complete bornological algebras, and whose morphisms are smooth homotopy classes of maps, denoted $[-,-]$.\footnote{note that \cite{CMR} denotes it by $\langle-,- \rangle$.} Since we only work with smooth homotopies throughout this article, we often shorten smooth homotopies and smooth homotopy invariance to homotopy invariance.
\end{definition}
We now record a result which will be relevant for later use. 
\begin{lemma}\label{lem: HoBAlg}
    The canonical functor 
    \[
    \BAlg \rightarrow \mathbf{Ho}(\BAlg)
    \]
    presents $\mathbf{Ho}(\BAlg)$ as the localisation of $\BAlg$ at the class of smooth homotopy equivalences in the sense of $1$-categories.
\end{lemma}
\begin{proof}
    This is true in the much more general context of a category with a path object, see \cite[Section 1.3]{PathObject}.
\end{proof}

\subsection{Stability and stabilisation}
We now discuss a notion of stabilisation in the bornological land. Here, Cuntz-Meyer-Rosenberg \cite[Definition 7.19]{CMR} offer multiple choices of stabilisations for which the resulting $KK$-category is triangulated. 

We begin by defining the algebra of \emph{smooth compact operators} on a bornological algebra.
\begin{definition}
    \hfill
    \begin{enumerate}
        \item Let $V$ be a (complete) bornological vector space. A function $\N^2 \rightarrow V$ is \emph{rapidly decreasing} if there are a rapidly decreasing sequence of scalars $(\varepsilon_{ij})$ and a bounded subset $ T \subseteq V$ such that \[f(i,j) \in \varepsilon_{ij}T ~\forall~i,j \in \N^2.\] A set of such functions is \emph{uniformly rapidly decreasing} if the same $\varepsilon_{ij}$ and $T$ work for all its elements.
        \item Set $\s(\N^2,V)$ to be the bornological vector space of rapidly decreasing functions $\N^2 \rightarrow V$ with the bornology of uniformly rapidly decreasing subsets. Note that $\s(\N^2) \haotimes V \cong \mathcal{S}(\N^2,V)$ since \(\mathcal{S}(\N)\) is a nuclear bornological vector space.
        \item For a bornological algebra $A$, we set $\mathcal{K}_{\mathcal{S}}(A):=\s(\N^2) \haotimes A$ as a bornological vector space. The \emph{smooth stabilisation} of $A$ is $\mathcal{K}_\s(A)$ equipped with the multiplication extended from the dense subspace $M_\infty(A)$, where the multiplication is given by the usual formula
        \[
        f_1 *f_2(i,j)=\sum_{k \in \N}f_1(i,k)f_2(k,j).
        \]
    \end{enumerate}
\end{definition}

Next we consider the Schatten class stabilisation. For \(1\leq p <\infty\), let \(\mathcal{L}^p(l^2(\N))\) denote the ideal of Schatten class operators inside \(\mathcal{B}(l^2(\N))\). 

\begin{definition}
    The \emph{stabilisation of $A$ by the Schatten ideal $\mathcal{L}^p$} is defined by 
    \[
    \mathcal{L}^p(A):=\mathcal{L}^p(l^2) \haotimes A.
    \]
\end{definition}

For our applications we will be working with the combination of the $\sK$ and $\mathcal{L}^p$ stabilisations. At this stage, such a stabilisation may come off as unmotivated; for this we direct the reader to Remark \ref{rem: stability explanation}.
\begin{definition}
    For \(p\geq 1\), denote by $\cK_p(A)$ the stabilisation
    \[
    \cK_p(A):= \mathcal{K}_\s \haotimes \mathcal{L}^p \haotimes A .
    \]
\end{definition}
\begin{remark}
    For the purposes of the $KK$-theory of bornological algebras (which we describe in the next section), the stabilisation with the trace class operators as above yields the same theory as with any Schatten ideal $\mathcal{L}^p$ with $1 \leq p < \infty$; see \cite[Corollary 7.4]{CMR}.
\end{remark}

For each of the stabilisations (jointly denoted by \(\cK^{?}\)) considered above, there is a canonical bounded algebra homomorphism \[\iota_{A,\cK^?} \colon A \to \cK^?(A), \quad a \mapsto e_{11} \otimes a,\] where \(e_{11}\) is a rank \(1\)-operator.

\begin{definition}\label{def:stabilisation}
    We call a functor \(F \colon \BAlg \to \mathcal{C}\) \emph{stable with respect to \(\cK^?\)} it maps \(\iota_{A,\cK^?}\) to an equivalence for each \(A\). 
\end{definition}

\subsection{Exactness}

The class of extensions \[0 \rightarrow A \to B \to C \rightarrow 0\] 
of complete bornological algebras for which we require excision in bivariant \(K\)-theory are those where the map \(B \to C\) has a \emph{bounded} \(\C\)-linear section (this is a property, not part of the data). We shall refer to such extensions as being \emph{semi-split}. 

We next mention a couple of semi-split extensions of bornological algebras which are crucial for us. Another example is contained in \Cref{lem: toep}.
\begin{example}\label{ex: ex of ext}
    \hfill
    \begin{enumerate}
        \item \emph{(The cone extension, see \cite[Section 6.1]{CMR})} For any bornological algebra $A$, we have a semi-split extension 
        \[
        0 \rightarrow C_{\operatorname{flat}}^{\infty}((0,1),A) \rightarrow C_{\operatorname{flat}}^{\infty}((0,1],A) \xrightarrow{\operatorname{ev}_1} A  \rightarrow 0,
        \]
        where the suffix ``flat" denotes all positive derivatives of the functions vanish at the end-points. We call the bornological algebras $C_{\operatorname{flat}}^{\infty}((0,1),A) $  and $C_{\operatorname{flat}}^{\infty}((0,1],A)$ the \emph{suspension} and \emph{cone} of $A$ respectively, and denote them by $\Omega A$ and $CA$. By iterating the constructions, we obtain algebras $\Omega^nA$ and $C^nA$. Furthermore, $C^nA$ is contractible for all $n \in \N$. Note that the notation $\Omega$ will be justified shortly, when we show in \ref{lem: loop descends} that it is indeed the loop functor on the stable \(\infty\)-categorical refinement of the triangulated category contructed in \cite{CMR}.

        More generally, define the mapping cone of a map $f: A \rightarrow B$ in $\BAlg$ as
        \[
        C(f):=\{(a,b) \in A \oplus CB \mid f(a)=b(1)\}.
        \]
        See \cite[Section 6.4]{CMR} for further properties of this construction, including the definition of \emph{mapping cone triangle}, which we will be using in \Cref{lem: maps equivalent to fibrations}.3.
        \item (\emph{The tensor algebra extension, see \cite[Section 6.2]{CMR}}) The forgetful functor $\BAlg \rightarrow \CBorn$ has been shown to have a left-adjoint in \Cref{prop: alg prop}. We give a concrete of this left-adjoint for the convenience of the reader.

        Let $V \in \CBorn$. Let $V^{\haotimes n}$ be the completed projective tensor product of $n$ copies of $  V$. Let 
        \begin{equation}\label{eqn: TALG description}
            TV:= \bigoplus_{n \in \N} V^{\haotimes n}
        \end{equation}
        be equipped with the direct sum bornology. We define a multiplication on $TV$ by concatenation of tensors
        \[
        (v_1 \otimes ... \otimes v_n)\cdot(v_{n+1}\otimes... \otimes v_{n+m}) :=v_1 \otimes ... \otimes v_n \otimes v_{n+1}\otimes... \otimes v_{n+m}.
        \]
        This defines a bounded bilinear map $TV \times TV \rightarrow TV$ which is clearly associative. $TV$ thus becomes a complete bornological algebra with respect to this multiplication, and is called \emph{the tensor algebra of $V$}. Note that there is a canonical bounded \emph{linear} map $V \rightarrow TV$ identifying $V$ with the first direct summand in \Cref{eqn: TALG description}. 
        
        For each $A \in \BAlg$, we have a universal extension 
        \[
        0 \rightarrow JA \rightarrow TA \rightarrow A \rightarrow 0,
        \]
        where $JA:= \ker(TA \rightarrow A)$. The universality is in the sense that given a semi-split exact sequence
    \[
    \mathfrak{S}: 0 \rightarrow A \rightarrow B \rightarrow C \rightarrow 0
    \]
    and a morphism $f:D \rightarrow C$ in $\BAlg$,
    the \emph{classifying map} of $\mathfrak{S}$ with respect to $f$ is the unique (up to homotopy, see \cite[Lemma 6.17]{CMR}) homomorphism $\gamma:JD \rightarrow A$ rendering the obvious diagram
\[\begin{tikzcd}
	0 & JD & TD & D & 0 \\
	0 & A & B & C & 0
	\arrow[from=1-1, to=1-2]
	\arrow[from=1-2, to=1-3]
	\arrow["\gamma", from=1-2, to=2-2]
	\arrow[from=1-3, to=1-4]
	\arrow[from=1-3, to=2-3]
	\arrow[from=1-4, to=1-5]
	\arrow["f", from=1-4, to=2-4]
	\arrow[from=2-1, to=2-2]
	\arrow[from=2-2, to=2-3]
	\arrow[from=2-3, to=2-4]
	\arrow[from=2-4, to=2-5]
\end{tikzcd}\]
commutative. Thus, the tensor algebra extension is a classifier of length $1$ extensions.\footnote{recall that the length of an extension is the number of non-trivial algebras in the extension minus $2$!} More generally, the extension (where all the maps are the canonical ones)
\begin{equation}\label{eqn: length k extension}
    0 \rightarrow J^kD \rightarrow TJ^{k-1}D \rightarrow TJ^{k-2}D \rightarrow ... \rightarrow TJ^0 D \rightarrow D \rightarrow 0
\end{equation}
classifies length $k$-extensions in the obvious sense.
    \end{enumerate}
\end{example}
In the next remark, we discuss some more aspects of the cone and the tensor algebra extensions, and set up some handy notation along the way.
\begin{remark}\label{rem: lambda and kappa}
     Let $f:D \rightarrow B$ be a morphism in $\BAlg$.
     \begin{enumerate}
     \item By combining features of the two examples above, we get an induced morphism $JD \xrightarrow{\Lambda f} \Omega B$. Iterating the construction, we obtain maps $J^n D \xrightarrow{\Lambda^nf} \Omega^nB$. When $D=B$ and $f$ is the identity, we denote $\Lambda f$ by $\Lambda_B: JB \rightarrow \Omega B$. Note that $\Lambda: J \to \Omega$ is a natural transformation on $\mathbf{Ho}(\BAlg)$, since the classifying map of an extension is well-defined only up to smooth homotopy. 
    \item Denote by $\kappa^k_{D,B}: J^k(D \haotimes B) \rightarrow (J^kD )\haotimes B$ the classifying map of the length $k$-extension obtained by tensoring \Cref{eqn: length k extension} with $B$.
    \item Since $\Omega^n B \cong \Omega^n \C \haotimes B$ as mentioned in \cite[Page 92]{CMR}, the last point yields distinguished homotopy class of map
    \[
    \kappa^{n,m}: J^n\Omega^mB \rightarrow \Omega^mJ^nB.
    \]
     \end{enumerate}
\end{remark}
We next mention a technical lemma, which streamlines the effect of applying the natural transformation $\Lambda$ iteratively. This will be extremely useful in the next section. Roughly speaking, the lemma makes rigorous the plausible fact that the order of application of $\Lambda$ to a string of $J$'s and $\Omega$'s does not really matter, and yields the same natural transformation. For example, there are two potential natural transformations from $J^2$ to $\Omega^2$ given by
\[
J^2 \overset{\operatorname{id}_J(\Lambda)}{\to} J\Omega \overset{\Lambda(\operatorname{id}_\Omega)}{\to} \Omega^2
\]
and 
\[
J^2 \overset{\Lambda(\operatorname{id}_J)}{\to } \Omega J \overset{\operatorname{id}_{\Omega}(\Lambda)}{\to} \Omega^2.
\]
It can be directly checked that they are equivalent natural transformations of $\mathbf{Ho}(\BAlg)$. Indeed, naturality of $\Lambda$ yields for any morphism $f: A \rightarrow B \in \BAlg$, $\Lambda_B \circ J(f) \simeq \Omega(f) \circ \Lambda_A$. Now take $f=\Lambda_A : JA \rightarrow \Omega A$.

The authors thank Melvin Weiss for formulating the lemma in its present form.
\begin{lemma}\label{lem: melvin}
    Denote by $\langle J,\Omega \rangle_{\operatorname{free}}$ the free monoid generated by the symbols $J$ and $\Omega$ as a poset category with respect to the ordering $\geq$, where we declare $J \geq \Omega$, and have the product order on two strings of the same length. Then, for each $A \in \BAlg$, there is a functor
    \[
    \langle J,\Omega \rangle_{\operatorname{free}} \rightarrow \mathbf{Ho}(\BAlg),
    \]
    which applies a string of $J$ and $\Omega$ to the algebra $A$ on objects, and applies the natural transformation induced by applying $\Lambda$ from \emph{right to left} on the morphisms.
\end{lemma}
\begin{proof}
    The proof follows from a direct induction argument on the length of the word.
\end{proof}

Of umpteen importance for us will be the \emph{smooth Toeplitz extension}, which we discuss now.
\begin{lemma}\label{lem: toep}
    There is a semi-split extension of bornological algebras
    \begin{equation}\label{eqn: BP}
        0 \rightarrow \mathcal{K}_\s(A) \rightarrow \mathcal{T}^0_{\s}(A) \rightarrow \Omega A \rightarrow 0,
    \end{equation}
    called the smooth Toeplitz extension of $A$. The classifying $J\Omega A \rightarrow \mathcal{K}_\s(A)$ is called the \emph{Bott map}.
\end{lemma}
\begin{proof}
    All of this is contained in \cite[Section 4.1]{CMR}. We begin with the familiar Toeplitz extension of $C^*$-algebras
    \begin{equation}\label{eqn: c star toeplitz}
        0 \rightarrow \mathcal{K}_{C^*} \rightarrow \mathcal{T}_{C^*} \rightarrow C(S^1) \rightarrow 0 .
    \end{equation}
    Recall that the Toeplitz algebra $\mathcal{T}_{C^*}$ is the universal unital $C^*$-algebra on a single isometry $\iota$. We let $\mathcal{T}_{\operatorname{alg}}$ denote the $\C$-algebra generated by $\iota$ inside of $\mathcal{T}_{C^*}$, and get by restricting \Cref{eqn: c star toeplitz} an extension of $\C$-algebras
    \[
    0 \rightarrow M_\infty \rightarrow \mathcal{T}_{\operatorname{alg}} \rightarrow \C[\Z] \rightarrow 0,
    \]
    where $\C[\Z]$ denotes the complex group ring of $\Z$, i.e., the algebra of Laurent polynomials. The extension has a linear split. Thus for any bornological algebra $A$, we get a semi-split extension
    \[
    0 \rightarrow M_\infty \otimes A \rightarrow \mathcal{T}_{\operatorname{alg}}\otimes A \rightarrow \C[\Z] \otimes A \rightarrow 0,
    \]
    and an associated decomposition $\mathcal{T_{\operatorname{alg}}} \otimes A \cong M_\infty\otimes A \oplus \C[\Z]\otimes A$ as $\C$-vector spaces. Upon taking smooth completions of the right hand side\footnote{note that the idenitification of $\mathbb{C}[\Z]$ with a dense subalgebra of $C^\infty(S^1)$ is via the Fourier transform.}, we obtain a decomposition 
    \[
    \mathcal{K}_\s(A) \oplus C^{\infty}(S^1,A) =: \mathcal{T}_{\s}(A),
    \]
    which yields the extension
    \[
    0 \rightarrow \mathcal{K}_\s(A) \rightarrow \mathcal{T}_\s(A) \rightarrow C^{\infty}(S^1,A) \rightarrow 0.
    \]
    The rest follows from a routine de-unitalization trick.
\end{proof}

\begin{remark}\label{rem: stability explanation}
   In this remark we motivate our choices of stabilizations. On the one hand, \Cref{lem: toep} demonstrates the necessity of the smooth compact operators $\mathcal{K}_\s$ to have a Bott periodic $KK$-category. On the other hand, from a $C^*$-algebraic point of view, it is desirable to enforce $\mathcal{L}^p$-stability. The \(\cK_p = \sK \haotimes \mathcal{L}^p\) stabilisation therefore yields a Bott periodic $KK$-category which resembles aspects of classical Kasparov's $KK$-category. 
\end{remark}

\section{Bivariant \(K\)-theory for bornological algebras}\label{sec:bivariant-K}
In this section, we first recapitulate the bivariant $K$-theory for bornological algebras developed by Cuntz, Meyer and Rosenberg in \cite{CMR} (denoted $\KKzero$ from now on\footnote{Actually \cite{CMR} define different \(KK\)-categories which vary depending on the stabilisation used. For the moment, we suppress the choice of stabilisation in the notation.}), and prove that (modulo some set-theoretic considerations; see the next paragraph) it is a Dwyer-Kan localisation at a certain class of morphisms between bornological algebras. From this, we deduce that there exists a \emph{stable} $\infty$-categorical refinement of $\KKzero$, which we denote by $\KK$. Finally, we construct a \emph{presentable refinement} of $\KK$, thereby obtaining the category $\bKK$ which has some better formal properties than the category $\KK$.

We now discuss the technically correct approach to constructing $\KK$. We fix a large enough regular cardinal $\kappa$, and consider the essentially small category $\BAlgk$ consisting of complete bornological algebras whose cardinality is less than $\kappa$. For concreteness, take $\kappa=\aleph_2$. It can be seen that the operations described in the last section preserve $\kappa$-smallness. We denote the category of $\kappa$-small complete bornological algebras by $\BAlg_{\kappa}$.
\subsection{Recapitulation of $\KKzero$}
We first briefly recall the construction of bornological bivariant \(K\)-theory as a triangulated category as in \cite{CMR}. 

First, one forms the \textit{suspension-stable homotopy category} of bornological algebras, denoted $\Sigma\!\operatorname{Ho}$  (\cite[Section 6.3]{CMR}). It can be equipped with the structure of a triangulated category, and it is in fact the universal triangulated \textit{homology theory} for bornological algebra, i.e., the universal functor from the category of bornological algebras to a triangulated category which is smoothly homotopy invariant and excisive; this is the content of \cite[Section 6.7]{CMR}.
\begin{definition}\label{defi: Sigma Ho}
    Denote by $\Sigma\operatorname{Ho}$ the category whose objects are objects of $\BAlgk$, and whose morphisms are given by
    \[
    \Sigma\operatorname{Ho}(A,B):=\operatorname{colim}_{n \in \N} [J^nA,\Omega^n B],
    \]
    where the colimit is formed with respect to the operator $\Lambda$ as defined in \Cref{rem: lambda and kappa}.1. The composition of $J^nA \xrightarrow{f} \Omega^n B$ and $J^mB \xrightarrow{g} \Omega^m C$ is given by the homotopy class of the composite
    \begin{equation}\label{eqn: sigma ho comp}
        J^{m+n}A=J^mJ^nA \xrightarrow{J^mf} J^m\Omega^nB \xrightarrow{(-1)^{nm}.\kappa_A^{m,n}} \Omega^nJ^mB \xrightarrow{\Omega^n g} \Omega^n\Omega^mC=\Omega^{n+m}C,
    \end{equation}
    and the composition of two general morphisms in $\Sigma\operatorname{Ho}(A,B)$ and $\Sigma\operatorname{Ho}(B,C)$ is then deduced from \Cref{eqn: sigma ho comp} using the operator $\Lambda$ in the obvious way. 
    That $\Sigma\operatorname{Ho}$ is indeed a category, which is furthermore additive, is contained in \cite[Section 6.3]{CMR}. Note that there is a canonical functor
    $
    \BAlgk \rightarrow \Sigma\operatorname{Ho}.
    $
\end{definition}

We are now in a position to define the category $\KKzero$.
\begin{definition}\label{defi: KK CMR}{(cf. \cite[Definition 7.19]{CMR})}
    For any \(p \geq 1\), denote by $\KKzero_p$ the category whose objects are objects of $\BAlgk$, and whose morphisms are given by
    \[
    \KKzero_p(A,B):=\Sigma\operatorname{Ho}(\cK_p(A),\cK_p(B)),
    \]
    and the composition is given by the composition in $\Sigma\operatorname{Ho}$. There is a canonical functor
    \[
    \kkzero_p:\BAlgk \rightarrow \KKzero_p
    \]
    given by the composite $\BAlgk \rightarrow \Sigma\operatorname{Ho} \xrightarrow{{\cK_p}_*} \KKzero_p$. Similarly, we define \(\KKzero_\s\) with bornological algebras as objects, and morphisms \[
    \KKzero_\s(A,B):=\Sigma\operatorname{Ho}(\sK(A),\sK(B)).
    \] For statements and definitions that apply to \(\KKzero_p\) with \(p \geq 1\) and \(\KKzero_{\s}\), we use the notation \(\cK_?\) for the stabilisation and \(\KKzero_?\) for the associated category. A morphism in $\BAlgk$ is said to be a \emph{$\KKzero_?$-equivalence} if it is mapped to an isomorphism by $\kkzero_?$.
\end{definition}
\begin{remark}
    Actually, Cuntz-Meyer-Rosenberg refine both $\Sigma\operatorname{Ho}$ and $\KKzero_?$ to a triangulated category, see \cite[Chapters 6 and 7]{CMR}. We omit the details of their triangulated structure here, and by abusing the notation a little, use the same notation to also denote the corresponding triangulated categories. This will not cause any confusion because we will need it in one place anyway (\Cref{thm: ho KK is triangulated KKCMR}).
\end{remark}
We next single out certain morphisms in $\BAlgk$ which are $\KKzero_?$-equivalences. We use notation from the Section 3.
\begin{proposition}\label{prop: morphisms which map to isos in KK}
    The following morphisms in $\BAlgk$ are $\KKzero_?$-equivalences.
    \begin{enumerate}
        \item Any smooth homotopy equivalence $f:A \rightarrow B$.
        \item Every stabilisation homomorphism $\iota_{A,\cK^?}: A \rightarrow \mathcal{K}^{?}(A)$.
        \item Classfying maps $JC \rightarrow A$ of semi-split extensions $0 \rightarrow A \rightarrow B \rightarrow C \rightarrow 0$ whenever $B$ is contractible.
        \item The classifying maps $\operatorname{clas}(\mathcal{T}_A): J\Omega A \rightarrow \mathcal{K}_{\mathcal{S}}(A)$ of the Toeplitz extensions as in \Cref{eqn: BP}.
    \end{enumerate}
\end{proposition}
\begin{proof}
    Since $\kkzero$ factors through $\Sigma\operatorname{Ho}$, it is clearly smoothly homotopy invariant. The stability assertion follows from \cite[Lemma 7.20]{CMR}. Since the functor $\KKzero_?(D,-): \BAlgk \rightarrow \operatorname{Ab}$ is half-exact (\cite[Theorem 7.23]{CMR}), the usual five lemma argument shows Assertion 3. Note that for this argument, it was only essential that \[\KKzero_?(D,B) \cong 0 \]
    for every $D \in \BAlgk$. This is true for the reduced Toeplitz algebra (tensored with an auxilliary algebra) by \cite[Theorem 7.24]{CMR}, whence Assertion 4 follows.
\end{proof}
\begin{notation}
    In the notation of \Cref{prop: morphisms which map to isos in KK}, denote the (essentially small) classes of maps described in Assertions (1),(2),(3),(4) therein by $\mathcal{W}_h,\mathcal{W}_{\cK_?},\mathcal{W}_{SS},\mathcal{W}_{\operatorname{Bott}}$ respectively. Denote by
    \[
    \mathcal{W}_{\operatorname{min}}:=\mathcal{W}_h\cup\mathcal{W}_{\cK_?}\cup\mathcal{W}_{SS}\cup\mathcal{W}_{\operatorname{Bott}}.
    \]
    Denote by $\mathcal{W}_{KK_?}$ the (essentially small) class of $\KKzero_?$-equivalences. 
\end{notation}
\Cref{prop: morphisms which map to isos in KK} shows that $\mathcal{W}_{\operatorname{min}} \subseteq \mathcal{W}_{\KKzero_?}$. At this stage, we record an important fact about the $\KKzero_?$ category.
\begin{proposition}\label{prop: zig zag}
    Every morphism in $\KKzero_?$ is a zig-zag of (images under $\kkzero_?$ of) morphisms in $\BAlgk$, where the wrong-way maps belong to $\mathcal{W}_{\operatorname{min}}$.
\end{proposition}
\begin{proof}
    This is shown in the context of algebraic $KK$-theory in \cite[Appendix A.4]{EE}. The proof with the obvious modifications (which are in fact simplifications) work for $\KKzero_?$.

    We remark that the proof in loc. cit. takes into account the full triangulated structure of the $KK$-category, which comes handy in bookkeeping certain algebras and morphisms.
\end{proof}

Our first main result is the following.
\begin{theorem}\label{thm: mainDM}
    Let $\mathcal{W}$ be an essentially small class of morphisms such that $\mathcal{W}_{\operatorname{min}} \subseteq \mathcal{W} \subseteq \mathcal{W}_{KK_?}$. Then the functor $\kkzero_?: \BAlgk \rightarrow \KKzero_?$ presents $\KKzero_?$ as the localisation of $\BAlgk$ at $\mathcal{W}$ in the sense of ordinary categories.
\end{theorem}
The next subsection is devoted to the proof of this theorem.
\subsection{Proof of \Cref{thm: mainDM}}
We begin by manipulating the hom sets in $\KKzero_?$. By a direct cofinality argument on \Cref{defi: KK CMR}, we conclude the following:
\begin{proposition}\label{prop: kk CMR groups}
    For $A,B$ in $\BAlgk$, we have 
    \begin{align}\label{eqn: new kk formula}
   \KKzero_?(A,B) \cong \operatorname{colim}_{n \in \N}[J^{2n}\cK_?(A),\Omega^{2n}\cK_?(B)].
    \end{align}
   Given \[f_1:J^{2n_1}\cK_?(A) \rightarrow \Omega^{2n_1}\cK_?(B)\in {\KKzero}(A,B)\] and \[f_2:J^{2n_2}\cK_?(B) \rightarrow \Omega^{2n_2}\cK_?(C) \in \KKzero(B,C),\]
   the composition $\KKzero_?(A,B) \times \KKzero_?(B,C) \rightarrow \KKzero_?(A,C)$ maps it to
    \begin{multline*}
    J^{2n_2+2n_1}\cK_?(A) \xrightarrow{J^{2n_2}f_1} J^{2n_2}\Omega^{2n_1}\cK_?(B) \\
    \xrightarrow{\kappa_{\cK_?(B)}^{2n_2,2n_1}} \Omega^{2n_1}J^{2n_2}\cK_?(B)  \xrightarrow{\Omega^{2n_1}f_2} \Omega^{2n_1+2n_2}\cK_?(C).
    \end{multline*}
\end{proposition}
With \Cref{prop: kk CMR groups} at our hand, we can describe the action of the $\Lambda$ operator on the morphism sets of $\KKzero_?$ in a different way, which will help us in proving \Cref{thm: mainDM}. We begin with a preparatory lemma.
\begin{lemma} \label{lem: Bott natural}
    The following diagram in $\mathbf{Ho}(\BAlgk)$ (keeping in mind \Cref{lem: HoBAlg} and \cite[Lemma 6.17]{CMR}) commutes:
\[\begin{tikzcd}
	A && B \\
	{\mathcal{K}_{\mathcal{S}}(A)} && {\mathcal{K}_{\mathcal{S}}(B)} \\
	{J\Omega A} && {J\Omega B} \\
	{J^2A} && {\Omega^2B}
	\arrow["f", from=1-1, to=1-3]
	\arrow["{\iota_A}"', from=1-1, to=2-1]
	\arrow["{\iota_B}", from=1-3, to=2-3]
	\arrow["{\mathcal{K}_{\mathcal{S}}(f)}", from=2-1, to=2-3]
	\arrow["{\operatorname{clas}(\mathcal{T}_A)}", from=3-1, to=2-1]
	\arrow["{J\Omega(f)}", from=3-1, to=3-3]
	\arrow["{\operatorname{clas}(\mathcal{T}_B)}"', from=3-3, to=2-3]
	\arrow["{\Lambda(\operatorname{id}_{\Omega B})}", from=3-3, to=4-3]
	\arrow["{J(\Lambda(\operatorname{id}_A))}", from=4-1, to=3-1]
	\arrow["{\Lambda^k(f)}", from=4-1, to=4-3]
\end{tikzcd}\]
\end{lemma}
\begin{proof}
    The top square evidently commutes. The middle square commutes by the naturality of the Toeplitz extension of $A$. The bottom square commutes by invoking the naturality of the (iterated) path-loop extension. More precisely, since the diagram of extensions
\[\begin{tikzcd}
	{J\Omega A} & {T\Omega A} & TA & A \\
	{J\Omega B} & {T\Omega B} & PB & B \\
	{\Omega^2B} & {P\Omega B} & PB & B
	\arrow[from=1-1, to=1-2]
	\arrow["{J\Omega(f)}"', from=1-1, to=2-1]
	\arrow[from=1-2, to=1-3]
	\arrow["{T\Omega (f)}"', from=1-2, to=2-2]
	\arrow[from=1-3, to=1-4]
	\arrow[from=1-3, to=2-3]
	\arrow["f", from=1-4, to=2-4]
	\arrow[from=2-1, to=2-2]
	\arrow["{\Lambda(\operatorname{id}_{\Omega B})}"', from=2-1, to=3-1]
	\arrow[from=2-2, to=2-3]
	\arrow[from=2-2, to=3-2]
	\arrow[from=2-3, to=2-4]
	\arrow[from=2-3, to=3-3]
	\arrow["{\operatorname{id}_B}", from=2-4, to=3-4]
	\arrow[from=3-1, to=3-2]
	\arrow[from=3-2, to=3-3]
	\arrow[from=3-3, to=3-4]
\end{tikzcd}\]
    evidently commutes in $\mathbf{Ho}(\BAlgk)$.
\end{proof}

Throughout the remainder of this subsection, let $F:\BAlgk \rightarrow \mathcal{C}$ be a functor to any category which inverts $\mathcal{{W}}$. In what follows, we use the notations described in the last section. Furthermore, for endofunctors $\mathcal{F},\mathcal{G}: \mathcal{D} \rightarrow \mathcal{D}$ on any category $\mathcal{D}$ and a natural transformation $\eta: \mathcal{F} \implies \mathcal{G}$, we denote by $\eta^n: \mathcal{F}^n \implies \mathcal{G}^n$ the induced natural transformation from the monoidal structure of the endofunctor category.
\begin{definition}\label{defi: higher bott}
     We have morphisms in $\mathcal{C}$ from $F(A)$ to $F(J^2A)$, and from $F(B)$ to $F(\Omega^2B)$, given by
    \[
    \operatorname{Bott}^J_A:F(A) \xrightarrow{F(\iota_A)} F(\mathcal{K}_{\mathcal{S}}(A)) \xrightarrow{F(\operatorname{clas}(\mathcal{T}_A))^{-1}}) F(J \Omega A) \xrightarrow{F(J(\Lambda(\operatorname{id}_A)))^{-1}} F(J^2A),\]
    \[
    \operatorname{Bott}^\Omega_B:F(B) \xrightarrow{F(\iota_B)} F(\mathcal{K}_{\mathcal{S}}(B)) \xrightarrow{F(\operatorname{clas}(\mathcal{T}_B))^{-1}}) F(J \Omega B) \xrightarrow{F(\Lambda(\operatorname{id}_{\Omega B}))} F(\Omega^2B). 
    \]
    For $n \geq 1$, we also define
    \[
    \operatorname{Bott}^{J^n}_A: F(A) \xrightarrow{\operatorname{Bott}^J_A} F(J^2A) \xrightarrow{\operatorname{Bott}^J_{J^2A}} F(J^4A) \xrightarrow{\operatorname{Bott}^J_{J^4A}} ... \xrightarrow{\operatorname{Bott}^J_{J^{2n-2}A}} F(J^{2n}A),
    \]
    \[
    \operatorname{Bott}^{\Omega^n}_B: F(B) \xrightarrow{\operatorname{Bott}^\Omega_B} F(\Omega^2B) \xrightarrow{\operatorname{Bott}^\Omega_{\Omega^2B}} F(\Omega^4B) \xrightarrow{\operatorname{Bott}^\Omega_{\Omega^4B}} ... \xrightarrow{\operatorname{Bott}^\Omega_{\Omega^{2n-2}B}} F(\Omega^{2n}B).
    \]
\end{definition}
Note that there is a priori another way of defining the higher Bott maps, namely as the composite
\[
\operatorname{Bott}^{J^n}_A:F(A) \xrightarrow{F(\iota_A^n)} F(\mathcal{K}_{\mathcal{S}}^n(A)) \xrightarrow{F(\operatorname{clas}(\mathcal{T}_A)^n)^{-1}} F((J \Omega)^n A) \xrightarrow{F(J(\Lambda(\operatorname{id}_A))^n)^{-1}} F(J^{2n}A)
\]
\[
\operatorname{Bott}^{\Omega^n}_B:F(B) \xrightarrow{F(\iota_B^n)} F(\mathcal{K}_{\mathcal{S}}^n(B)) \xrightarrow{F(\operatorname{clas}(\mathcal{T}_B)^n)^{-1}} F((J \Omega)^n B) \xrightarrow{F(\Lambda(\operatorname{id}_{\Omega B})^n)} F(\Omega^{2n}B).
\]
We claim that they are the same. 
\begin{lemma}\label{lem: well defined Bott}
    The two definitions of $\operatorname{Bott}^{J^n}_A, \operatorname{Bott}^{\Omega^n}_B$ given above coincide.
\end{lemma}
\begin{proof}
    We prove it by induction on $n$. The case $n=1$ is clear.
    For the general case, consider the following diagram in $\mathcal{C}$. We omit mentioning the functor $F$ for simplicity. Also, some arrows have to be inverted to make sense of some shape commuting and to compose some morphisms; note that all the arrows in the below diagram are isomorphisms. 
\[\begin{tikzcd}[row sep=large, column sep=huge]
	A && A \\
	{\mathcal{K}_{\mathcal{S}}^{n+1}(A)} && {\mathcal{K}_{\mathcal{S}}^n(A)} \\
	{(J\Omega)^{n+1}A} & {\mathcal{K}_{\mathcal{S}}((J\Omega)^nA)} & {(J\Omega)^nA} \\
	{(J\Omega)J^{2n}A} & {\mathcal{K}_{\mathcal{S}}(J^{2n}A)} & {} \\
	{J^{2n+2}A} && {J^{2n}A}
	\arrow["{=}", from=1-1, to=1-3]
	\arrow[from=1-1, to=2-1]
	\arrow["{\operatorname{Bott}^{J^{n+1}}_A}"', bend right=55, looseness=1.25, from=1-1, to=5-1]
	\arrow[from=1-3, to=2-3]
	\arrow["{\operatorname{Bott}^{J^n}_A}", bend left=55, looseness=1.25, from=1-3, to=5-3]
	\arrow["{\mathcal{K}_{\mathcal{S}}(-)}"', from=2-3, to=2-1]
	\arrow[from=3-1, to=2-1]
	\arrow[from=3-1, to=3-2]
	\arrow[from=3-2, to=2-1]
	\arrow["{(*)}"{description}, draw=none, from=3-3, to=2-1]
	\arrow[from=3-3, to=2-3]
	\arrow["{\mathcal{K}_{\mathcal{S}}(-)}"', from=3-3, to=3-2]
	\arrow[from=4-1, to=3-1]
	\arrow[from=4-1, to=4-2]
	\arrow[from=4-2, to=3-2]
	\arrow[from=5-1, to=4-1]
	\arrow["{(**)}"', draw=none, from=5-3, to=3-2]
	\arrow[from=5-3, to=3-3]
	\arrow["{\mathcal{K}_{\mathcal{S}}(-)}"{description}, from=5-3, to=4-2]
	\arrow["{\operatorname{Bott}^J_{J^{2n}A}}"', from=5-3, to=5-1]
\end{tikzcd}\]
By induction hypothesis, the right vertical composite is the iterated composite of the $\operatorname{Bott}^J_-$ functors. Thus, it suffices to show the the outer rectangle commutes. The top rectangle evidently commutes. The left triangle commutes by the definition of $\operatorname{clas}(\mathcal{T}_A^n)$. The rectangle on the left commutes by the naturality of the classifying map of the Toeplitz extension. The bottom left trapezoid commutes by the definition of $\operatorname{Bott}^J_{J^{2n}A}$. The trapezoids marked $(*)$ and $(**)$ commute by the functoriality of $\mathcal{K}_{\mathcal{S}}(-)$.

This shows the assertion for $\operatorname{Bott}^{J^n}_A$. The assertion about $\operatorname{Bott}^{\Omega^n}_B$ follows analogously.
\end{proof}
Immediately from \Cref{lem: Bott natural} and \Cref{lem: well defined Bott}, we conclude:
\begin{corollary}\label{cor: Lambda^2 formulation}
    In $\mathcal{C}$, we have \[F(\Lambda^{2n}(f)) = \operatorname{Bott}^{\Omega^n}_B \circ F(f) \circ \operatorname{Bott}^{J^n}_A{^{-1}}:F(J^{2n}A) \rightarrow F(\Omega^{2n}B).\] In particular, this holds for the functor $\kkzero_?$.
\end{corollary}
We now present a technical result involving the $\operatorname{Bott}^{J^n}_-$ and $\operatorname{Bott}^{\Omega^n}_-$ functors, which will come handy later. We need a preparatory lemma.
\begin{lemma}\label{lem: technical kappa}
    The diagram
\[\begin{tikzcd}
	& {J^{n+m}B} \\
	{J^n\Omega^mB} && {\Omega^mJ^nB}
	\arrow["{J^n(\Lambda^m_B)}"'{pos=0.8}, from=1-2, to=2-1]
	\arrow["{(-1)^{mn}\Lambda^m_{J^nB}}"{pos=0.8}, from=1-2, to=2-3]
	\arrow["{\kappa_B^{n,m}}", from=2-1, to=2-3]
\end{tikzcd}\]
commutes in $\KKzero_?$.
\end{lemma}
\begin{proof}
    Consider the following diagram, where each of the rectangles commute by \cite[Eqn. 6.33]{CMR}
\[\begin{tikzcd}
	{J^{n+m}B} && {J^{n+m}B} \\
	{J^{n+m-1}\Omega B} && {\Omega J^{n+m-1}B} \\
	{...} && {...} \\
	{J^{n+2}\Omega^{m-2}B} && {\Omega^{m-2}J^{n+2}B} \\
	{J^{n+1}\Omega^{m-1}B} && {\Omega^{m-1}J^{n+1}B} \\
	{J^n\Omega^mB} && {\Omega^mJ^nB}
	\arrow["{\kappa^{n+m,0}_B=\operatorname{id}}", from=1-1, to=1-3]
	\arrow["{J^{n+m-1}\Lambda_B}"', from=1-1, to=2-1]
	\arrow["{(-1)^{n+m-1}\Lambda_{J^{n+m-1}_B}}", from=1-3, to=2-3]
	\arrow["{\kappa^{n+m-1,1}}", from=2-1, to=2-3]
	\arrow[from=2-1, to=3-1]
	\arrow[from=2-3, to=3-3]
	\arrow[from=3-1, to=4-1]
	\arrow[from=3-3, to=4-3]
	\arrow["{\kappa^{n+2,m-2}_B}", from=4-1, to=4-3]
	\arrow["{J^{n+1}\Lambda_{\Omega^{m-2}B}}"', from=4-1, to=5-1]
	\arrow["{(-1)^{n+m-1}\Omega^{m-2}\Lambda_{J^{n+1}B}}", from=4-3, to=5-3]
	\arrow["{\kappa_B^{n+1,m-1}}", from=5-1, to=5-3]
	\arrow["{J^n(\Lambda_{\Omega^{m-1}B})}"', from=5-1, to=6-1]
	\arrow["{(-1)^{n+m-1}\Omega^{m-1}\Lambda_{J^nB}}", from=5-3, to=6-3]
	\arrow["{\kappa^{n,m}_B}", from=6-1, to=6-3]
\end{tikzcd}\]
Now, the left vertical composite is $J^n(\Lambda^m_B)$, whereas the right vertical composite is $(-1)^{m(n+m-1)}\Lambda^m_{J^nB}$, as required. 
\end{proof}
\begin{proposition}
    \label{lem: technical Bott}
    The diagram
\[\begin{tikzcd}
	{F(\Omega^{2n}B)} && {F(J^{2n}B)} \\
	{F(J^{2n}\Omega^{2n}B)} && {F(\Omega^{2n}J^{2n}B)}
	\arrow["{\operatorname{Bott}^{J^n}_{\Omega^{2n}B}}"', from=1-1, to=2-1]
	\arrow["{\Lambda^{2n}(\operatorname{id})}"', from=1-3, to=1-1]
	\arrow["{\operatorname{Bott}^{\Omega^n}_{J^{2n}B}}", from=1-3, to=2-3]
	\arrow["{\kappa^{2n,2n}_B}", from=2-1, to=2-3]
\end{tikzcd}\]
commutes.
\end{proposition}
\begin{proof}
     Unraveling the Bott maps in the diagram, we get the following diagram with the obvious arrows, each of which are isomorphisms
\[
\begin{tikzcd}[row sep=large, column sep=huge]
	{F(\Omega^{2n}B)} && {F(J^{2n}B)} \\
	{F(\mathcal{K}_{\mathcal{S}}^n(\Omega^{2n}B))} && {F(\mathcal{K}_{\mathcal{S}}^n(J^{2n}B))} \\
	{F((J\Omega)^n\Omega^{2n}B)} && {F((J\Omega)^nJ^{2n}B)} \\
	& {F(J^{4n}B)} \\
	{F(J^{2n}\Omega^{2n}B)} && {F(\Omega^{2n}J^{2n}B)}
	\arrow[from=1-1, to=2-1]
	\arrow["{\operatorname{Bott}^{J^n}_{\Omega^{2n}B}}"', bend right=55, looseness=1.25, from=1-1, to=5-1]
	\arrow["{\Lambda^{2n}_B(\operatorname{id})}"', from=1-3, to=1-1]
	\arrow[from=1-3, to=2-3]
	\arrow["{\operatorname{Bott}^{\Omega^n}_{J^{2n}B}}", bend left=55, looseness=1.25, from=1-3, to=5-3]
	\arrow[from=2-3, to=2-1]
	\arrow[from=3-1, to=2-1]
	\arrow[from=3-3, to=2-3]
	\arrow[from=3-3, to=3-1]
	\arrow[from=3-3, to=5-3]
	\arrow[from=4-2, to=3-1]
	\arrow[from=4-2, to=3-3]
	\arrow[from=4-2, to=5-1]
	\arrow[from=4-2, to=5-3]
	\arrow[from=5-1, to=3-1]
	\arrow["{\kappa^{2n,2n}_B}"', from=5-1, to=5-3]
\end{tikzcd}
\]
It suffices to show that each of the triangles involving $F(J^{4n}B)$ commute. The bottom triangle commutes from \Cref{lem: technical kappa}. The other triangles commute as a consequence of \Cref{lem: melvin}.
\end{proof}
We can finally prove \Cref{thm: mainDM}.
\begin{proof}[Proof of \Cref{thm: mainDM}]
    Let $\mathcal{D}$ be any $1$-category, and $F:\BAlgk \rightarrow \mathcal{D}$ be a functor which sends $\mathcal{W}$ to isomorphisms. We have to show that there is a unique functor $\overline{F}: \KKzero_? \rightarrow \mathcal{D}$ such that $\overline{F} \circ \KKzero_? = F$ (see \cite[Lemma 1.2.2]{HOVEY}, \cite[Lemma 1.2]{GZ}).

    Since $\KKzero_?: \BAlgk \rightarrow \KKzero_?$ is identity on objects, we are forced to define \(\overline{F}(A):=F(A)\) for all \(A \in \BAlgk\). For a discussion of morphisms, we use \Cref{prop: kk CMR groups}. Given \[f:J^{2n}\cK_?(A) \rightarrow \Omega^{2n}\cK_?(A)\] (viewed as a morphism in $\KKzero_?(A,B)$), we assign to it the morphism in $\mathcal{C}$ given as the composite
    \begin{multline}\label{eqn: definition of Fbar}
    \small
        F(A) \xrightarrow{F(\iota_A)} F(\cK_?(A)) \xrightarrow{\operatorname{Bott}^{J^n}_{\cK_? (A)}} F(J^{2n} \cK_? (A)) \\
        \xrightarrow{F(f)} F(\Omega^{2n}\cK_?(B))  \xrightarrow{(\operatorname{Bott}^{\Omega^n}_{\cK_? (B)})^{-1}} F(\cK_?(B)) \xrightarrow{F(\iota_B)^{-1}} F(B).
    \end{multline}

    Since the functor $F$ is smoothly homotopy invariant, this yields a well defined map
    \begin{equation}\label{eqn: functor kk to c}[J^{2n}\cK_?(A),\Omega^{2n}\mathcal{K}_?(B)] \rightarrow \mathcal{C}(A,B).
    \end{equation}
    Note that the assignment in \Cref{eqn: functor kk to c} is compatible with the inductive system of \Cref{eqn: new kk formula}, as readily follows from Lemma \ref{lem: Bott natural}
         and Corollary \ref{cor: Lambda^2 formulation}. Thus, we get a well-defined assignment \begin{equation}\label{eqn: functor out of kk}
             \overline{F}_{A,B}:\KKzero_?(A,B) := \operatorname{colim}_{n \in \N}[J^{2n}\mathcal{K}_?(A),\Omega^{2n}\mathcal{K}_?(B)] \rightarrow \mathcal{C}(A,B).
         \end{equation}
        It remains to prove that the assignment in \ref{eqn: functor out of kk} yields a functor $\overline{F}:\BAlgk  \rightarrow \mathcal{D}$.

        We begin by showing that \[\forall f \in \KKzero_?(A,B), \forall g \in \KKzero_?(B,C), \overline{F}_{B,C}(g) \circ \overline{F}_{A,B}(f)=\overline{F}_{A,C}(g \circ f).\] Indeed, let \[f:J^{2n}\mathcal{K}_?(A) \rightarrow \Omega^{2n}\mathcal{K}_?(B)\] and \[g:J^{2m}\mathcal{K}_?(B) \rightarrow \Omega^{2m}\mathcal{K}_?(C)\] be any two elements in $\KKzero_?(A,B)$ and $\KKzero_?(B,C)$ respectively. We can assume $n=m$ without loss of generality. Consider the following diagram
\[\begin{tikzcd}[column sep=small]
	{F(A)} && {F(B)} && {F(B)} && {F(C)} \\
	{F(\cK_?(A))} && {F(\cK_?(B))} && {F(\cK_?(B))} && {F(\cK_?(C))} \\
	{F(J^{2n}\cK_?(A))} && {F(\Omega^{2n}\cK_?(B))} && {F(J^{2n}\cK_?(B))} && {F(\Omega^{2n}\cK_?(C))} \\
	&&& {} \\
	{F(J^{4n}\cK_?(A))} && {F(J^{2n}\Omega^{2n}\cK_?(B))} && {F(\Omega^{2n}J^{2n}\cK_?(B))} && {F(\Omega^{4n}\cK_?(C))}
	\arrow["{:=\overline{F}_{A,B}(f)}", color={rgb,255:red,61;green,150;blue,245}, dashed, from=1-1, to=1-3]
	\arrow["{F(\iota_A)}", color={rgb,255:red,61;green,150;blue,245}, from=1-1, to=2-1]
	\arrow["{=}", no head, from=1-3, to=1-5]
	\arrow["{F(\iota_B)}", color={rgb,255:red,61;green,150;blue,245}, from=1-3, to=2-3]
	\arrow["{:=\overline{F}_{B,C}(g)}", color={rgb,255:red,214;green,92;blue,92}, dashed, from=1-5, to=1-7]
	\arrow["{F(\iota_B)}", color={rgb,255:red,214;green,92;blue,92}, from=1-5, to=2-5]
	\arrow["{F(\iota_C)}"', color={rgb,255:red,214;green,92;blue,92}, from=1-7, to=2-7]
	\arrow["{\operatorname{Bott}^{J^n}_{\cK_?(A)}}", color={rgb,255:red,61;green,150;blue,245}, from=2-1, to=3-1]
	\arrow["{\operatorname{Bott}^{\Omega^n}_{\cK_?(B)}}", color={rgb,255:red,61;green,150;blue,245}, from=2-3, to=3-3]
	\arrow["{\operatorname{Bott}^{J^n}_{\cK_?(B)}}", color={rgb,255:red,214;green,92;blue,92}, from=2-5, to=3-5]
	\arrow["{\operatorname{Bott}^{\Omega^n}_{\cK_?(C)}}"', color={rgb,255:red,214;green,92;blue,92}, from=2-7, to=3-7]
	\arrow["{F(f)}", color={rgb,255:red,61;green,150;blue,245}, from=3-1, to=3-3]
	\arrow["{\operatorname{Bott}^{J^n}_{J^{2n}\cK_?(A)}}"{description}, from=3-1, to=5-1]
	\arrow["{\operatorname{Bott}^{J^n}_{\Omega^{2n}\cK_?(B)}}"{description}, from=3-3, to=5-3]
	\arrow["{F(\Lambda^{2n}(\operatorname{id}_{\cK_?(B)}))}", from=3-5, to=3-3]
	\arrow["{F(g)}", color={rgb,255:red,214;green,92;blue,92}, from=3-5, to=3-7]
	\arrow["{\operatorname{Bott}^{\Omega^n}_{J^{2n}\cK_?(B)}}"{description}, from=3-5, to=5-5]
	\arrow["{\operatorname{Bott}^{\Omega^n}_{\Omega^{2n}\cK_?(C)}}"{description}, from=3-7, to=5-7]
	\arrow["{F(J^{2n}f)}"', from=5-1, to=5-3]
	\arrow["{F(\kappa^{2n,2n}_{\cK_?(B)})}"', from=5-3, to=5-5]
	\arrow["{F(\Omega^{2n}g)}", from=5-5, to=5-7]
\end{tikzcd}\]
As described in \Cref{eqn: definition of Fbar}, the dashed arrows are so defined so as to make the rectangles with blue and red edges commute. In order to prove our thesis, we have to show that the outer square commutes (note that the bottom row is precisely the formula of composition in $\KKzero_?$). The bottom-left and bottom-right rectangles commute by the naturality of $\operatorname{Bott}^{J^{n}}_-$ and $\operatorname{Bott}^{\Omega^{n}}_-$ respectively. The bottom-middle rectangle commutes by \Cref{lem: technical Bott}. The top-middle rectangle commutes by \Cref{cor: Lambda^2 formulation}.

The only check remaining in order to verify that $\overline{F}:\KKzero_? \rightarrow \mathcal{D}$ is indeed a functor is unitality. But that is clear from \Cref{eqn: definition of Fbar} (take $n=0$ therein). 

It remains to prove that $\overline{F}:\KKzero_? \rightarrow \mathcal{D}$ is the unique functor such that $F=\overline{F} \circ \KKzero_?$. But this follows from Proposition \ref{prop: zig zag}.
\end{proof}
\begin{remark}
    We mention that the complexity of our proof of Theorem \ref{thm: mainDM} is strictly between that of the analogous proof for $C^*$-algebraic $KK$ (as done in \cite[Proposition 2.1]{BEL}) and algebraic $KK$ (as done in \cite[Theorem 4.24]{EE}). 
    
    On the one hand, we had to carefully handle the colimit definition of the hom sets in $\KKzero_?$ (for which we crucially use Lemma \ref{lem: Bott natural}), which was not a problem in $C^*$-case, thanks to the Cuntz-picture of Kasparov's $KK$-theory. On the other hand, we had at our disposal the Bott map $J\Omega A \rightarrow \mathcal{K}_\s(A)$, which evidently got exploited in the proof. This is not available in Corti\~nas-Thom's algebraic $KK$-theory as homotopy \(K\)-theory does not satisfy Bott periodicity.
\end{remark}

\subsection{A stable $\infty$-categorical refinement of $\KKzero$}
With \Cref{thm: mainDM} at hand, we can now proceed to construct $\infty$-categorical refinements of $\KKzero_?$. Our treatment is parallel to \cite[Section 2]{BEL}.

We begin by defining the $\infty$-category $\KK_?$.
\begin{definition}
    Denote by $\KK_p$ the small $\infty$-category obtained by Dwyer-Kan localising $\BAlgk$ at $\mathcal{W}_{\KKzero_?}$, i.e.,
    \[
    \KK_? \defeq \BAlgk[\mathcal{W}_{\KKzero_?}^{-1}].
    \]
    Let \[\kk_?:\BAlgk \rightarrow \KK_?\] denote the canonical functor.
\end{definition}
\comment{\begin{remark}
    We note that $\KK$ could equivalently defined as the Dwyer-Kan localisation of $\BAlgk$ at any of the class of morphisms ${\cW}$ such that $\cW_{\operatorname{min}} \subseteq {\cW} \subseteq \cW_{\KKzero_?}$. This follows from the \Cref{thm: mainDM}. However, as we will see down the road, it is \emph{more readily} seen that $\BAlgk[\cW_{\KKzero_?}^{-1}]$ is a stable $\infty$-category. Compare with \cite{bunke2024kk}, where Bunke had to perform the stepwise localisations beginning with the category of separable $C^*$-algebras in a specific order to prove that the small $KK$-category he constructs is indeed stable.
\end{remark}}
Since every ordinary category can be viewed as an $\infty$-category via the nerve functor, by the defining universal property of $\KK_?$, we have a functor $h:\KK_? \rightarrow \KKzero_?$ rendering the following diagram commutative.
\[
\begin{tikzcd}
	\mathbf{\BAlgk} & & \KKzero_? \\
	\KK_p
	\arrow["\kk_?", from=1-1, to=1-3]
	\arrow["\kk_?", from=1-1, to=2-1]
	\arrow["h"', curve={height=12pt}, dashed, from=2-1, to=1-3]
\end{tikzcd}
\]
From basic properties of homotopy categories, we get that $h$ descends to a functor between ordinary categories 
$
\operatorname{ho}(h):\operatorname{ho}(\KK_?) \rightarrow \KKzero_?.
$
\begin{proposition}\label{lem: ho of KKCMR INFTY is KKCMR0}{(cf. \cite[Proposition 2.15.2]{BEL})}
    The functor
    \[
    \operatorname{ho}(h):\operatorname{ho}(\KK_?) \rightarrow \KKzero_?
    \]
    is an equivalence.
\end{proposition}
\begin{proof}
    This follows from the following general fact: let $\mathcal{C}$ be an $\infty$-category and $W$ a class of morphisms. Then we have a commutative square 
\[\begin{tikzcd}
	\mathcal{C} & \mathcal{C}[W^{-1}] \\
	\operatorname{ho}(\mathcal{C}) & {\operatorname{ho}(\mathcal{C}[W^{-1}])}
	\arrow["l", from=1-1, to=1-2]
	\arrow[from=1-1, to=2-1]
	\arrow[from=1-2, to=2-2]
	\arrow["{\operatorname{ho}(l)}", from=2-1, to=2-2]
\end{tikzcd}\]
where $\operatorname{ho}(l): \operatorname{ho}(\mathcal{C}) \rightarrow \operatorname{ho}(\mathcal{C}[W^{-1}])$ presents $\operatorname{ho}(\mathcal{C}[W^{-1}]$ as the localisation of $\operatorname{ho}(\mathcal{C})$ in the sense of \(1\)-categories. In our context, we take $\mathcal{C}=\BAlgk$ and $W=\cW_{\KKzero_?}$. Since $\BAlgk$ is an ordinary category, the morphism 
\[
\mathcal{C} \rightarrow \operatorname{ho}(\mathcal{C})
\]
in the above diagram becomes an equivalence, whence we have that the composite
\[
\BAlgk \xrightarrow{\kk_?} \KK_? \rightarrow \operatorname{ho}(\KK_?) 
\]
presents $\operatorname{ho}(\KK_?)$ as the localisation of $\BAlgk$ at $\cW_{\KKzero_?}$-equivalences in the sense of ordinary categories. We are now done by the universal property of $\KKzero_?$ that was shown in \Cref{thm: mainDM}.
\end{proof}

\comment{\begin{corollary}\label{cor: kk eq saturated}
    The functor
    \[
    \kk_?:\BAlgk \rightarrow \KK_?
    \]
    precisely inverts the maps in $\mathcal{W}_{\KKzero_?}$
\end{corollary}

\begin{proof}
    Let $f:A \rightarrow B$ be a morphism in $\BAlgk$ which gets mapped to an equivalence in $\KK_?$. Then by definition, $f$ gets mapped to an isomorphism in $\operatorname{ho}(\KK_?)$, and from Lemma \ref{lem: ho of KKCMR INFTY is KKCMR0}, also to an isomorphism in $\KKzero_?$. But then by definition, $f \in \mathcal{W}_{\KKzero_?}$.
\end{proof}}

\begin{proposition}\label{prop: functor to addivtie cat with adjectives inverts kk eq}
    A functor from $\BAlgk$ to an additive $\infty$-category admitting finite limits which is reduced, smoothly homotopy invariant, $\cK_?$-stable and \emph{semi-exact} (i.e., sends semi-split exact sequences in $\BAlgk$ to fiber sequences) sends $\KKzero_?$-equivalences to equivalences.
\end{proposition}
\begin{remark}
    Note that a functor as above is automatically reduced; see \cite[Remark 2.3]{BEL}.
\end{remark}
\begin{proof}[Proof of \Cref{prop: functor to addivtie cat with adjectives inverts kk eq}]
    Let 
    \[
    F: \BAlgk \rightarrow \mathcal{C}
    \]
    be a functor as in the statement of the proposition.
    We show that each of the classes of morphisms $\mathcal{W}_h,\mathcal{W}_{\cK_p}, \mathcal{W}_{SS}, \mathcal{W}_{\operatorname{Bott}}$ are sent to equivalences by $F$. Morphisms in $\mathcal{W}_h, \mathcal{W}_{\cK_p}$ are sent to equivalences by definition. A standard 5-lemma argument shows that morphisms in $\cW_{SS}$ are also sent to equivalences. That $\mathcal{W}_{\operatorname{Bott}}$ is inverted follows since $F$ is automatically additive, using a direct $\infty$-categorical version of \cite[Theorem 7.24]{CMR} (cf. \cite[Section 6]{bunke2024kk}). 

    Thus, from \Cref{thm: mainDM}, we get a unique factorisation of $F$ through $\KKzero_?$ as depicted by the dashed arrow below.
\[\begin{tikzcd}
	\mathbf{\BAlgk} & \mathcal{C} \\
	\KKzero_? & \operatorname{ho}(\mathcal{C})
	\arrow["F", from=1-1, to=1-2]
	\arrow["{\kk_?}", from=1-1, to=2-1]
	\arrow["{\operatorname{ho}(-)}", from=1-2, to=2-2]
	\arrow["{\overline{F}}", dashed, from=2-1, to=2-2]
\end{tikzcd}
\]
But since the functor $\kkzero_?$ also inverts $\KKzero_?$-equivalences (by definition), we get that both the functors $\BAlgk \rightarrow \operatorname{ho}(\mathcal{C})$ invert $\KKzero_?$-equivalences. Since the canonical functor $\mathcal{C} \rightarrow \operatorname{ho}(\mathcal{C })$ detects equivalences, the result follows.
\end{proof}
\begin{remark}\label{rem: inf loc kk satisfies properties}
    At this point, it is not clear if the localisation functor 
    \[
    \kk_? : \BAlgk \rightarrow \KK_?
    \]
    is a functor satisfying the requirements of \Cref{prop: functor to addivtie cat with adjectives inverts kk eq} (that its smoothly homotopy invariant and $\cK_?$-stable are clear though). Even more, it is not clear if the category $\KK_?$ has finite limits. We begin preparations to demonstrate this now.
\end{remark}
We next want to prove that $\KK_?$ is a stable $\infty$-category. To that end, we consider the functor
\begin{align*}
    \Omega: \BAlgk &\rightarrow \BAlgk \\
    A &\mapsto C^{\infty}_{\operatorname{flat}}(0,1) \haotimes A,
\end{align*}
 which turns out to be an explicit model for the loops functor in $\KK_?$.
\begin{lemma} \label{lem: loop descends}
    \hfill
    \begin{enumerate}
        \item The functor $\Omega$ descends essentially uniquely to an equivalence of additive categories such that 
\[\begin{tikzcd}
	{\BAlgk} & {\BAlgk} \\
	\KKzero_? & \KKzero_?
	\arrow["\Omega", from=1-1, to=1-2]
	\arrow["\kkzero_?", from=1-1, to=2-1]
	\arrow["\kkzero_?", from=1-2, to=2-2]
	\arrow["{\Omega_0}", dotted, from=2-1, to=2-2]
\end{tikzcd}\]
commutes. Furthermore, $\Omega_0:\KKzero_? \rightarrow \KKzero_?$ is self-inverse.
\item The functor $\Omega$ descends essentially uniquely  to a functor (also denoted by) $\Omega$ rendering the square
\[\begin{tikzcd}
	{\BAlgk} & {\BAlgk} \\
	\KK_? & \KK_?
	\arrow["\Omega", from=1-1, to=1-2]
	\arrow["\kk_?"', from=1-1, to=2-1]
	\arrow["\kk_?", from=1-2, to=2-2]
	\arrow["\Omega", dotted, from=2-1, to=2-2]
\end{tikzcd}\]
commutative. $\Omega: \KK_? \rightarrow \KK_?$ is self-inverse too.
    \end{enumerate}
\end{lemma}
\begin{proof}
    For Assertion $1$, we first show that \[
    \Omega(\mathcal{W}_h \cup \mathcal{W}_{\cK_?} \cup \mathcal{W}_{SS} \cup \mathcal{W}_{\operatorname{Bott}}) \subseteq \mathcal{W}_{\KKzero_?}.
    \]
    Clearly, $\Omega$ preserves the classes $\mathcal{W}_h$ and $\mathcal{W}_{\cK_?}$, and in particular, maps them into $\mathcal{W}_{\KKzero_?}$. Now let 
    \begin{equation}\label{eqn: extn}
          0 \rightarrow A \rightarrow B \rightarrow C \rightarrow 0
    \end{equation}
    be a semi-split extension with $B$ contractible, and let $f:JC\rightarrow A$ be its classifying map. We have to show that the induced map $\Omega (JC) \xrightarrow{\Omega f} \Omega A$ is a $\KKzero_?$-equivalence. To that end, consider the following commutative diagram
\[\begin{tikzcd}
	0 & {J(\Omega C)} & {T(\Omega C)} & {\Omega  C} & 0 \\
	0 & {\Omega  (J C)} & {\Omega (TC)} & {\Omega C} & 0 \\
	0 & {\Omega A} & {\Omega B} & {\Omega C} & 0
	\arrow[from=1-1, to=1-2]
	\arrow[from=1-2, to=1-3]
	\arrow[from=1-2, to=2-2]
	\arrow[from=1-3, to=1-4]
	\arrow[from=1-3, to=2-3]
	\arrow[from=1-4, to=1-5]
	\arrow["{=}", no head, from=1-4, to=2-4]
	\arrow[from=2-1, to=2-2]
	\arrow[from=2-2, to=2-3]
	\arrow["{\Omega(f)}"', from=2-2, to=3-2]
	\arrow[from=2-3, to=2-4]
	\arrow[from=2-3, to=3-3]
	\arrow[from=2-4, to=2-5]
	\arrow["{=}", no head, from=2-4, to=3-4]
	\arrow[from=3-1, to=3-2]
	\arrow[from=3-2, to=3-3]
	\arrow[from=3-3, to=3-4]
	\arrow[from=3-4, to=3-5]
\end{tikzcd}\]
The composite left vertical composite $J(\Omega C) \rightarrow \Omega A$ is the classifying map of the loop of the extension in \Cref{eqn: extn}, and top left vertical arrow $J(\Omega C) \rightarrow \Omega JC$ is the classifying map of the extension in the middle row\footnote{we are implicitly using the universal property of the tensor algebra to assert that this diagram commutes.}. Thus, these two morphisms are in $\mathcal{W}_{SS}$, and in particular in $\mathcal{W}_{\KKzero_?}$. Since $\mathcal{W}_{\KKzero_?}$ is closed under 2-out-of-3, we conclude that the map $\Omega (JC) \xrightarrow{\Omega f} \Omega A$ is also $\KKzero_?$-equivalence. A similar argument also applies to the maps in $\mathcal{W}_{\operatorname{Bott}}$.

This shows the existence of the functor $\Omega_0:\KKzero_? \rightarrow \KKzero_?$. It is additive because both $\Omega$ and $\kkzero_?$ are additive. That it is also an equivalence with it being its own inverse follows from Bott periodicity. Assertion 2 follows directly from Assertion 1. 
\end{proof}

Recall the definition of a category fibrant objects as in \cite[Definition 2.7]{BEL}. Our next result is crucial for better understanding of the category $\KK_?$.
\begin{proposition}\label{prop:fibration-kk}
    The triple \((\BAlgk,\cF,\cW_{\KKzero_?})\) forms a category of fibrant objects structure on $\BAlgk$, where $\cF$ denotes the class of bounded algebra homomorphisms that admit bounded linear sections. 
\end{proposition}

\begin{proof}
The proof of \cite[Proposition 3.11]{EE} goes through verbatim, in light of \cite[Section 6.6]{CMR}.
\end{proof}
\begin{proposition}\label{prop: KKCMR stable and related topics}
\hfill
    \begin{enumerate}
        \item $\KK_?$ is a pointed, left-exact $\infty$-category.
        \item 
        The loop functor on $\KK_?$ (in the sense of the loop functor on a pointed left-exact $\infty$-category, see \cite[Section 1.1.2]{lurie:HA} ) is naturally equivalent to the functor $\Omega$ as constructed in \Cref{lem: loop descends}.2.
        \item The functor 
        \[
        \kk_?: \BAlgk \rightarrow \KK_?
        \]
        is semi-exact.
        \item The category $\KK_?$ is stable, and the loop functor on $\KK_?$ is self-inverse.
    \end{enumerate}
\end{proposition}
\begin{proof}
    Assertion 1 is an immediate consequence of consequence of the fact that $\KK_?$ is the $\infty$-category associated to a category of fibrant objects, see \cite[Remark 7.1.15 and Proposition 7.5.6]{cisinski2020higher}.  Assertion 2 follows from the existence of the path-loop extension
    \[
    C_{\operatorname{flat}}^{\infty}((0,1),A) \rightarrow C_{\operatorname{flat}}^{\infty}((0,1],A) \xrightarrow{\operatorname{ev}_1} A,
    \]
    and using the fact that the path algebra $C_{\operatorname{flat}}^{\infty}((0,1],A)$ is contractible and $\operatorname{ev}_1$ admits a bounded linear section, whence is a fibration in our chosen category of fibrant objects structure on $\BAlgk$.  Assertion 3 is again an immediate consequence of the fact that $\KK_?$ is the $\infty$-category associated to a category of fibrant objects structure on $\BAlgk$ whose fibrations are the semi-split surjections.
    We now prove Assertion 4. To show that $\KK_?$ is stable, it suffices to show (\cite[Proposition 1.4.2.27]{lurie:HA}) that it has finite limits, is pointed, and the loop functor on $\KK_?$ is an autoequivalence. We are just left with showing the last assertion. From \cite[Lemma 3.4(2)]{land2018}, it suffices to show that \[\operatorname{ho}(\Omega): \operatorname{ho}(\KK_?) \rightarrow \operatorname{ho}(\KK_?)\]
    is an equivalence of ordinary categories. We are now done from \Cref{lem: ho of KKCMR INFTY is KKCMR0} and \Cref{lem: loop descends}.1.
\end{proof}
Note that the above proposition resolves the issue mentioned in \Cref{rem: inf loc kk satisfies properties}. We now present a lemma containing some technical corollaries to the existence of the category of fibrant objects structure on $\BAlgk$ as considered in \Cref{prop:fibration-kk}.
\begin{lemma}\label{lem: maps equivalent to fibrations}
    \hfill
    \begin{enumerate}
        \item Every morphism in $\KK_?$ is equivalent to a morphism $\kk_?(f)$ for some morphism $f$ in $\BAlgk$.
        \item Every binary product in $\KK_?$ is equivalent to the image under $\kk_?$ of a product in $\BAlgk$.
        \item Every fiber sequence in $\KK_?$ is equivalent to the image under $\kk_?$ of a mapping cone sequence.
    \end{enumerate}
\end{lemma}
\begin{proof}
    Proofs of Assertion 1 and 2 follow mutatis-mutandis from the relevant proofs in \cite[Proposition 5.9]{EE}. The proof of Assertion 3 is analogous to the corresponding proof in \cite[Proposition 2.17.3]{BEL}.
\end{proof}
\begin{theorem}\label{thm: ho KK is triangulated KKCMR}
    The functor
        \[
        \operatorname{ho}(h):\operatorname{ho}(\KK_?) \rightarrow \KKzero_?
        \]
        is an equivalence of triangulated categories.
\end{theorem}
\begin{proof}
    We have already seen in \Cref{lem: ho of KKCMR INFTY is KKCMR0} that the homotopy category of $\KK_?$ is equivalent to $\KKzero_?$ as \(1\)-categories. Since both the functors \[\kkzero_?: \BAlgk \rightarrow \KKzero_?\] and \[\kk_?:\BAlgk \rightarrow \KK_?\]preserve finite products, so does the comparison functor \[ \operatorname{ho}(h):\operatorname{ho}(\KK_?) \rightarrow \KKzero_?.\] By Lemma \ref{lem: loop descends}, we conclude that $\operatorname{ho}(h)$ also commutes with the shift functor.

    On the one hand, the exact triangles in $\KKzero_?$ are generated by the mapping cone sequences for the morphisms in $\BAlgk$. The same is also true in $\operatorname{ho}(\KK_?)$ in light of Lemma \ref{lem: maps equivalent to fibrations}. We conclude that $\operatorname{ho}(h)$ is compatible with the collection of distinguished triangles, and this finishes the proof.
\end{proof}
The next result strengthens \Cref{prop: functor to addivtie cat with adjectives inverts kk eq}, and will also be important for us in the next section. Before proceeding, we need to set up some notation.
\begin{notation}
    For any $\infty$-category $\mathcal{D}$, let 
    \[
    \mathbf{Fun}^{r,h,\cK_?,SS}(\BAlgk,\mathcal{D})
    \]
    be the full subcategory of $\mathbf{Fun}(\BAlgk,\mathcal{D})$ comprising of functors which are reduced, smoothly homotopy invariant, $\cK_?$-stable and semi-split exact. 
\end{notation}
\begin{proposition}\label{prop: Fun lex KKCMR C}\label{cor: uni prop stable inf cat KKcmr}
    Pullback along the functor 
    \[
    \kk_?: \BAlgk \rightarrow \KK_?
    \]
    induces an equivalence of functor categories
    \[
    \mathbf{Fun}^{\operatorname{lex}}(\KK_?,\mathcal{C}) \xrightarrow[\simeq]{{\kk}^*} \mathbf{Fun}^{r,h,\cK_?,SS}(\BAlgk,\mathcal{C})
    \]
    for every additive $\infty$-category $\mathcal{C}$ admitting finite limits.
\end{proposition}
\begin{proof}
    We build the following diagram
\begin{equation}\label{eqn: diagram with adjectives}
    \begin{tikzcd}
	{\mathbf{Fun}^{\operatorname{lex}}(\KK_?,\mathcal{C})} & {\mathbf{Fun}^{r,h,\mathcal{K}_?,SS}(\BAlgk,\mathcal{C})} \\
	{\mathbf{Fun}^{\operatorname{}}(\KK_?,\mathcal{C})} & {\mathbf{Fun}^{\kk_?}(\BAlgk,\mathcal{C})}
	\arrow[dashed, from=1-1, to=1-2]
	\arrow[hook, from=1-1, to=2-1]
	\arrow[hook, from=1-2, to=2-2]
	\arrow["\simeq", from=2-1, to=2-2]
\end{tikzcd}
\end{equation}
The vertical arrows are inclusions of full subcategories (using \Cref{prop: functor to addivtie cat with adjectives inverts kk eq}) and the horizontal ones are induced by precomposition with $\kk_?$. The bottom horizontal arrow is an equivalence from the definition of $\KK_?$ as the Dwyer-Kan localisation of $\BAlgk$ at $\mathcal{W}_{\KKzero_?}$. We now justify the existence of the dashed arrow rendering the diagram commutative. To that end, note that the functor 
\[
\kk_?:\BAlgk \rightarrow \KK_?
\]
is reduced, smoothly homotopy invariant, $\cK_?$-stable and semi-exact, as mentioned after the proof of \Cref{prop: KKCMR stable and related topics}. The first three of these properties are clearly preserved by postcomposing with a left-exact functor $\overline{F}: \KK_? \rightarrow \mathcal{C}$. To see that semi-exactness would also be preserved, we use that semi-split extensions give rise to fiber sequences in $\KK_?$ (\Cref{lem: maps equivalent to fibrations}.3), and then use the left-exactness of $\overline{F}$.

We now proceed to show that the dashed arrow is an equivalence. Since all the bold arrows are fully faithful, so is the dashed arrow. So we just need to show that the dashed arrow is essentially surjective. So let 
\[
F \in \mathbf{Fun}^{r,h,\cK_?,SS}(\BAlgk,\mathcal{C})
\]
be any functor. In view of the  bottom horizontal arrow, there exists an $\overline{F} \in \mathbf{Fun}(\KK_?,\mathcal{C})$ such that $\overline{F} \circ \kk_? \simeq F$. It remains to prove that $\overline{F}$ is left-exact. 

It suffices to show that $\overline{F}$ preserves binary products, fiber sequences, and is reduced. Since both $\kk_?$ and $F$ are reduced, so is $\overline{F}$. That $\overline{F}$ preserves binary products follows since both $\kk$ and $F$ preserve binary products (the latter because smoothly homotopy invariant and reduced automatically implies additivity of $F$). By \Cref{lem: maps equivalent to fibrations}.3, every fiber sequence in $\KK$ is the image of a semi-split exact sequence in $\BAlgk$. By assumption, $F$ maps this sequence to a fiber sequence. In view of the relation $\overline{F} \circ \kk_? \simeq F$, we conclude $\overline{F}$ maps the image of this sequence under $\kk$ to a fiber sequence. Consequently, $\overline{F}$ preserves fiber sequences.
\end{proof}
\begin{proposition}\label{prop:symmetric-monoidal}
    There is a symmetric monoidal structure on \(\KK\) such that the localisation \(j \colon \BAlg \to \KK\) has a symmetric monoidal refinement. Moreover, for any symmetric monoidal \(\infty\)-category \(\bD\), we have an equivalence \[j^* \colon \mathbf{Fun}_{\otimes/\mathrm{lax}}(\KK, \bD) \xrightarrow[j^*]{\simeq} \mathbf{Fun}_{\otimes/\mathrm{lax}}^{\cW_{\CMR}}(\KK, \bD),\]
    where the notation $\mathbf{Fun}_{\otimes / \mathrm{lax}}$ stands for two formulae, one for symmetric monoidal and the other for lax monoidal functors.
\end{proposition}

\begin{proof}
    Since \(\KK\) is a Dwyer-Kan localisation of \(\BAlg\) at the collection \(\cW_{\CMR}\), to see that it has a symmetric monoidal structure and that the localisation functor \(j\) has a symmetric monoidal refinement, it suffices by \cite[Definition 3.2.2]{hinich2013dwyer} to check that \(A \haotimes -\) preserves maps in \(\cW_\CMR\). To see this, we have for the map \(B \to B \haotimes C^\infty([0,1])\) in \(\cW_\CMR\), \(A \haotimes B \to A \haotimes B \haotimes C^\infty([0,1]) \in \cW_\CMR\), so that \(A \haotimes -\) preserves smooth homotopy equivalences. For the map \(B \to \cK \haotimes B\), we again have \(A \otimes B \to A \haotimes B \haotimes \cK \in \cW_\CMR\), so that \(A \haotimes - \) preserves corner inclusions into smooth compact operators. Finally, if \(B \to C \to D\) is a semi-split extension of bornological algebras, then since \(C \cong B \oplus D\) as bornological vector spaces, we have a semi-split extension \(A \haotimes B \to A \haotimes C \oplus A \haotimes D\) of bornological algebras.  
\end{proof}
\subsection{A presentable refinement of $\KK$}\label{subsec:presentable-bivariant-K}

It is not clear whether $\KK_?$ has countable coproducts. Indeed, if an $\infty$-category admits coproducts indexed by a set $I$, then so does its homotopy category (\cite[Proposition 6.5.5]{cisinski2020higher}). Since 
\[
\operatorname{ho}(\KK_?) \simeq \KKzero_?,
\]
we would also have countable coproducts in $\KKzero_?$. But this is not clear, as mentioned at the end of Section 13.1.1 in \cite{CMR}\footnote{strictly speaking, loc. cit. states that its not clear whether countable direct sums of algebras are a model of countable coproducts in $\KKzero$.}.

In this section, we construct  a \emph{presentable refinement} of \(\KK_?\) as follows. We first take the Yoneda embedding into the ind-completion of $\KK_?$  \[y\colon \KK_? \to \mathbf{Ind}_{\aleph_0}(\KK_?).\] 

The $\infty$-category \(\mathbf{Ind}_{\aleph_0}(\KK_?)\) has all colimits (since $\KK_?$ is stable). However, the composition \[\BAlgk \xrightarrow{\kk_?} \KK_? \overset{y}\to  \mathbf{Ind}_{\aleph_0}(\KK_?)\] need not preserve countable filtered colimits, so this needs to be forced.
\begin{notation}
    Let $\mathcal{C},\mathcal{D}$ be $\infty$-categories. $\mathbf{Fun}_L(\mathcal{C},\mathcal{D})
    $ denote the full subcategory of $\mathbf{Fun}(\mathcal{C},\mathcal{D})$ spanned by those functors which are left adjoints. For a small set of morphisms $S$ in $\mathcal{C}$, $\mathbf{Fun}^S_L(\mathcal{C},\mathcal{D})
    $ denotes the full subcategory of $\mathbf{Fun}_L(\mathcal{C},\mathcal{D})
    $ spanned by those functors which send morphisms in $S$ to equivalences in $\mathcal{D}$.
\end{notation}
\begin{proposition}\label{prop: BL at a set of maps}{\cite[Proposition 5.5.4.20]{lurie2009higher}}
    Let $\mathcal{C}$ be a presentable $\infty$-category, and $S$ a small set of morphisms of $\mathcal{C}$. Then there exists a (unique upto equivalence) presentable $\infty$-category $\mathcal{C}[S^{-1}]_{BL}$, called \emph{the Bousfield localisation of $\mathcal{C}$ at $S$}, together with a cocontinuous functor
    \[
    L:\mathcal{C} \rightarrow \mathcal{C}[S^{-1}]_{BL}
    \]
    satisfying the following universal property: for any $\infty$-category $\mathcal{D}$, we have an equivalence of $\infty$-categories
    \[
    \mathbf{Fun}_L(\mathcal{C}[S^{-1}]_{BL},\mathcal{D}) \xrightarrow[\simeq]{L^*} \mathbf{Fun}_L^{S}(\mathcal{C},\mathcal{D})
    \]
\end{proposition}
\comment{\begin{proof}
    In the notation of \cite[Proposition 5.5.4.20]{lurie2009higher}, the category $\mathcal{C}[S^{-1}]_{BL}$ can be taken to be $S^{-1}\mathcal{C}$, the full subcategory of $\mathcal{C}$ spanned by \textit{$S$-local} objects, i.e., objects $Z \in \mathcal{C}$ such that for every morphism $s: X \rightarrow Y$ in $S$, the induced map
    \[
    s^*:\operatorname{Map}_{\mathcal{C}}(Y,Z) \rightarrow \operatorname{Map}_{\mathcal{C}}(X,Z)
    \]
    is an equivalence in Anima, see \cite[Definition 5.5.4.1]{lurie2009higher} and the paragraph above \cite[Proposition 5.5.4.16]{lurie2009higher}. $S^{-1}{\mathcal{C}}$ is a presentable $\infty$-category and the inclusion $S^{-1}\mathcal{C \hookrightarrow \mathcal{C}}$ admits a left-adjoint, as follows from \cite[Proposition 5.5.4.15]{lurie2009higher}. This gives us the sought after cocontinuous functor
    \[
    L:\mathcal{C} \rightarrow \mathcal{C}[S^{-1}]_{BL}.
    \]
    That it also satisfies the mentioned universal property is precisely the content of \cite[Proposition 5.5.4.20]{lurie2009higher}.
\end{proof}}

\begin{remark}\label{rem: BL of pointed cat}
    If $\mathcal{C}$ is a pointed, presentable $\infty$-category, then the Bousfield localisation of $\mathcal{C}$ at any small set of maps $S$ is also a pointed, presentable $\infty$-category. Indeed, \Cref{lem: bunke seven} shows that the Bousfield localisation of $\mathcal{C}$ at $S$ has an initial and a terminal object, given by $L(0_{\mathcal{C}})$ and $0_{\mathcal{C}}$ respectively. We also have a morphism $0_{\mathcal{C}} \rightarrow L(0_{\mathcal{C}})$ in $S^{-1}\mathcal{C}$ courtesy $0_{\mathcal{C}}$ being the initial object in $\mathcal{C}$. We can now conclude by \cite[Remark 1.1.1.2]{HA}.
\end{remark}
\begin{definition}\label{def:new-KK}
    We define \emph{bivariant \(K\)-theory for bornological algebras} as the Bousfield localisation \[\bKK_? \defeq \cInd(\KK_?)[{{\cW_{\aleph_0}}}^{-1}]_{BL},\]
    where $\cW_{\aleph_0}$ is the (essentially small) class of maps in $\cInd(\KK_?)$ of the form
    \[
    \operatorname{colim}_{n \in I} y(\kk_?(A_n)) \rightarrow y(\kk_?(\operatorname{colim}_{n \in I}A_n)),
    \]
    where \(I\) is a countable directed set and $A_n$'s form a filtered system of ($\kappa$-small) bornological algebras with injective structure maps. 
Denote by
\[
\bkk_? \colon \BAlgk\overset{\kk_?}\to \KK_? \overset{y}\to \cInd(\KK_?) \overset{L}\to \bKK_? \]
 the canonical functor defined as the composition of the universal functor, the Yoneda embedding and the Bousfield localisation at ${{\cW_{\aleph_0}}}$.
\end{definition}
 
 To formulate the universal property of \(\bKK_?\), we introduce the following notation: for an additive $\infty$-category $\mathcal{D}$, denote by \[\mathbf{Fun}^{r,h,\cK_?,SS,\aleph_0}(\BAlgk, \mathcal{D})\] the full subcategory of $\mathbf{Fun}^{r,h,\cK_?,SS}(\BAlgk, \mathcal{D})$ spanned by functors which preserve countable filtered colimits of bornological algebras along injective structure maps. \comment{To avoid notational clutter, we also use the shorthand 
 \[
 \BAlgk \xrightarrow{j} \KK_?
 \]
 to denote $\kk_?$.}

\begin{proposition}\label{prop:presentable-stable}
We have the following:
\hfill
  \begin{enumerate}
        \item $\cInd(\KK_?)$ is a presentable, stable $\infty$-category, and the loop functor on $\cInd(\KK_?)$ is given by 
        \[
        \cInd(\KK_?) \xrightarrow{\cInd(\Omega_{\KK_?})} \cInd(\KK_?).
        \]
        Furthermore, it is self-inverse.
      \item \(\bKK_?\) is a presentable stable \(\infty\)-category, and the loop functor on $\bKK_?$ is self-inverse.
      \item The pullback \((\bkk)^{*}\) induces an equivalence \[\mathbf{Fun}_{L}(\bKK_?, \mathcal{D}) \overset{\simeq}\to \mathbf{Fun}^{{r,h,\cK_?,SS,\aleph_0}}(\BAlgk, \mathcal{D})\] for all stable cocomplete \(\infty\)-categories \(\mathcal{D}\). 
  \end{enumerate}  
\end{proposition}
\begin{proof}

We begin with Assertion 1. Presentability of $\cInd(\KK_?)$ is clear; stability follows from \cite[Proposition 1.1.3.6]{lurie:HA}. Now, the proof of loc. cit. shows that the suspension functor on $\cInd(\KK_?)$ is 
    \[
    \cInd(\KK_?) \xrightarrow{\cInd(\Sigma_{\KK_?})} \cInd(\KK_?).
    \]
    But since $\Omega_{\KK_?}$ is self-inverse by \Cref{prop: KKCMR stable and related topics}.4, stability of $\KK_?$ yields $\Sigma_{\KK_?}$ is also self-inverse. Functoriality of the $\cInd(-)$ functor shows that $\cInd(\Sigma_{\KK_?}) \simeq \Sigma_{\cInd(\kk_?)}$ is self-inverse as well. So from the stability of $\cInd(\KK_?)$, we have equivalences 
    \[
    \Omega_{\cInd(\KK_?)} \simeq \Sigma_{\cInd(\KK_?)} \simeq \cInd(\Sigma_{\KK_?}) \simeq \cInd(\Omega_{\KK_?}),
    \]
    as desired. That $\Omega_{\cInd(\KK_?)}$ is self-inverse follows from the preceding discussion. This completes the proof of Assertion 1. Presentability of $\bKK_?$ follows from \Cref{prop: BL at a set of maps}. To show that $\bKK_?$ is stable, we use \cite[Corollary 1.4.2.27(3)]{lurie:HA}, and show that $\bKK_?$ is pointed, left-exact, and 
    \[
    \Omega: \bKK_? \rightarrow \bKK_?
    \]
    is an equivalence. That $\bKK_?$ is pointed follows from \Cref{rem: BL of pointed cat}. $\bKK_?$ is presentable and therefore in particular left-exact. It remains to show that the loops functor on $\bKK_?$ is an equivalence. We note that $\Omega_{\cInd(\KK_?)}$ descends to $\bKK_?$, as follows from \Cref{lem: bunke seven}. But since Assertion 1 shows that \[\Omega^2_{\cInd(\KK_?)} \simeq \operatorname{id}_{\cInd(\KK_?),}\]
    the same conclusion continues to hold in $\bKK_?$. Thus Assertion 2 follows from the trivial fact that if the composition of an endofunctor on a category with itself is equivalent to the identity functor, then the functor is an equivalence. To show Assertion 3, consider the following diagram
    
    \[
\begin{adjustbox}{max width=\textwidth}
\begin{tikzcd}
	{\mathbf{Fun}_L(\bKK_?,\mathcal{D})}
	&
	{\mathbf{Fun}_L^{{{\mathcal{W}_{\aleph_0}}}}
	(\mathbf{Ind}_{\aleph_0}(\KK_?,\mathcal{D})}
	&&
	{\mathbf{Fun}^{r,h,\cK_?,SS,\aleph_0}
	(\mathbf{BAlg}_\kappa,\mathcal{D})}
	\\
	&
	{\mathbf{Fun}_L(\mathbf{Ind}_{\aleph_0}(\KK_?,\mathcal{D})}
	&
	{\mathbf{Fun}^{\operatorname{lex}}
	(\KK_?,\mathcal{D})}
	&
	{\mathbf{Fun}^{r,h,\mathcal{K}_?,SS}(\mathbf{BAlg}_\kappa,\mathcal{D})}
	\arrow["{(!)}", from=1-1, to=1-2]
	\arrow[dashed, from=1-2, to=1-4]
	\arrow[hook', from=1-2, to=2-2]
	\arrow[hook', from=1-4, to=2-4]
	\arrow["{(!!)}", from=2-2, to=2-3]
	\arrow["{(!!!)}", from=2-3, to=2-4]
\end{tikzcd}
\end{adjustbox}
\]
where all the bold arrows are obvious ones; $(!)$ is an equivalence by \Cref{prop: BL at a set of maps}, $(!!)$ is an equivalence by basic properties of ind-completions of stable $\infty$-categories (see \cite[Section 8]{bunke2024kk}), $(!!!)$ is an equivalence by \Cref{prop: Fun lex KKCMR C}, and the vertical arrows are inclusion of full subcategories. We now justify the existence of the dashed arrow rendering the diagram commutative. Observe that the lower horizontal composition is induced by the composite
\[
\BAlgk \xrightarrow{\kk_?} \KK \xrightarrow{y} \cInd(\KK_?).
\]
Let $\tilde{F} \in \mathbf{Fun}_L^{{{\mathcal{W}_{\aleph_0}}}}(\cInd(\KK_?),\mathcal{D})$ be any functor. Then, for any diagram \[A: I \rightarrow \BAlg_{\kappa}\] in ${{\mathcal{W}_{\aleph_0}}}$, we have equivalences
\begin{equation}\label{eqn: defining prop bkk}
    \operatorname{colim}_{n \in N}\tilde{F}yj(A_n) \xrightarrow{\simeq} \tilde{F}(\operatorname{colim}_{n \in N}yj(A_n)) \xrightarrow{\simeq} \tilde{F}yj\operatorname{colim}_{n \in N}A_n,
\end{equation}
 where the first equivalence comes from $\tilde{F}$ being a left-adjoint (so preserves all colimits), and the second equivalence follows from the fact that $\tilde{F}$ inverts ${{\mathcal{W}_{\aleph_0}}}$. The existence of the dashed arrow follows easily from this. 

We now proceed to show that the dashed arrow is an equivalence. Since all the bold arrows in the diagram are fully faithful, so is the dashed one. So we just need to show that it is essentially surjective. So let 
\[
F \in \mathbf{Fun}^{r,h,\mathcal{K}_?,SS,\aleph_0}(\BAlgk,\mathcal{D})
\]
be any functor. In view of the  bottom horizontal composite, there exists an $\tilde{F} \in \mathbf{Fun}_L(\cInd(\KK_?),\mathcal{D})$ such that $\tilde{F} \circ yj \simeq F$. It remains to prove that $\tilde{F}$ inverts the maps in ${{\mathcal{W}_{\aleph_0}}}$. Let $A: I \rightarrow \BAlg_\kappa$ be an arbitrary member of ${{\mathcal{W}_{\aleph_0}}}$. Consider the commutative triangle in $\mathcal{D}$ with the obvious morphisms
\begin{equation}\label{eqn: zigzag in KKborn}\begin{tikzcd}
	{\tilde{F}(\operatorname{colim}_{n \in I}yj (A_n))} && {\tilde{F}(yj(\operatorname{colim}_{n \in I} A_n))} \\
	{\operatorname{colim}_{n \in I}\tilde{F}yj(A_n)}
	\arrow[from=1-1, to=1-3]
	\arrow[from=2-1, to=1-1]
	\arrow[from=2-1, to=1-3].
\end{tikzcd}
\end{equation}
We have to show that the top horizontal morphism is an equivalence in $\mathcal{D}$. To that end, it suffices to show that the other two morphisms are equivalences. The vertical morphism is an equivalence because $\tilde{F}$ preserves colimits, and the diagonal morphism is an equivalence since $\tilde{F}\circ yj \simeq F$. This completes the proof of Assertion 3.
\end{proof}
\begin{remark}
    Since $\KK_?$ has a symmetric monoidal refinement, so does the category $\mathbf{Ind}_{\aleph_0}(\KK_?)$. However, it is not clear whether $\bKK$ has a reasonable symmetric monoidal structure, since $\cW_{\aleph_0}$ is not closed under tensoring with an arbitrary auxilliary bornological algebra. This will be fixed in the following section.
\end{remark}

\section{Symmetric monoidal structures and the nuclear bivariant category}
In this section, we show that the adjunction in \Cref{eqn: NBA adj} descends to an adjunction \(\bKK_\s \leftrightarrows \nKK\), whereby the category \(\nKK\) has the last feature we desire, namely, a symmetric monoidal structure descended from algebras. As before, we work with the full subcategory \(\NBA_\kappa \subset \BAlg_\kappa\) of \(\kappa\)-small strongly nuclear bornological algebras.

\begin{theorem}\label{thm:universal-nuc}
    There is a presentable, stable \(\infty\)-category \(\nKK\) and a functor \(\nkk \colon \NBA_\kappa \to \nKK\) such that pullback along \(\nkk\)  
    \[
    (\nkk)^* \colon \Fun^{\colim}(\nKK, \bD) \overset{\simeq}\to \Fun^{r,h,SS, \cK_{\mathcal{S}},\aleph_0}(\NBA_\kappa, \bD) 
    \] is an equivalence of \(\infty\)-categories for every presentable, stable $\infty$-category $\mathcal{D}$.
\end{theorem}

\begin{proof}{(cf. \cite[Theorem 2.10]{bunke2026theory})}
We obtain $\nKK$ by a series of Bousfield localisations of $\mathbf{PSh}(\NBA_{\kappa},\Sp)$, the stabilisation of the category of (Anima valued) presheaves on $\NBA_{\kappa}$. For any presentable, stable $\infty$-category $\mathcal{D}$, we have equivalences
\begin{multline*}
     \mathbf{Fun}^{\operatorname{colim}}(\mathbf{PSh}(\NBA_{\kappa},\Sp),\mathcal{D}) 
    \xrightarrow[\simeq]{{\Sigma^\infty}^*} \mathbf{Fun}^{\operatorname{colim}}(\mathbf{PSh}(\NBA_{\kappa}),\mathcal{D}) \\
    \xrightarrow[\simeq]{y^*} \mathbf{Fun}(\NBA_{\kappa},\mathcal{D}),
\end{multline*}
where the first equivalence follows from \cite[Corollary 1.4.4.5]{HA} and the second one from \cite[Theorem 5.1.5.6]{HTT}.

In $\mathbf{PSh}(\NBA_{\kappa},\Sp)$, consider the small class of morphisms given by
\begin{enumerate}
    \item $\mathcal{W}_1:= (\Sigma^\infty \circ y)\{f \in \NBA_{\kappa} \mid  f ~\text{is a smooth homotopy equivalence}\}$.
    \item $\mathcal{W}_2:=(\Sigma^\infty \circ y)\{A \rightarrow A \hat\otimes \mathcal{K}_{\mathcal{S}} \mid A \in \NBA_{\kappa}\}$
    \item $\mathcal{W}_3:=\{0 \rightarrow (\Sigma^\infty \circ y)(0)\}$
\end{enumerate}
We now take the Bousfield localisation of $\mathbf{PSh}(\NBA_{\kappa},\Sp)$ at $\mathcal{W}:=\mathcal{W}_1 \cup \mathcal{W}_2 \cup \mathcal{W}_3$. Denote by \[L_0:\mathbf{PSh}(\NBA_{\kappa},\Sp) \rightarrow \mathbf{PSh}(\NBA_{\kappa},\Sp)[\mathcal{W}^{-1}]_{BL}.\]
To avoid notational clutter, we denote $\mathbf{PSh}(\NBA_{\kappa},\Sp)[\mathcal{W}^{-1}]_{BL}$ by $\mathbf{LOC}$ henceforth\footnote{The notation $\mathbf{LOC}$ stands for \emph{local objects}.}. Then we have equivalences
\begin{align*}
    \mathbf{Fun}^{\operatorname{colim}}(\mathbf{LOC}) \xrightarrow[L_0^*]{\simeq} \mathbf{Fun}^{\operatorname{colim},\mathcal{W}}(\mathbf{PSh}(\NBA_{\kappa},\Sp), \mathcal{D}) \simeq \mathbf{Fun}^{h,\mathcal{K}_\mathcal{S},\operatorname{red}}(\NBA_{\kappa},\mathcal{D}),
\end{align*}
where the last equivalence can be deduced analogous to \Cref{eqn: diagram with adjectives}.

We next consider the small set of maps 
\begin{multline*}
\mathcal{W}_{SS}:=\{L_0\Sigma^{\infty} y(A) \longrightarrow \operatorname{fib}(L_0\Sigma^{\infty} y(B) \rightarrow L_0\Sigma^{\infty} y(C)) \\
\mid 
0 \rightarrow A\rightarrow B \rightarrow C \rightarrow 0 ~\text{is semi-split exact in } \NBA_{\kappa}\},
\end{multline*}
and take the Bousfield localisation
\[
L_1:\mathbf{LOC} \rightarrow \mathbf{LOC}[\mathcal{W}_{SS}^{-1}]_{BL}
\]
Then we have equivalences
\[
\mathbf{Fun}^{\operatorname{colim}}(\mathbf{LOC}[\mathcal{W}_{SS}^{-1}]_{BL},\mathcal{D}) \simeq \mathbf{Fun}^{\operatorname{colim},\mathcal{W}_{SS}}(\mathbf{LOC},\mathcal{D}) \simeq \mathbf{Fun}^{h,\mathcal{K}_{\mathcal{S}} ,\operatorname{red},SS}(\NBA_{\kappa},\mathcal{D}),
\]
where the last equivalence again follows from \Cref{eqn: diagram with adjectives}.

For the final step, we consider the small set of maps
\[
\mathcal{W}_{\aleph_0}:= \operatorname{colim}_{n \in \N} L_1L_0\Sigma^{\infty} y(A_n) \rightarrow L_1L_0\Sigma^{\infty} y(\operatorname{colim}_{n \in \N}A_n)),
\]
where $A_n$'s are in $\NBA_{\kappa}$, and the structure maps are injective. We now define \[
\nKK:= \mathbf{LOC}[(\mathcal{W}_{SS} \cup \mathcal{W}_{\aleph_0})^{-1}]_{BL},~~~~\nkk:= L_2L_1L_0\Sigma^{\infty} y:\NBA_{\kappa} \rightarrow \nKK, 
\]
where
\[
L_2: \mathbf{LOC}[\mathcal{W}_{SS}^{-1}]_{BL} \rightarrow \mathbf{LOC}[(\mathcal{W}_{SS} \cup \mathcal{W}_{\aleph_0})^{-1}]_{BL}.
\]
Similar arguments as above show that $\nKK$ and $\nkk$ indeed satisfy the universal property mentioned in the statement of the theorem. 

\end{proof}

\begin{theorem}\label{thm:bousfield-nuc}
    The right Bousfield localisation \((-)^\rig :\BAlg \leftrightarrows \NBA_{\kappa}: \iota\) descends to a right Bousfield localisation 
    \[
    (-)^\rig \colon \bKK_\s \leftrightarrows \nKK: \iota.
    \]
\end{theorem}

\begin{proof}
We only need to show that the functors in the adjunction \((\iota, (-)^\rig)\) descend to functors between \(\bKK_\s\) and \(\nKK\). To see this, we first need to show that the functor \(\nkk \circ (-)^\rig\) is homotopy invariant, \(\cK_{\mathcal{S}}\)-stable, excisive and filtered-colimit-preserving. Consider the maps \(A \to C^\infty([0,1],A)\) and \(A \to \cK_\s(A)\). By Proposition \ref{prop:infinitely-nuc}, since \(C^\infty[0,1]\) and \(\cK_\s\) are strongly nuclear, we have \begin{multline*}
A^\rig \to (C^\infty([0,1],A)^\rig \cong C^\infty([0,1],A^\rig) \cong C^\infty[0,1] \haotimes A^\rig,\\
A^\rig \to (\cK_{\s} \haotimes A)^\rig \cong \cK_\s \haotimes A^\rig = \sK(A^\rig)
\end{multline*} which are mapped to equivalences by \(\nkk\). Excision for semi-split extensions is clear. Finally, filtered-colimit preservation again follows from \ref{prop:infinitely-nuc}. So by Theorem \ref{prop:presentable-stable}, the functor \(L\) descends to \(\bKK_\s\).  

In the other direction, the functor \(\bkk \circ \iota \colon \NBA_{\kappa} \to \bKK\) is homotopy invariant, \(\mathcal{K}_{\mathcal{S}}\)-stable and filtered colimit-preserving by \ref{prop:infinitely-nuc}. Excision for semi-split extension follows from the observation that kernels and cokernels in \(\mathsf{Nuc}^\infty\) are formed underlying in \(\CBorn\), so that the former is an extension closed subcategory of the latter.
\end{proof}

We finally conclude with the promised symmetric monoidal structure on our bivariant nuclear category. 

\begin{proposition}\label{prop:symmetric-monoidal}
    There is a symmetric monoidal structure on \(\nKK\) such that the localisation \(\nkk \colon \NBA_{\kappa} \to \nKK\) has a symmetric monoidal refinement. Moreover, for any symmetric monoidal \(\infty\)-category \(\bD\), we have an equivalence \[(\nkk)^* \colon \mathbf{Fun}_{\otimes/\mathrm{lax}}(\nKK, \bD) \overset{\simeq}\to \mathbf{Fun}_{\otimes/\mathrm{lax}}^{r,h,SS,\mathcal{K}_{\mathcal{S}},\kappa}(\NBA_\kappa, \bD)\] induced by restriction.  
\end{proposition}

\begin{proof}
This follows directly from the steps in the proof of \Cref{thm:universal-nuc}, together \cite[Proposition 3.3.2]{HIN} and \cite[Lemma 14.5.2]{Cnossen2026}.
\end{proof}

\section{Relationship with \(K\)-theories for topological algebras}

In this section, we study the functor corepresented by \(\C\) in our category \(\bKK_{?}\) and $\nKK$. It turns out that in general, this functor is the extension of the diffeotopy \(K\)-theory functor defined in \cite{cortinas2008comparison} to general bornological algebras. We then clarify in Remark \ref{subsec:remarks} how the corepresented functor \(\KKzero(\C,-)\) relates to other \(K\)-theories for topological algebras.

\comment{\textcolor{purple}{In this section, we want to study the functor corepresented by $\C$ in the different $KK$-categories studied before. It turns out that in general, this functor is (a stabilised version of) the so-called diffeotopy $K$-theory functor....... Our main results of this section are \Cref{thm:comparision-kd} and \Cref{cor:corepresentable}....At the end of this section, we will see that when $?=\mathcal{K}_{\mathcal{L}^p}$ and the algebra is an LBA, then we have an equivalence of spectra 
\[
\KK_{\cK^p}(\C,A) \simeq \K^{\operatorname{top}}(A),
\]
the familiar topological K-theory spectrum of a local Banach algebra.
}}
\subsection{Diffeotopy \(K\)-theory for bornological algebras}

Consider again the cone extenion \[\Omega(A) \to C(A) \to A\] from Example \ref{ex: ex of ext}. By excision in negative algebraic \(K\)-theory \cite[Proposition 2.5.2 (ii)]{Cortinas2011}, we have natural maps \begin{equation}\label{eq:negative-K}
K_n(A) \to K_{n-1}(\Omega(A)),
\end{equation} for all \(n \leq 0\), where \(K_{-n}(A) = K_0(\Sigma_{\mathrm{alg}}^n(A))\). Here \(\Sigma_{alg}\) is part of the algebraic Calkin extension 
\[
M_\infty \to \Gamma \to \Sigma_{\mathrm{alg}}.
\] We equip \(\Sigma_{\mathrm{alg}}\) with the fine bornology, so that \(\Sigma_{\mathrm{alg}}(A) \cong \Sigma_{\mathrm{alg}} \haotimes A = \Sigma_{\mathrm{alg}} \otimes A\), as a consequence of which we have \[\Omega\Sigma_{\mathrm{alg}}(A) \cong \Sigma_{\mathrm{alg}}\Omega(A).\]

By iterating the maps in Equation \ref{eq:negative-K}, we get a nonconnective version of topological \(K\)-theory in our setting:
\begin{definition}\label{def:KH-theory}
The \textit{diffeotopy \(K\)-theory} of a complete bornological algebra \(A\) is defined as the abelian groups 
\[K_n^\an(A) \defeq \colim_{m:n+m\geq 0} K_0(\Sigma_{\mathrm{alg}}^m \Omega^{n+m}(A)),\] for \(n \in \Z\). 
\end{definition}

\comment{We can also express the groups \(K_n^\an(A)\) in terms of the diffeotopic \(KV\)-groups. To see this, we replace \(A\) by \(\coma{\Sigma}(A)\) in Equation \ref{eqref:KV-K0} and compose with the map \(K_0(-) \onto KV_1^\an(\coma{\Sigma}(-))\) to get \[KV_1^\an(\coma{\Sigma}(A)) \to K_0(\coma{\Sigma}S(A)) \to KV_1^\an(\coma{\Sigma}^2 S(A)) = KV_2^\an(\coma{\Sigma}^2(A)).\] Iterating, we get an isomorphic inductive system whose colimit is again 
\[K_n^\an(A) \cong \colim_m KV_1^\an(\coma{\Sigma}^{m+1} S^{n+m}(A)) = \\colim_m KV_{n+m}^\an( \coma{\Sigma}^m(A))\] for \(n \in \Z\).}

We now express \(K_n^\an\) as the homotopy groups of a spectrum-valued functor. Consider the simplices \[\mathbb{A}^n = \setgiven{(x_0, \dotsc, x_n) \in \R^{n+1}}{\sum_{i=0}^n x_i = 1} \supset \Delta^n = \setgiven{(x_0,\dotsc, x_n) \in \mathbb{A}^n}{x_i \geq 0, 0 \leq i \leq n}.\] Let \(C^\infty(\Delta^n)\) be the locally convex vector space of smooth functions on \(\Delta^n\) that are restrictions of smooth functions on \(\mathbb{A}^n\), equipped with the precompact bornology. Then for a complete bornological vector space \(V\), we define \[\Delta_n^{\an}(V) \defeq  C^\infty(\Delta^n, V).\] For \(M = A\) a complete, unital bornological algebra, \(\Delta^{\an}(A) \defeq [n] \mapsto \Delta_n^{\an}(A)\) is a simplicial unital bornological algebra, and we define the spectrum \[\K^\an(A) \defeq \abs{\K(\Delta^{\an}(A))} = \colim_{[n] \in \Delta^{\op}}\K(\Delta_n^\an(A)).\] For non-unital algebras, define \[\K^\an(A) \defeq \mathsf{fib}(\K^\an(A^{+}) \to \K^{\an}(\C)).\] It is then easy to see (for example by the same argument as that in \cite[Theorem 10.3.2]{Cortinas2011} that \(\pi_n(\K^\an(A)) \cong K_n^\an(A)\). We now summarise the properties of diffeotopy \(K\)-theory as follows:

\begin{theorem}\label{thm:KH-theory-properties}
The functor \(\K^\an \colon \BAlg \to \mathbf{Sp}\) satisfies the following properties:
\begin{enumerate}
\item Smooth homotopy invariance;
\item Excision for semi-split extension of complete bornological algebras;
\item Filtered colimit-preservation for filtered systems of bornological algebras with injective structure maps.
\end{enumerate}
\end{theorem}

\begin{proof}
    Smooth homotopy invariance and excision follow from the construction of \(\K^\an\). For filtered colimit-preservation, let \(A = \colim A_i\) be a filtered system with injective structure maps. Since filtered colimits are termwise, and \[C^\infty(\Delta^n, \colim_i A_i) \cong C^\infty(\Delta^n) \haotimes \colim_i A_i \cong \colim_i C^\infty(\Delta^n, A_i)\] by nuclearity, we have \(\Delta^{\mathrm{diff}}(A) \cong \colim_i \Delta^{\mathrm{diff}}(A_i)\). Moreover, again by nuclearity, tensoring by \(C^\infty(\Delta^n)\) preserves monomorphisms, so that for each \([n] \in \Delta\) and \(i \leq j\), the maps \(\Delta_n^{\mathrm{diff}}(A_i) \to \Delta_n^{\mathrm{diff}}(A_j)\) are monomorphisms, so that the colimit \(\colim_i \Delta^{\mathrm{diff}}(A)\) is essentially reduced, and is therefore also a filtered colimit of underlying rings. Consequently, we have an equivalence \[\K(\Delta^{\mathrm{diff}}(\colim_i A_i)) \simeq \colim_i \K(\Delta^{\mathrm{diff}}(A_i)).\] Since geometric realisation commutes with colimits, the conclusion follows.  
\end{proof}

We now compare diffeotopy \(K\)-theory for bornological algebras with the functor \(\KD \colon \mathsf{LCAlg} \to \mathbf{Sp}\) on complete locally convex algebras considered in \cite{cortinas2008comparison}. Recall that we have a functor \[\mathsf{Cpt} \colon \mathsf{LC} \to \mathsf{CBorn}\] that equips a complete locally convex topological vector space with the precompact bornology. Recall that this is always lax monoidal so that it in particular yields a functor \(\mathsf{Cpt} \colon \mathsf{LCAlg} \to \BAlg\). 

\begin{theorem}\label{thm:KD-main}
    Let \(A\) be a Fr\'echet algebra. We have a natural equivalence \[\K^\an(\Cpt(A)) \simeq \KD(A)\] of spectra.
\end{theorem}
\begin{proof}
    For a unital Fr\'echet algebra \(A\) with the precompact bornology, we have \[C^\infty(M, \mathsf{Cpt}(A)) \cong \mathsf{Cpt}(C^\infty(M) \haotimes_\pi A)\] by \cite[Theorem 1.87]{Meyer:HLHA} so that we have a natural equivalence of functors
    \[\Alg_1(\mathsf{Fre}) \to \mathsf{Fun}(\Delta^{\op},\mathsf{Rings}_1), \quad A \mapsto U_b(\Delta^\an(\Cpt(A))), \text{ and } A \mapsto U_t(\Delta^{\mathrm{lc}}(A)),\] where \(U_t\) and \(U_b\) are the extensions of the functors that forget the topology and bornology, respectively. As a consequence, we have a natural equivalence \(\K^\an(A) = |\K(\Delta^{\an}(\Cpt(A))| \simeq |\K(\Delta^{\mathrm{lc}}(A))| =  \KD(A)\) of spectra. The extension of this equivalence to non-unital Fr\'echet algebras is now formal.
\end{proof}

We now compare homotopy algebraic \(K\)-theory in a manner similar to \cite{cortinas2008comparison}. Let \(X \colon \BAlg \to \Sp\) be a functor, \(\iota : A \to \Delta^{\mathrm{alg}}(A)\) and \(j: A \to \Delta^{\mathrm{diff}}(A)\) be the natural maps to the algebraic and smooth simplices. Then there is a commuting diagram 
\begin{equation}
\begin{tikzcd}\label{eq:KH=KD}
    X(A) \arrow{r}{} \arrow{d}{X(j)} & X(\Delta^{\mathrm{alg}}(A)) \arrow{d}{X(\Delta^{alg}(j))} \\
    X(\Delta^{\mathrm{diff}}(A)) \arrow{r}{X(\iota(\Delta^{\mathrm{diff}}))} & X(\Delta^{\mathrm{alg}}\Delta^{\mathrm{diff}}(A)),
\end{tikzcd}
\end{equation} in spectra. 

\begin{proposition}\label{prop:Higson}
    Let \(X \colon \BAlg \to \bC\) be a functor with values in a stable \(\infty\)-category, such that \(X(\Delta^{\mathrm{alg}}(-))\) is split exact, \(M_2\)-stable and \(\mathcal{L}^p\)-stable for \(p>1\). Then \(X(\Delta^{\mathrm{alg}}(A)) \simeq X(\Delta^{\mathrm{diff}}(A))\).
\end{proposition}

\begin{proof}
    By making obvious modifications of the proof of \cite[Lemma 4.3.1]{cortinas2008comparison} the bottom horizontal map of (\ref{eq:KH=KD}) is an equivalence with no assumptions on \(X\). The point is to observe that since for each \(n\), the algebraic simplex   \(\Delta^{\mathrm{alg}}_n(\C) = \C[x_0,\dotsc, x_n]/\{\sum_{i=0}^n x_i - 1\}\) has the fine bornology and is therefore complete, we have \(\Delta_n^{\mathrm{alg}}(\C) \haotimes A \cong \Delta_n^{\mathrm{alg}}(\C) \otimes A \cong A[x_0,\dotsc, x_n]/\{\sum_{i=0}^n x_i - 1\} = \Delta_n^{\mathrm{alg}}(A)\) for any complete bornological algebra \(A\). Furthermore, for each \(n\), \(m\), we have \[\Delta_n^{\mathrm{diff}}\Delta_n^{\mathrm{alg}}(A) \cong \Delta_m^{\mathrm{alg}}\Delta_n^{\mathrm{diff}}(A),\] so that it only remains to show that the smoothly homotopy invariant functor \(X(\Delta^{\mathrm{diff}}(-))\) maps polynomially homotopic maps to equivalences. The required smooth homotopy between \(s_0d_0 \colon A[t] \to A[t]\) and the identity is given by \(H \colon A[t] \to C^\infty[0,1] \otimes A[t]\), \(Hf(s,t) = f(st)\). To see that the functor \(X(\Delta^{\mathrm{alg}}(-))\) satisfies smooth homotopy invariance, we observe that the argument of \cite[Lemma 3.26]{CMR} (see also \cite[6.1.7]{cortinas2008comparison}) applies to stabilisations by \(p\)-Schatten class operators for \(p>1\) to get the desired result. 
\end{proof}

\begin{corollary}\label{cor:KH=KU}
    For \(p>1\), we have an equivalence \[\KH(\mathcal{L}^p(A)) \simeq \K^\an(\mathcal{L}^p(A)).\] In particular, we have an equivalence \(\mathbf{KU} \simeq \KH(\mathcal{L}^p) \simeq \K^\an(\mathcal{L}^p)\) of ring spectra.  
\end{corollary}
\begin{remark}
    In this remark, we mention how to view $\KH(\mathcal{L}^p)$ and $\K^\an(\mathcal{L}^p)$ as ring spectra. We only discuss the lax monoidality of \(\K^\an\) as the argument for \(\KH\) is similar: since the \(K\)-theory spectrum is lax monoidal, we have a natural map \[\K(\Delta^\an(A)) \wedge \K(\Delta^\an(B)) \to \K(\Delta^\an(A) \haotimes \Delta^\an(B))\] of spectra. The completion functor \(\Delta^\an(A) \otimes \Delta^\an(B) \to \Delta^\an(A) \haotimes \Delta^\an(B) \cong \Delta^\an(A \haotimes B)\) induces a natural map \(\K(\Delta^\an(A) \otimes \Delta^\an(B)) \to \K(\Delta^\an(A \haotimes B))\) of spectra, which then induces the desired lax monoidal structure \[\KD(A) \wedge \KD(B) \to \KD(A \haotimes B).\] 
    
    Note that we have a map of operator ideals 
    \[
    \mathcal{L}^p \haotimes \mathcal{L}^p \rightarrow \mathcal{L}^p,
    \]
    as described in \cite[Example 6.2]{Perrot2008}. Using the lax monoidality of $\KD$, we get a morphism of spectra
    \[
    \K^\an(\mathcal{L}^p) \wedge \K^\an(\mathcal{L}^p) \rightarrow \K^\an(\mathcal{L}^p \haotimes \mathcal{L}^p).
    \]
    Since $\mathcal{L}^p$ is an operator ideal, the ring structure on $\KD$ follows now from \cite[Proposition 3.16]{CMR}.
\end{remark}
\begin{proof}
    For the first part, we apply Proposition \ref{prop:Higson} to \(X = \K(\mathcal{L}^p\haotimes -)\), using that homotopy algebraic \(K\)-theory is \(M_2\)-stable and split exact. For the second part, the second equivalence follows from \cite[Theorem 6.2.1(2)]{cortinas2008comparison}. The main part left to argue is the equivalence between \(\KU\) and \(\KH(\mathcal{L}^p)\). To see this, we first note that the ideal inclusion \(i: \mathcal{L}^p \subseteq \csK\) induces a map 
    \begin{equation}\label{eqref:KH-comparison}
    \KH(\mathcal{L}^p) \to \KH(\csK)
    \end{equation} of spectra. Since the ring structure on both sides of \Cref{eqref:KH-comparison} comes from the tensor product of the operator ideals, we conclude that loc.cit. is a ring map (compare \cite[Lemma 3.2.2.6]{HA}).
    Furthermore, since \(\csK\) being a Fr\'echet operator ideal is \(\mathbf{K}^{\mathrm{inf}}\)-regular by \cite[Theorem 6.5.3(i)]{cortinas2008comparison}, the secondary Chern character \[\mathbf{K}^{\mathrm{nil}}(\csK) \to \tau\mathbf{K}^{\nil}(\csK) \simeq \mathbf{HC}(\csK)[-1]\] is an equivalence, where the left hand side is the fibre of the canonical map \(\mathbf{K}(-) \to \mathbf{KH}(-)\), and the right hand side is a shift of cyclic homology. Since the cyclic homology of \(\csK\) vanishes by \cite[Theorem 8.2.5(i)]{cortinas2008comparison}, \(\mathbf{K}^{\mathrm{inf}}(\csK) \simeq 0\), so that \(\mathbf{K}(\csK) \simeq \mathbf{KH}(\csK)\), and the algebraic \(K\)-theory of \(\csK\) is equivalent to \(\mathbf{KU}\). The right hand side of Equation \ref{eqref:KH-comparison} is thus equivalent \(\KU\). 
    
    It remains to show therefore that the maps induced by \eqref{eqref:KH-comparison} at the level of homotopy groups \(KH_n(\mathcal{L}^p) \to KH_n(\csK)\) are isomorphisms for all \(n \in \Z\). As in the proof of \cite[Theorem 6.5.3 (ii)]{cortinas2008comparison}, since \(\mathcal{L}^1\) multiplies \(\mathcal{L}^p\), \(KH_*(\mathcal{L}^p)\) is a graded module over the graded ring \(KH_*(\mathcal{L}^1) \cong \Z[u,u^{-1}]\), which is used in loc.cit to show that \(KH_*(\mathcal{L}^p)\) is \(2\)-periodic. For exactly the same reason, \(KH_*(\csK)\) is \(2\)-periodic, so we are reducing to proving that \(KH_*(\mathcal{L}^p) \to KH_*(\csK)\) is an isomorphism for \(* = 0\), \(-1\).  Since the ideals \(\mathcal{L}^p\) and \(\csK\) are \(K^{\mathrm{inf}}\)-regular, it follows that the natural maps \(K_*(\csK) \to KH_*(\csK)\) and \(K_*(\mathcal{L}^p) \to KH_*(\mathcal{L}^p)\) are isomorphisms for \(* = 0\), \(-1\). By naturality, we have a commuting diagram 
    \[
    \begin{tikzcd}
        K_0(\mathcal{L}^p) \arrow{r}{K_0(i)} \arrow{d}{\cong} & K_0(\csK) \arrow{d}{\cong} \\
        KH_0(\mathcal{L}^p) \arrow{r}{} & KH_0(\csK).
    \end{tikzcd}
    \] By \cite[Theorem 6.4.1]{cortinas2008comparison}, we have isomorphisms \(K_0(\mathcal{L}^p) \cong \Z\) and \(K_0(\csK) \cong \Z\) given by taking the rank of a projection. The top horizontal map is then given by mapping any rank \(1\)-projection \([e] \in K_0(\mathcal{L}^p)\) to \([i(e)] \in K_0(\csK)\), which is an isomorphism, so that \(KH_0(\mathcal{L}^p) \to KH_0(\csK)\) is an isomorphism by two-out-of-three. Finally, \(K_{-1}(\csK) = KH_{-1}(\csK) = 0\) and \(K_{-1}(\mathcal{L}^p) = KH_{-1}(\mathcal{L}^p) = 0\) by \cite[Theorem 6.4.1]{cortinas2008comparison}, so that we again have \(KH_{-1}(\mathcal{L}^p) \cong KH_{-1}(\csK) \cong 0\). 
\end{proof}

Consider the stabilisation \(\cK_p = \cK_{\cS} \haotimes \cK_{\mathcal{L}^p}\) for any \(p\geq 1\). 

\begin{corollary}\label{cor:KD=KU}
    For \(p > 1\), we have an equivalence \(\KU \simeq \KD(\cK_p) \simeq \K^\an(\Cpt(\cK_p))\) of ring spectra. 
\end{corollary}

\begin{proof}
    As in the proof of \cite[Proposition 5.2.1]{cortinas2008comparison}, the natural map \(M_2 \to \cK_{\cS}\) induces a smooth homotopy equivalence \[j: \cJ \haotimes M_2 \to \cJ \haotimes \cK_{\cS}\] for any Fr\'echet operator ideal \(\cJ\). Since \(\KD\) is smoothly homotopy invariant and \(M_2\)-stable, we have an equivalence \[\KU \simeq \KD(\cJ) \to \KD(\cJ \haotimes \cK_{\mathcal{S}})\] of ring spectra by Corollary \ref{cor:KH=KU}. Setting \(\cJ = \mathcal{L}^p\) for any \(p > 1\) yields the desired result.
\end{proof}

\subsection{Corepresentability of diffeotopy \(K\)-theory}

Consider the functor \(K_0^\an(- \haotimes \cK_?) \colon \BAlg \to \mathsf{Ab}\). By Theorem \ref{thm:KH-theory-properties}, it is smoothly homotopy invariant and excisive. When \(\cK_? = \cK_{\mathcal{S}}\), it is \(\cK_{\mathcal{S}}\)-stable (and \(M_n\)-stable for all \(n \in \N \cup \{\infty\}\)) by \cite[Lemma 3.38]{CMR}. When \(\cK_{?} = \mathcal{L}^p\), it is \(\cK_p\)-stable by \cite[Proposition 7.16]{CMR}. Consequently, by \cite[Proposition 7.26]{CMR}, it factors through \(\KKzero_{^?}\), yielding a natural map \[\beta_{0,A}^{\cK_?} \colon {\KKzero}(\C, A) \to \mathrm{Hom}_{\mathsf{Mod}_{\Z}}(K_0^\an(\cK_?), K_0^\an(A \haotimes \cK_?)).\] When \(\cK_? = \cK_p\), we may use Corollary \ref{cor:KD=KU} to identify \(\Z \cong \pi_0(\KU) \cong K_0^\an(\cK_p)\), which induces a map 
\[
\beta_A^{\cK_p} \colon \KKzero(\C,A) \to K_0^\an(\cK_p(A)).
\]

When \(\cK_? = \cK_{\mathcal{S}}\), consider the inclusion \(l: \C \to \cK_{\mathcal{S}}\). Then composition with \(\Z \cong K_0^\an(\C) \to K_0^\an(\cK_{\mathcal{S}})\) induces a homomorphism \(\mathrm{Hom}_{\mathsf{Mod}_{\Z}}(K_0^\an(\cK_{\mathcal{S}}), K_0^\an(\cK_{\mathcal{S}}(A)) \to \mathrm{Hom}_{\mathsf{Mod}_{\Z}}(K_0^\an(\C), K_0^\an(\cK_{\mathcal{S}}(A)) \cong K_0^\an(\cK_{\mathcal{S}}(A)).\) Denote the composition with \(\beta_{0,A}^{\cK_{\mathcal{S}}}\) by \[\beta_{A}^{\sK} \colon \KKzero(\C,A) \to K_0^\an(\sK(A)).\]

\begin{theorem}\label{thm:comparision-kd}
    The natural map \(\beta_{A}^{\cK^?} \colon \KKzero_{\cK^?}(\C,A) \to K_0^\an(A \haotimes \cK^?)\) is an isomorphism. 
\end{theorem}

\begin{proof}
    The same proof as \cite[Theorem 8.2.1]{CortinasThom} goes through with obvious modifications. We only provide the proof for \(\cK^? = \sK\) as this is what we mostly use later. We first describe the map \[\beta_{0,A} \colon \KKzero(\C,A) \to \mathrm{Hom}_{\Z}(K_0^\an(\sK), K_0^\an(A \haotimes \sK))\] more concretely. Let \([e] \in K_0^\an(\sK) \cong \Z\) be the canonical generator. Then to an algebra homomorphism \(\theta \colon J^n
    \sK\to  S^n(A \haotimes \sK)\) representing an element in \(\KKzero(\C,A)\), we associate the image of \([e]\) under the map \[\delta_l^{-n} \circ \delta_c^n \circ K_0^\an(\coma{\Sigma}^n(\theta)) \circ \delta_c^{-1} \circ \delta_u^n \colon K_0^\an(\sK) \to K_0^\an(A \haotimes \cK)\] where \(\delta_c\), \(\delta_u\), \(\delta_l\) are the connecting maps associated to the cone, universal and loop extensions. 
    
    In the other direction, consider the set \(qq(\C,A)\) of quasi-homomorphisms \[\C \to M_\infty(A_+) \unrhd M_\infty(A)\] and the map \(qq(\C,A) \to K_0(A_+)\) taking a quasi-homomorphism \((e_0,e_1) \mapsto [e_0] - [e_1]\). The image of this map is precisely \(K_0(A) = \ker(K_0(A_+) \to K_0(\C))\). By the universal property of quasi-homomorphisms \cite[Definition 8.48]{CMR}, the pair \[(e_0,e_1) \colon \C \rightrightarrows M_\infty(A_+) \unrhd M_\infty(A)\] induces an algebra homomorphism \(q(\C) \to M_\infty(A)\), and thereby a map in the category \(\KKzero_{\s}(\C,A)\). This yields a map \(qq(\C,A) \to \KKzero_{\s}(\C,A)\) that factors through a group homomorphism \(K_0(A) \to \KKzero_{\s}(\C,A)\). The latter map may be described as follows: let \(e \in M_n(A)\) be an idempotent; this induces an algebra homomorphism \(\C \to M_n(A)\), \(\lambda \mapsto \lambda e\). Denote the class of this homomorphism by \(\gen{e} \in \KKzero_{\s}(\C,A)\). Consider the assignment \(\mathbf{V}(A) \to \KKzero_{\s}(\C,A)\), mapping \([e] \mapsto \gen{e}\). To see that this is well-defined, we see that \(e\) and \(f\) are equivalent idempotents, then there is an invertible \(u\) such that \(e = u^{-1}f u\), so that \(\phi_e = \mathsf{Ad}_{u} \circ \phi_f\). Now use that the canonical functor \(j\) is \(M_2\)-stable by \cite[Proposition 7.22]{CMR}, and that by \cite[Proposition 3.16]{CMR}, \(j(\mathsf{Ad}_u) = \mathrm{id}\). Then by the split exactness of the universal functor, the map above is the homomorphism \(\mathbf{V}(A) \to \KKzero_{\s}(\C,A)\) \([e] \mapsto \gen{e}\) of monoids, which by the universal property of group completion induces a map \[K_0(A) \to \KKzero_{\s}(\C,A), [e] - [e'] \mapsto \gen{e} - \gen{e'}\] of groups. By split exactness, we may extend the map above to non-unital algebras. Replacing \(A\) by \(\coma{\Sigma}^n S^n(\sK(A))\) and taking colimits, we get a natural map \begin{multline*}\alpha_{0,A} \colon K_0^\an(\sK(A)) \to \KKzero_{\s}(\C,\Sigma^n S^n (\sK(A))) \\
    \cong \KKzero_{\s}(\C,\sK(A)) \cong \KKzero_{\s}(\C,A).\end{multline*} The same computation as in \cite[Theorem 8.2.1]{CortinasThom} now shows that \(\alpha_{0,A}\) and \(\beta_{0,A}\) are mutual inverses. 
\end{proof}

\begin{remark}
    Combining \ref{thm:KD-main} with \cite[Remark 8.3.3]{Cortinas2011}, we see that \(K^\an_n(\sK) \cong kk_n^{lc}(\C,\C)\), where the right hand side is Cuntz' bivariant \(K\)-theory for locally convex algebras defined in \cite{cuntz2004bivariant}. We will return to the relation between our bivariant \(K\)-theory and other topological \(K\)-theories in the next subsection \ref{subsec:remarks}.
\end{remark}

By Theorem \ref{thm:KH-theory-properties}, the functor \(\K^\an \colon \BAlg \to \mathbf{Sp}\) is smoothly homotopy invariant, excisive. It is also reduced as algebraic \(K\)-theory is reduced. By \cite[Proposition 7.16]{CMR}, the functor \(\K^\an(- \haotimes \cK_?)\) is  \(\cK_?\)-stable, so that by Proposition \ref{prop: Fun lex KKCMR C}, it extends to a left exact functor on \(\KK_?\).

\begin{corollary}\label{cor:corepresentable}
    The functors \(\KK_?(\C,-)\) and \(\K^\an (- \haotimes \cK^?) \colon \KK_{?} \to \mathbf{Sp}\) are equivalent. 
\end{corollary}
\begin{proof}
    We first show that the corepresented functor \(\KK_{?}(\C,-) \colon \KK_{?} \to \Sp\) is equivalent to \(\K^\an( - \haotimes \cK^?) \colon \KK_{?} \to \Sp\). By the stable Yoneda lemma, we have an equivalence 
\[
\Hom_{\Fun(\KK_?,\Sp)}(\KK_?(j(\C),-), \K^\an(-\haotimes \cK^?) \simeq \Omega^\infty\K^\an(\cK^?),
\] leading to a commuting diagram 

\[
\begin{tikzcd}
\pi_0(\Hom_{\Fun(\KK_?,\Sp)}(\KK_?(\C,-), \K^\an(-\haotimes \cK^?)) \arrow{r}{\simeq} \arrow{d}{\pi_0} & \pi_0(\Omega^\infty \K^{\an}(\cK^?)) \arrow{d}{\pi_0} \\
\Hom_{\Fun(\KKzero,\mathrm{Ab})}(\KKzero(\C,-), K^\an_0(-\haotimes \cK^?)) \arrow{r}{\simeq} & K^\an_0(\cK^?),
\end{tikzcd}
\] where the bottom row is an equivalence by the Yoneda lemma, and the second vertical map is an equivalence by the fact that \(\pi_0(\K^\an(\cK^?)) = K_0^\an(\cK^?)\). As a consequence, the left vertical map is an equivalence. By Theorem \ref{thm:comparision-kd}, we have a natural isomorphism \(\beta_0 \colon \KKzero_?(\C,-) \overset{\cong}\to K_0^\an(- \haotimes \cK^?)\). The equivalence of the left vertical map now implies that there is a natural transformation \(\beta \colon \KK_?(\C,-) \to \K^\an(-\haotimes \cK^?)\) such that \(\pi_0(\beta) = \beta_0 \colon \KKzero_?(\C, -) \overset{\cong}\to K_0^\an(-\haotimes \cK^?)\). The conclusion now follows from the exactness of the functors \(\KK_?(\C,-)\) and \(K^\an(- \haotimes \cK^?)\), and by replacing \(A\) by \(S^n(A)\).
\end{proof}
\begin{corollary}\label{cor:corep2}
The functors \(\bKK_{?}(\C,-)\), \(\K^\an(- \haotimes \cK^?) \colon \bKK_{?} \to \mathbf{Mod}_{\K^{\an}(\cK^?)}\) are equivalent. 
\end{corollary}

\begin{proof}
    Both \(\mathcal{L}^p\) and \(\sK\) have the Grothendieck approximation property, so that  by \cite[Lemma 1.170]{Meyer:HLHA}, tensoring \(\cK^?\) with a filtered system of bornological algebras with injective structure maps remains a filtering system with injective structure maps. The result now follows from Theorem \ref{thm:KH-theory-properties}.
\end{proof}
\begin{corollary}\label{cor: endo of unit in the diff KKs}
    We have equivalences of \emph{ring} spectra
    \[
    \KU \simeq \KK_p(\C,\C) \simeq \bKK_p(\C,\C), \qquad p>1,
    \]
    and 
    \[
    \K^{\operatorname{diff}}(\cK_\s) \simeq \KK_{\s}(\C,\C) \simeq \bKK_{\s}(\C,\C) \simeq \nKK(\C,\C).
    \]
\end{corollary}

\begin{proof}
    Follows by combining Corollary \ref{cor:KD=KU}, \ref{cor:corepresentable} and  Theorem \ref{thm:bousfield-nuc}.
\end{proof}

We end this section by mentioning some more structural properties of the functors corepresented by $\C$ in $\bKK_p$ and $\nKK$.
\begin{proposition}\label{cor: bKK K connection}
    The functors
    \[
    \operatorname{map}_{\bKK_p}(\C,-): \bKK_p \rightarrow \mathbf{Mod}_{\KU}
    \]
    and
    \[
    \operatorname{map}_{\nKK}(\C,-): \nKK\rightarrow \mathbf{Mod}_{\K^{\operatorname{diff}}(\cK_\s)}
    \]
    preserve all limits and all colimits.
\end{proposition}
\begin{proof}
    The assertion about the preservation of limits is clear since the functors under consideration are representable. The colimit preservation claim follows from Proposition \ref{prop:presentable-stable}(3) and Theorem \ref{thm:universal-nuc}.
\end{proof}

\subsection{Some remarks about topological \(K\)-theory}\label{subsec:remarks}




    Recall that by \cite[Theorem 7.33]{CMR}, for \(p>1\), there is an isomorphism \[\KKzero_0(\C,A) \cong K_0(\mathcal{L}^p(A)),\] where the right hand side is \emph{defined as} topological \(K\)-theory for a bornological algebra. This definition is justified by virtue of the fact that if \(A\) is a local Banach algebra (i.e., a directed union of Banach algebras), then it coincides with the topological \(K\)-theory of \(A\). By applying the suspension functor, we get \[\KKzero_n(\C,A) \cong K_n^{top}(A)\] for all local Banach algebras \(A\) and \(n \in \Z\). 
    
    We note that there are several other abelian group-valued functors referred to in the literature as topological \(K\)-theory, defined for more general topological algebras. In what follows, we summarise known comparisons, though none of them seem to be useful towards providing a general spectrum level equivalence between topological \(K\)-theory for Banach algebras and \(\KK(\C,-)\). 
    
    \begin{enumerate}
        \item In \cite{phillips1991k}, Phillips defines a topological \(K\)-theory for Fr\'echet algebras, which recovers topological \(K\)-theory for Banach algebras. But it is unclear whether Phillips' functor may be upgraded to a spectrum-valued functor that coincides with the spectrum-valued topological \(K\)-theory functor for Banach algebras.
        \item In \cite{cuntz1997bivariant}, Cuntz defines topological \(K\)-theory for pro-Banach algebras (or \(m\)-locally convex algebras in the sense of loc. cit) as the functor \(kk_0^{lc}(\C,-)\) corepresented by his bivariant \(K\)-theory. By \cite[Theorem 7.4]{cuntz1997bivariant}, there is a natural isomorphism \(kk_*^{top}(\C,A) \cong K_*^{top}(A)\) for Fr\'echet algebras, where the right hand side is Phillips' \(K\)-theory. In a manner analogous to the construction of \(\KK\), one could construct a stable \(\infty\)-category \(\mathbf{kk}^{top}\) with homotopy category \(kk^{top}\). One could then deduce that for a representation of an \(m\)-locally convex algebra \(A \cong \lim A_i\) as a cofiltered limit, \[\mathbf{kk}^{top}(\C, A) \simeq \lim \mathbf{kk}^{top}(\C,A_i),\] but such a functor would not interact well with filtered colimits.
        
        \item In \cite{cuntz2004bivariant}, Cuntz defines a bivariant \(K\)-theory \(kk^{lc}\) for all locally convex algebras. By \cite[Remark 8.3.3]{Cortinas2011}, one has \(kk_*^{lc}(\C,-) \cong KD_*(\sK \haotimes -)\). If one promotes Cuntz' \(kk^{lc}\) to a stable \(\infty\)-category as discussed above, we expect an equivalence \(\mathbf{kk}^{lc}(\C,-) \simeq \KD(\sK \haotimes -)\). It seems plausible that the right hand side is equivalent to a nonconnective version of the topological \(K\)-theory functor constructed in \cite{aretz2026super}, but we have not investigated this.
    \end{enumerate}

\section{Applications}

\subsection{The bootstrap class}
We now introduce the bootstrap class of $\nKK$ and \(\bKK_?\). The following result, which is well-known in the case of $C^*$-algebras (see \cite[Section 9]{bunke2024kk}), has been elusive in the bornological setting due to the absence of colimits in \(\KK\). We include some details which are also applicable in the $C^*$-algebraic context.

\begin{proposition}\label{prop:bootstrap}
    \hfill
    \begin{enumerate}
        \item The functor
        \[
        R:=\operatorname{map}_{\nKK}(\C,-): \nKK \rightarrow \mathbf{Mod}_{\K^\an(\cK_\s))}
        \]
        has a left adjoint denoted by $\iota$;
        \item $\iota$ is fully faithful and symmetric monoidal. Its essential image is called the \emph{bootstrap class}, and denoted by $\mathcal{B}^{\nKK}$;
        \item $\mathcal{B}^{\nKK}$ is the localising subcategory of $\nKK$ generated by $\C$. The Künneth formula also holds: for $A,B \in \NBA_{\kappa}$, if $B \in \mathcal{B}$, then we have an equivalence of $\K^\an(\sK)$-module spectra
        \[
        \K^\an(A \haotimes \cK_\s) \otimes_{\K^\an(\sK)} \K^\an(B \haotimes \cK_\s) \rightarrow \K^\an(A \haotimes B \haotimes \cK_\s);
        \]
        \item If $A: I \rightarrow \NBA_{\kappa}$ is a diagram with \(I\) a countable directed set, and injective bonding maps such that each $A_i$ belongs to $\mathcal{B}$, then so does $\operatorname{colim}_{\NBA_{\kappa}}A_n$;
        \item{(Universal Coefficient Theorem)} For any \(B \in \bdd^{\nKK}\), we have an equivalence \[\nKK(B,-) \simeq \mathbf{Mod}_{\K^\an(\cK_\s)}(R(B), R(-)).\]
    \end{enumerate}
\end{proposition}
\begin{proof}
    Since $R$ preserves limits and $\nKK$ and $\mathbf{Mod}(\K^\an(\sK))$ are presentable, (1) follows immediately from the adjoint functor theorem \cite[5.5.2.9(2)]{lurie2009higher}. Note that we used \Cref{cor: bKK K connection} to verify the \emph{accessibility} hypothesis of loc. cit. 
    We now show Assertion 2. We begin with the fully faithfulness assertion. We have to show that the counit of the adjunction $M \rightarrow R\iota M$ is an equivalence for all $M \in \mathbf{Mod}(\K^\an(\sK))$. Since $\mathbf{Mod}(\K^\an(\sK))$ is generated under colimits {and finite limits} by $\K^\an(\sK)$ and both $\iota$ and $R$ preserve these, it suffices to show that the counit at $\K^\an(\sK)$
    \[\K^\an(\sK) \rightarrow R \iota (\K^\an(\sK))\] is an equivalence. To that end, we note that \begin{multline*}
    \operatorname{map}_{\bKK}(\iota \K^\an(\sK),-) \simeq \operatorname{map}_{\mathbf{Mod}(\K^\an(\sK))}(\K^\an(\sK),R(-)) \\
    \simeq R(-) = \operatorname{map}_{\nKK}(\C,-),
    \end{multline*}
    whence Yoneda implies that $\C \simeq \iota \K^\an(\sK)$ in $\nKK$. Applying $R$ to this equivalence yields our assertion, as the equivalences can be directly traced to be the unit map. A similar treatment shows symmetric monoidality of $\iota$; see also \cite[Section 2.4]{Mathew2016}. 
    Assertion (3) is immediate from Assertion (2). Finally, Assertion (4) follows from the defining property of $\nKK$, see also \Cref{eqn: defining prop bkk}. Assertion (5) follows from the adjunction between \(R\) and \(\iota\), and the fully faithfulness of the latter. 
\end{proof}

\begin{remark}
    A completely analogous result holds in the category \(\bKK_{?}\), that is, the functor \(R=\operatorname{map}_{\bKK_?}(\C,-)\) has a left adjoint whose essential image is the localising subcategory \(\bKK_?\) generated by \(\C\). When \(? = p\), in light of \ref{cor:KD=KU}, the functor \(R\) takes values in \(\mathbf{Mod}_{\KU}\). The only difference is that the symmetric monoidality claims no longer hold.

    We also mention that the bootstrap classes of $\bKK_\s$ and $\nKK$ are naturally identified, since both are the localizing subcategory generated by $\C$.
\end{remark}

In what follows, we show that several interesting classes of algebras lie in the bootstrap class of \(\bKK\).

\subsubsection{Smooth functions on a manifold}

\begin{theorem}\label{thm: Cc in B}
    For a smooth manifold $M$ (with or without boundary), the algebra $C_c^{\infty}(M)$ belongs to \(\mathcal{B}^{\nKK}\) (and $\mathcal{B}^{\bKK_{?}}$).
\end{theorem}
The proof will be carried out in several stages. We first need a preparatory lemma from the general theory of smooth manifolds.
\begin{lemma}\label{lem: technical manifold}
    Let $M$ be compact, connected smooth manifold of dimension n with or without boundary. Then, there exists a handle decomposition of $M$ as 
    \[
    D^n=M_0\subset M_1 \subset ... \subset M_j=M.
    \]
    For all $0 \leq i < j$ , we have 
    \begin{equation}\label{eqn: functions decomposition}
        C_{\operatorname{flat}}^{\infty}(M_{i+1},M_i) \cong C_{\operatorname{flat}}^{\infty}(D^k \times D^{n-k},\partial(D^k) \times D^{n-k}),
    \end{equation}
    \footnote{The notation $C^{\infty}_{\operatorname{flat}}(M,N)$ denotes smooth functions on $M$ which vanish and are flat on $N$. This should not be confused the notation used in \Cref{ex: ex of ext}.1! }where $k$ depends on $i$. Furthermore, the extension of Fréchet algebras for all $0 \leq i <j$
    \begin{equation}\label{eqn: mfld extn}
         0 \rightarrow C^{\infty}(M_{i+1},M_i) \rightarrow C^{\infty}(M_{i+1}) \rightarrow C^\infty(M_i) \rightarrow 0
    \end{equation}
    possesses a bounded linear split. 
\end{lemma}
\begin{proof}
    The first two results follow from the general theory of handle decompositions, as for example is discussed in Chapters 6 and 7 of \cite{KOS}. The splitting assertion follows from a standard partition of unity argument applied to Seeley's extension theorem (\cite{SEELEY}), see also Theorem A in \cite{ALEXANDER} and the paragraph following it.
\end{proof}
\begin{proposition}\label{prop: mfld bootstrap}
    For any compact $n$-manifold $M$ with or without boundary, \(C^{\infty}(M) \in \bdd^{\nKK}\) (and \(\mathcal{B}^{\bKK_{?}}\)).
\end{proposition}

\begin{proof}
    It suffices to show that result for connected manifolds. Taking off from \Cref{lem: technical manifold}, we first note that (we are implicitly applying the precompact bornology functor to the context of \Cref{lem: technical manifold}) $C_{\operatorname{flat}}^{\infty}(M_{i+1},M_i) \in \mathcal{B}^{\bKK_{?}}$ (and \(\bdd^{\nKK}\)). Indeed, this follows from 
    \[
    C_{\operatorname{flat}}^{\infty}(M_{i+1},M_i) \cong C_{\operatorname{flat}}^{\infty}(D^k \times D^{n-k},\partial(D^k) \times D^{n-k}) \cong C^{\infty}(D^{n-k},C^{\infty}_{\operatorname{flat}}(D^k,\partial(D^k)))
    \]
    where the first isomorphism is \Cref{eqn: functions decomposition} and the second one is a standard currying isomorphism. Since $D^{n-k}$ is smoothly contractible, the claim now follows from the observation that \[C_{\operatorname{flat}}^{\infty}(D^k,\partial D^k) \cong \Omega^n \C,\] and the latter belongs to $\mathcal{B}^{\nKK}$ (and \(\bdd^{\bKK_?}\)) since it is closed under (co)fibers.

    Since the inclusion $C_{\operatorname{flat}}^{\infty}(M_{i+1},M_i) \hookrightarrow C^{\infty}(M_{i+1},M_i)$ is a smooth homotopy equivalence as shown in (the proof of) \cite[Proposition 8.4]{puschnigg2003diffeotopy}, the proof is now complete by a direct induction argument in light of \Cref{lem: technical manifold}.
\end{proof}
\begin{corollary}\label{cor: mfld bootstrap}
     For any compact $n$-dimensional manifold $M$ with (possibly empty) boundary $\partial M$, the algebra $C^{\infty}(M,\partial M) \in \mathcal{B}^{\nKK_{?}}$ (and \(\bdd^{\bKK_?}\)).
\end{corollary}
\begin{proof}
     {Using ideas similar to the proof of the existence of the semi-split exact sequence in \Cref{lem: technical manifold}}, we have an extension of Fréchet algebras 
     \[
     0 \rightarrow C^{\infty}(M,\partial M) \rightarrow C^\infty(M) \xrightarrow{\operatorname{res}} C^\infty(\partial M) \rightarrow 0,
     \]
     which furthermore has a bounded linear split. Since $\mathcal{B}^{\bKK_{?}}$ is closed in 2-out-of-3, we are now done from \Cref{prop: mfld bootstrap}.
\end{proof}
\begin{proof}[Proof of \Cref{thm: Cc in B}]
 Since \(C_c^{\infty}(M) =\bigcup_{n \in \N}C^{\infty}(M_n,\partial M_n)\), we are done by \Cref{cor: mfld bootstrap}, and as $\mathcal{B}^{\bKK_{?}}$ is closed under countable colimits.
\end{proof}

\begin{corollary}\label{cor:coefficient-bootstrap}
    For any \(A \in \BAlg\), the algebra \(C_c^\infty(M,A)\) lies in the localising subcategory of \(\bKK_{?}\) generated by \(A\). Similarly, for \(A \in \NBA_{\kappa}\), \(C_c^\infty(M,A)\) lies in the localising subcategory generated by \(A\) in \(\nKK\).
\end{corollary}
\begin{proof}
    The proof of Proposition \ref{prop: mfld bootstrap} goes through with obvious modifications. 
\end{proof}

\subsubsection{Leavitt path algebras}

Let \(E = (E^0, E^1, s,r: E^1 \rightrightarrows E^0)\) be a countable graph (that is, \(E^0\) and \(E^1\) are countable sets), and \(L(E)\) its Leavitt path algebra, which we view as a bornological algebra with the fine bornology. As a countable generated \(\C\)-vector space, it is automatically nuclear with the fine bornology. 

\begin{theorem}\label{thm:graph-boot}
    The Leavitt path algebra of a countable graph belongs to the bootstrap class of \(\bKK_{?}\) and  \(\nKK\).
\end{theorem}

\begin{proof}

Consider the semi-split extension \[0 \to \cK(E) \to C(E) \to L(E) \to 0\] of \(\C\)-algebras defined in \cite[Section 5]{cortinas2021algebraic}. As the fine bornology functor is exact, the above is an extension of complete bornological algebras. By \cite[Equation 4.13]{cortinas2021algebraic}, the algebra \(\cK(E)\) isomorphic to a countable coproduct \(\bigoplus_{v \in \mathrm{reg}(E^0)} M_{P_v}\), where \(M_{P_v}\) denotes \(P_v \times P_v\)-matrices, \(P_v\) being the set of finite paths with range \(v\). By \(M_\infty\)-stability and countable sum preservation, we have that \(\cK(E)\) is \(\bKK_{?}\) equivalent to \(\bigoplus_{\N} \bkk_{?}(\C^{(\mathrm{reg}(E^0))})\). The same proof as \cite[Theorem 4.2]{cortinas2021algebraic} shows that the map \(\varphi: \C^{(E^{0})} \to C(E)\), \(\chi_v \mapsto v\) induces a \(\KK_{?}\)-equivalence; we only need to observe that all algebra homomorphisms and quasi-homomorphisms are automatically bounded in the fine bornology. Now using that \(\bkk_{?}\) preserves countable direct sums (as \(\KK_?\) preserves finite sums), we get that \(\bigoplus_{E^{(0)}}\bkk_{?}(\C)\). The cofibre sequence \(\bkk_{?}(\mathcal{K}(E)) \to \bkk_{?}(C(E)) \to \bkk_{?}(L(E))\) is therefore equivalence to the cofibre sequence 
\[
\bigoplus_{\mathrm{reg}(E)}\bkk_{?}(\C) \overset{1 - A_E^{t}}\to \bigoplus_{E^0}\bkk_{?}(\C) \to \bkk_{?}(L(E))
\] by \cite[Theorem 5.2]{cortinas2021algebraic}. Since \(\mathcal{B}^{\bKK_?}\) is closed under countable coproducts and cofibres, the conclusion follows.\qedhere

  %
\end{proof}

\comment{\subsubsection{Overconvergent analytic functions}

Let \(r>0\), and \(\C\gen{\frac{x}{r}}\) be disc algebra of radius \(r\). More concretely, 
\[
\C\gen{\frac{x}{r}} \defeq \{\sum_{n=0}^\infty c_n t^n| \sum_{n=0}^\infty \abs{c_n}r^n < \infty\}
\] with pointwise multiplication and the weighted \(l^1\)-norm \(\norm{\sum_{n=0}^\infty c_nt^n}_r = \sum_{n=0}^\infty\abs{c_n}r^n\). Consider the \emph{dagger algebra} 
\[
\C\{\frac{x}{r}\}^\dagger = \colim_{\varrho>r} \C\{\frac{x}{\varrho}\}
\] of \emph{overconvergent analytic functions} on a disc of radius \(r\).}

\subsection{The Chern-Connes character}

Recall that analytic cyclic homology is a functor \[\mathbf{HA} \colon \BAlg \to \mathsf{Ch}(\mathsf{CBorn})\] to the homotopy category of chain complexes of complete bornological \(\C\)-vector spaces satisfying smooth homotopy invariance \cite[Theorem 5.29 and 5.45]{Meyer:HLHA}, \(\cK_{\mathcal{S}}\)-stability \cite[Theorem 5.65]{Meyer:HLHA}, and excisive with respect to semi-split extensions \cite[Theorem 5.76]{Meyer:HLHA}. By composing with the dissection functor  \(\mathrm{diss} \colon \mathsf{Ch}(\mathsf{CBorn}) \to \mathsf{Ch}(\mathsf{Ind}(\mathsf{Ban}))\), we view analytic cyclic homology as a functor valued in the homotopy category of ind-Banach spaces. On the category of chain complexes over \(\mathsf{Ind}(\mathsf{Ban})\), \cite{mukherjee2022quillen} constructs a model structure, with local chain homotopy equivalences in the sense of \cite[Section 2.3]{Meyer:HLHA} as weak equivalences. In loc.cit. it is also shown that this model structure comes from the locally split exact category structure defined in \cite[Definition 2.35]{Meyer:HLHA}, which by \cite[Lemma 2.10 and 2.11]{mukherjee2022quillen} has enough projectives and injectives. Since \(\mathsf{Ind}(\mathsf{Ban})\) is presentable, by \cite[Theorem 4.5.95]{kelly2024homotopy}, there exists an induced injective model category structure, which by \cite[Theorem 4.5.98]{kelly2024homotopy} is equivalent to the projective model structure, which by \cite[4.3.58]{kelly2024homotopy} is combinatorial. It is easily also seen to be stable. As a consequence, the underlying \(\infty\)-category \(\mathbf{D}^{\mathrm{loc}}(\C)\) is a stable, presentable \(\infty\)-category.

\comment{The category \(\bD^{\mathrm{born}}(\C) \defeq \bD(\mathsf{Ind}(\mathsf{Ban}))\) is a compactly generated, presentable stable symmetric monoidal \(\infty\)-category by \cite[Theorem 4.3.58]{kelly2024homotopy} with compact unit, we get a  natural, filtered colimit-preserving functor \[\mathbf{Hom}_{\bD(\C)}(\C, -) \colon \bD^{\mathrm{born}}(\C) \to \mathbf{Sp},\] using which we may view analytic cyclic homology as a spectrum valued functor \[\mathbf{HA} \colon \BAlg \to \bD^{\mathrm{born}}(\C) \overset{\mathbf{Hom}(\C,-)}\to \mathbf{Sp}.\]}

 There is a related functor called \emph{local cyclic homology} \(\mathbf{HL} \colon \Alg(\mathsf{Ind}(\mathsf{Ban})) \to \mathbf{D}^{\mathrm{loc}}(\C)\) that is smoothly homotopy invariant, excisive and \(\mathcal{K}_{\mathcal{S}}\)-stable. By \cite[Theorem 2.13]{Meyer:HLHA}, we have an isomorphism \(\mathbf{HL}(\mathrm{diss}(A)) \cong \mathrm{diss}\mathbf{HA}(A)\) for nuclear bornological algebras (and also Fr\'echet algebras with the precompact bornology, Banach algebras with the von Neumann bornology). For the purposes of this paper, we call the composition \(\mathrm{diss} \circ \mathbf{HA}\) local cyclic homology. 

\begin{lemma}\label{lem:HA-cocontinuity}
    Let \((A_n)\) be a countable filtered system of nuclear bornological algebras, with injective structure maps, and let \(A\) be their filtered colimit. Then \(\mathbf{HA}(A) \simeq \colim \mathbf{HA}(A_n)\).
\end{lemma}

\begin{proof}
    By Proposition \ref{prop: diss and classical together adjunctions}, we have an adjunction \[L: \mathsf{Ind}(\mathsf{Ban}) \leftrightarrows \mathsf{CBorn}: \mathrm{diss}.\] By \cite[Theorem 1.139]{Meyer:HLHA}, the functor \(L\) is symmetric monoidal and by Proposition \ref{prop: abstract symm mon loc}, the functor \(\mathrm{diss}\) is lax-monoidal, so that we have an adjunction of monoid objects 
    \[
    L : \Alg(\mathsf{Ind}(\mathsf{Ban})) \leftrightarrows \Alg(\mathsf{CBorn}): \mathrm{diss}.
    \] Let \(X_i = \mathrm{diss}(A_i) \in \Alg(\mathsf{Ind}(\mathsf{Ban}))\) and \(X = \colim X_i \in \Alg(\mathsf{Ind}(\mathsf{Ban}))\). We then have \[\colim (X) \cong \colim L(X_i) \cong \colim A_i = A.\] Then by \cite[Proposition 2.11]{Meyer:HLHA}, we have an isomorphism of chain complexes \begin{multline*} 
    \mathbf{HA}(A) \cong L(\mathbf{HL}(X)) \cong \colim_i L(\mathbf{HL}(X_i)) \\ \cong \colim_i \mathbf{HA}(\colim_i\mathrm{diss}(A_i)) \cong \colim_i \mathbf{HA}(A_i)\end{multline*} in complete bornological vector spaces, where in the second isomorphism, we have used \cite[Equation 5.69]{Meyer:HLHA}. Now since \(A\) and \(A_i\) are nuclear, we have \(\mathrm{diss}(\mathbf{HA}(A)) \cong \mathbf{HL}(\mathrm{diss}(A))\) and \(\mathrm{diss}(\mathbf{HA}(A_i)) \cong \mathbf{HL}(X_i)\) by Theorem \cite[Theorem 2.13]{Meyer:HLHA}. Since by hypothesis the filtering system \((A_i)\) has monomorphic structure maps, we have 
    \begin{equation*}
    \mathrm{diss}(A) \cong \mathrm{diss}(\colim_i A_i) \cong \mathrm{diss}(\colim_i L \mathrm{diss}(A_i)) \cong \mathrm{diss} (L X) \cong X,
    \end{equation*} where in the last isomorphism, we have used \cite[Proposition 1.144]{Meyer:HLHA}, so that we have 
    \begin{equation}\label{eqref:reduced-nuc}
    \mathrm{diss}(\mathbf{HA}(A)) \cong \mathbf{HL}(\mathrm{diss}(A)) \cong \mathbf{HL}(X).
    \end{equation} By Proposition \cite[Proposition 2.28]{Meyer:HLHA}, we have an equivalence \[\mathrm{hocolim}_i \mathrm{diss}\mathbf{HA}(A_i) \simeq \colim_i \mathrm{diss}\mathbf{HA}(A_i) \simeq \colim_i \mathbf{HL}(X_i) \simeq \mathbf{HL}(X) \cong \mathrm{diss} \mathbf{HA}(A),\] where the last isomorphism follows from Equation \ref{eqref:reduced-nuc}. 
    \end{proof}

   \begin{corollary}\label{lem:HL-born}
  Local cyclic homology descends to a colimit-preserving functor on \(\nKK\). 
\end{corollary}
\begin{proof}
    As already noted, \(\mathbf{HA}\) is homotopy invariant, excisive, \(\cK^?\)-stable, and therefore so is \(\mathrm{diss} \circ \mathbf{HA}\). By Lemma \ref{lem:HA-cocontinuity}, local cyclic homology commutes with filtered colimits of nuclear bornological algebras with injective structure maps. The descent to \(\nKK\) now follows from Theorem \ref{prop:presentable-stable} and Theorem \ref{thm:universal-nuc}, respectively. 
\end{proof}

In what follows, we abusively denote by \(\mathbf{HA}\) the composition \(\mathrm{diss} \circ \mathbf{HA}\) as we restrict analytic cyclic homology to strongly nuclear bornological algebras. Denote by \(\mathbf{HA}(A,B) \defeq \mathbf{Hom}_{\bD^{\mathrm{loc}}(\C)}(\mathbf{HA}(A), \mathbf{HA}(B))\). By the universal property of \(\nKK\), we get an essentially unique map \[\mathsf{Ch} \colon \nKK(A,B) \to \mathbf{HA}(A,B)\] of spectra that we call the \emph{bivariant Chern-Connes character}. Since \(\C[0]\) is a compact object in \(\mathbf{D}^{\mathrm{loc}}(\C)\) by \cite[Lemma 2.24]{Meyer:HLHA}, we get a filtered colimit-preserving functor \(\mathrm{map}_{\bD^{\mathrm{loc}}(\C)}(\C[0],-) \colon \bD^{\mathrm{loc}}(\C) \to \Sp\), whose composition with \(\mathbf{HA} \colon \NBA_{\kappa} \to \bD^{\mathrm{loc}}(\C)\) yields a spectrum valued functor that descends to \(\nKK\). By Theorem \ref{thm:bousfield-nuc} and Corollary \ref{cor:corep2}, we get an induced natural transformation \[\mathsf{Ch} \colon \mathbf{K}^\an(- \haotimes \sK) \to \mathbf{HA}(-) \colon \nKK \to \Sp\] of spectra. We now state the main theorem of this section:

\begin{theorem}[UCT for Local Cyclic Homology]\label{thm:UCT}
    The complexified Chern-Connes character \(\mathsf{Ch}_\C \colon \mathbf{K}^{\an}(- \haotimes \sK) \otimes \C \to \mathbf{HA}(-)\) above is a natural equivalence on the localising subcategory of \(\nKK\) generated by \((\mathcal{L}^p)^\rig\) for any \(p> 1\).
\end{theorem}

The proof of Theorem \ref{thm:UCT} requires some preparation. Following \cite[Chapter 3]{Meyer:HLHA}, the \emph{spectral radius} \(\varrho(S;A)\) of a bounded subset \(S \subseteq A\) in a complete bornological algebra is defined as \[\varrho_A(S) \defeq \inf\{r>0| \bigcup_{n=1}^\infty (r^{-1}S)^n \text{ is bounded}\}\] is bounded. A bounded homomorphism \(f: A \to B\) where \(B\) is a local Banach algebra is called \emph{isoradial} if \(\varrho(S;A) = \varrho(f(S);B)\). A subset \(X \subseteq V\) is said to be \emph{locally dense} if for any bounded subset \(S\), there is a bounded disk \(T\) such that the norm closure of \(X \cap T\) in \(V_T\) contains \(S\). 

\begin{proposition}\label{prop:Lp-isoradial}
    For any \(p \geq 1\), the canonical map \((\mathcal{L}^p)^\rig \to \mathcal{L}^p\) is isoradial and has locally dense range. 
\end{proposition}

\begin{proof}
    Let \(A = \mathcal{L}^p\) and \(S\) a bounded subset in \(A^\rig\). Then there is a sequence \((D_n)\) of complete bounded discs containing \(S\), such that \(v_n \colon A_{D_n} \to A_{D_{n+1}}\) is trace-class of order \(0\) for each \(n\). Let \(D\) be the closure of the disked hull of \(S\) under the norm \(\abs{-}_p\) of \(A\). Since in particular the operator \(A_{D_0} \to A_{D_1}\) is compact, the image of \(D_0\) is relatively compact in \(A_{D_1}\). Moreover, since \(D_1\) contains the image of \(D_0\), we have \(\overline{v_0(D_0)}^{A_{D_1}} \subseteq D_1\). Now consider the bounded inclusion \(A_{D_1} \to A\). Since \(\overline{v_0(D_0)}^{A_{D_1}}\) is compact in \(A_{D_1}\), its image is in particular closed in \(A\). Consequently, we have \(D \subseteq \overline{D_0}^{\abs{-}_p} \subseteq\overline{v_0(D_0)}^{A_{D_1}} \subseteq D_1\), which induces bounded inclusions \[u: A_D \to A_{D_1} \to A_{D_2} \to A,\] wherein the second map is trace-class of order zero, making the entire composition trace-class of order zero. A straightforward argument using the linearity of multiplication in conjunction with \cite[Lemma 3.8]{Meyer:HLHA} shows that \(S\) and \(D\) have the same spectral radius. 

    Let \(\varrho = \varrho(S;A) = \varrho(S;B)\) and let \(r>t>\varrho\). We claim that \(\bigcup_{n=1}^\infty (r^{-1}S)^n\) is bounded in \(A^\rig\). Let \(0< q \leq 1\), and  for two \(q\)-Schatten class operators \(X\) and \(Y\), consider the multiplication operator \(M_{X,Y} \colon \mathcal{L}^p \to \mathcal{L}^p\), \(T \mapsto XTY\). It is easy to see that \(M_{X,Y}\) is clearly trace-class of order \(q\), with trace-class quasi-norm \(v_{q}(M_{X,Y}) \leq \norm{X}_q \norm{Y}_q\). More generally, suppose \(T_1,\dotsc, T_k\) are \(p\)-Schatten class inside the unit ball of \(\mathcal{L}^p\), then their product \(T_1 \cdots T_k \in \mathcal{L}^{p/k}\) with \(\norm{T_1\cdots T_k}_{p/k} \leq 1\). So for any \(q>0\), we may choose \(k\) such that \(qk \geq p\). Then \(\mathcal{L}^{p/k} \subseteq \mathcal{L}^q\), so that 
    \begin{equation}\label{eq:q-schatten}
    \norm{T_1\cdots T_k}_{q}\leq 1
    \end{equation}
    Now for \(n \geq 1\), consider the \(n\)-fold multiplication map \[\mu_n \colon A_D^{\haotimes n} \to A, x_1 \otimes \cdots \otimes x_n \mapsto u(x_1)\cdots u(x_n).\] Then since \(\varrho(D;A) < t\), we have \begin{equation}\label{eqref:mu-estimate}
    \norm{\mu_n} \leq Ct^n
    \end{equation} for some \(C>0\). Let \(T_r(A_D) = \bigoplus_{n\geq 1}^{l^1} A_D^{\haotimes n}\) be the \(l^1\)-direct sum equipped with the norm \(\norm{(z_n)} = \sum r^n \norm{z_n}\), where each \(z_n \in A_D \haotimes_\pi \cdots \haotimes_{\pi}A_D \) is equipped with the projective tensor product norm. Consider the bounded map \[\phi \colon T_r(A_D) \to A \quad (z_n) \mapsto \sum_{n \geq 1}\mu_n(z_n).\] We then have \begin{equation}\label{eqref:specradius}
        \bigcup_{n\geq 1} r^{-n}D^n \subseteq \phi(B_{T_r(A_D)}).
    \end{equation}
    
    It remains to show that \(\phi\) is trace-class of order zero, that is, it is trace-class of order \(q\) for every \(0< q \leq 1\). For a fixed \(q\), choose \(k \geq 1\) such that \(kq \geq p\), and let \(n = 2k + m\) for \(m \geq 1\). Since \(u \colon A_D \to A\) is of order zero, we may write \begin{equation}\label{eq:standard-u}
    u(x) = \sum_{j \geq 1} \lambda_j \varphi_j(x) T_j
    \end{equation} where \((\lambda_j)\) is rapidly decreasing, \(\norm{\varphi_j} \leq 1\) and \(\norm{T_j}_p \leq 1\). We first estimate the multiplication map \(\mu_n \colon A_D^{\haotimes n} \to A\), \(x_1 \otimes \cdots \otimes x_n \mapsto u(x_1)\cdots u(x_n)\) in the quasi-norm \(v_q\). With \(n\) and \(k\) as above, we may decompose the product \[u(x_1)\cdots u(x_n) = (u(x_1)\cdots u(x_{k}))(u(x_{k-1})\cdots u(x_{k+m}))(u(x_{k+m+1})\cdots u(x_n)).\] 
    Setting \(\alpha = (\alpha_1,\dotsc, \alpha_k)\) and \(\beta = (\beta_1, \dotsc, \beta_k)\), we may use the standard form \ref{eq:standard-u} to write \begin{equation}\label{eqref:product-schatten}
    \mu_n(a \otimes b \otimes c) = \sum_{\alpha, \beta} \lambda_\alpha \lambda_\beta \varphi_\alpha(a) \varphi_\beta(c) X_\alpha \mu_m(b) Y_\beta,
    \end{equation} where \(a = x_1 \otimes \cdots \otimes x_k\), \(b = x_{k-1} \otimes \cdots x_{k+m}\) and \(c = x_{k+m +1} \otimes \cdots \otimes x_n\), \(\lambda_\alpha = \lambda_{\alpha_1}\cdots \lambda_{\alpha_k}\), \(X_\alpha = T_{\alpha_1}\cdots T_{\alpha_k}\) (and similarly for \(\lambda_\beta\) and \(X_\beta\)). By \ref{eq:q-schatten}, we have \(\norm{X_\alpha} \leq 1\) and \(\norm{X}_\beta \leq 1\). 
    To describe the term \(X_\alpha \mu_m(b) Y_\beta\) in \ref{eqref:product-schatten}, we write singular value decompositions \begin{align*}
    X_\alpha = \sum_i s_i(X_\alpha) \ket{a_i^\alpha}\bra{b_i^\alpha} \\
    Y_\beta = \sum_j s_j(Y_\beta) \ket{c_j^\beta} \bra{d_j^\beta},
    \end{align*} where \(s_i\) and \(s_j\) in the unit ball of \(l^q\). We then have \[X_\alpha \mu_m(b) Y_\beta = \sum_{i,j} s_i(X_\alpha) s_j(Y_\beta) \gen{b_i^\alpha, \mu_m(b) c_j^\beta} \ket{a_i^\alpha}\bra{d_j^\beta}.\] By \eqref{eqref:mu-estimate}, the functionals \[\psi_{\alpha, \beta, i, j} \colon A_D^{\haotimes n} \to \C, b \mapsto \gen{b_i^\alpha, \mu_m(b) c_j^\beta}\] have norm at most \(C t^m\). Since the functionals \(\norm{\varphi_j}\), \(\norm{\psi_j} \leq 1\), the map \[a \otimes b \otimes c \mapsto \varphi_j(a) \psi_{\alpha,\beta, i,j}(b) \varphi_j(c)\] extends to a bounded linear functional on \(A_D^{\haotimes n}\) of norm at most \(C t^n\). Normalising, we may assume that this functional has norm at most \(1\). So to check that the \(l^q\)-quasi-norm of \(\mu_m\) is bounded, it remains to analyse the coefficient term \(Ct^m \lambda_\alpha \lambda_\beta s_i(X_\alpha) s_j(Y_\beta)\). Using that \(s_i(X_\alpha)\) and \(s_j(X_{\beta})\) are in the unit ball of \(l^q\), we get that \(v_q(\mu_n)^q \leq (Ct^m)^q (\sum_{j} |\lambda_j|^q)^{2k}\), so that \[v_q(\mu_n) \leq C_qr^n,\] with \(C_q = Ct^{-2k} (\sum_{j}|\lambda_j|^{q})^{2k/q}\), which is independent of \(n>2k\). For \(1\leq n \leq 2k\), a similar argument shows that \(v_q(\mu_n) \leq \sum_{j \geq 1} (\abs{\lambda_j}^q)^{n/q} = (t^{-n}\sum_{j \geq 1}\sum_{j \geq 1} (\abs{\lambda_j}^q)^{n/q})t^n\). Setting \(\tilde{C}_q = \max_{1 \leq n \leq 2k} (t^{-n}\sum_{j \geq 1}\sum_{j \geq 1} (\abs{\lambda_j}^q)^{n/q})\), and \(D_q = \max\{C_q, \tilde{C}_q\}\), we get \[v_q(\mu_n) \leq D_q t^n\] for all \(n \geq 1\). 
    
    Recalling that \(T_r(A_D)\) has the weighted norm, we see that \[\sum_{n}v_q(r^{-n}\mu_n)^q \leq C_q^q \sum_n (\frac{t}{r})^{nq} < \infty.\] This shows that \(\phi\) is \(q\)-nuclear for every \(0<q <1\), from which we conclude that it is nuclear of order \(0\).

    Since \(\phi : T_r(A_D) \to A\) is of order \(0\), it factors through the space \(\s(\N)\) of rapidly decreasing sequences. Since \(\s(\N)\) is strongly nuclear by \ref{prop:infinitely-nuc}, we have that the map \(\s(\N) \to A\) factors through \(A^\rig\). As a consequence, the subset \(\phi(B_{T_r(A_D)})\) is bounded in \(A^\rig\). Since for any \(r> \varrho(S;A)\), we have that \(\bigcup_{n \geq 1}(r^{-1}S)^n \subseteq \phi(B_{T_r(A_D)})\) by \ref{eqref:specradius}, it follows that \[\varrho(S;A^\rig) \leq \varrho(S;A).\] The converse estimate \(\varrho(S,A) \leq \varrho(S;A^\rig)\) holds since \(A^\rig \to A\) is a bounded algebra homomorphism. This shows that \((\mathcal{L}^p)^\rig\) is a local Banach algebra. Finally, local density follows from the fact that counit map \((\mathcal{L}^p)^\rig \to \mathcal{L}^p\) is surjective (even an isomorphism) as a map of vector spaces.
\end{proof}

\begin{corollary}\label{cor:L^p-isoradial-HL}
   The natural map \((\mathcal{L}^p)^\rig \to \mathcal{L}^p\) induces a local chain homotopy equivalence \[\mathbf{HA}((\mathcal{L}^p)^\rig) \simeq \mathbf{HA}(\mathcal{L}^p) \simeq \C[0]\] for any \(p \geq 1\).
\end{corollary}

\begin{proof}
    Since \(\mathcal{L}^p\) has the precompact bornology and satisfies the Grothendieck approximation property, we have that the counit map \((\mathcal{L}^p)^{\rig} \to \mathcal{L}^p\) has approximable dense range in the sense of \cite[Definition 1.171]{Meyer:HLHA} since by Proposition \ref{prop:Lp-isoradial} it has locally dense range. Combining approximable density with the isoradiality assertion of Proposition \ref{prop:Lp-isoradial}, \cite[Corollary 6.22]{Meyer:HLHA} says that the natural map \((\mathcal{L}^p)^\rig \to \mathcal{L}^p\) induces a local chain homotopy equivalence \(\mathrm{diss} \circ \mathbf{HA}((\mathcal{L}^p)^\rig) \to \mathrm{diss} \circ \mathbf{HA}(\mathcal{L}^p)\)-equivalence.  
\end{proof}

\begin{proof}[Proof of Theorem \ref{thm:UCT}]
    It suffices to check that the Chern character is an equivalence on the generating object \((\mathcal{L}^p)^\rig\) for \(p>1\). Note that since \(\nKK\) is a coreflective subcategory of \(\bKK_\s\), and since \(\C\) is strongly nuclear, we have natural equivalences \begin{multline*}
    \K^\an(\sK \haotimes (\mathcal{L}^p)^\rig) \simeq  \bKK_{\s}(\C,(\mathcal{L}^p)^\rig) \simeq \nKK(\C,(\mathcal{L}^p)^\rig) \\
    \simeq \bKK_\s(\C,\mathcal{L}^p) \simeq \K^\an(\sK \haotimes \mathcal{L}^p) \simeq \KU,
    \end{multline*} where the second and last equivalences are Theorem \ref{thm:bousfield-nuc} and Corollary \ref{cor:KD=KU}, respectively. As a consequence, the Chern character \[\nKK(\C,(\mathcal{L}^p)^\rig) \to \mathbf{HA}((\mathcal{L}^p)^{\rig}) \simeq \mathbf{HA}(\mathcal{L}^p) \simeq \C[0]\] is equivalent to the usual Chern character \(\KU \to \mathbf{HA}(\C) \simeq \C[0]\), which is an equivalence upon tensoring with (the Eilenberg-Maclane spectrum of) \(\C\) by \cite[Theorem 5.63]{Meyer:HLHA}.\qedhere
\end{proof}

\begin{remark}
    We mention related work of Solleveld \cite{solleveld2005some}, which considers a suitable subcategory of Fr\'echet \(m\)-algebras, for which the complexified Chern character from Phillips' representable \(K\)-theory to periodic cyclic homology is an isomorphism. This is adequate to cover the case of smooth compact manifolds, but as both periodic cyclic homology and representable \(K\)-theory do not commute with filtered colimits (which do not in general exist for Fr\'echet algebras), the Chern character (or even Phillips' \(K\)-theory) for Fr\'echet algebras does not extend to larger classes of topological algebras.
\end{remark}

\begin{corollary}\label{cor:smooth-chern}
    Let \(p>1\). For any smooth manifold \(M\), we have isomorphisms \begin{multline*}
    \mathsf{Ch}_{\C,n} \colon KH_n(C_c^\infty(M,\cK_p)) \otimes \C \cong K_n^\an(C_c^\infty(M,\cK_p)) \otimes \C \\
    \to \mathsf{HA}_n(C_c^\infty(M)) \cong \mathsf{HA}_n(C_0(M)) \cong \bigoplus_{j \in \Z}H_c^{n+2j}(M,\C),
    \end{multline*} for \(n = 0,1\), where the right hand side denotes compactly supported sheaf cohomology with constant coefficients \(\C\).
\end{corollary}

\begin{proof}
    The first isomorphism follows from Corollary \ref{cor:KH=KU}. By Corollary \ref{cor:coefficient-bootstrap}, \(C_c^\infty(M, \cK_p)\) lies in the localising subcategory generated by \((\cK_p)^\rig\) in \(\nKK\). The second equivalence now follows from \ref{thm:UCT}. The second equivalence follows from \cite[6.27]{Meyer:HLHA}. The last isomorphism is \cite[Theorem 8.6]{puschnigg2003diffeotopy}.
\end{proof}

\subsection{Diffeotopy \(K\)-homology for smooth manifolds}

Let \(\mathsf{Mfld}\) denote the category of smooth manifolds and smooth maps. Consider the functor \[\K^\diff \colon \mathsf{Mfld}^\op \to \Sp, \quad M \mapsto \nKK(C_c^\infty(M), (\mathcal{L}^p)^\rig);\] we call this functor \emph{diffeotopy \(K\)-homology}. Using the bivariant Chern character, we get for each \(A \in \nKK\), a natural transformation 
\[
\mathsf{Ch}_{\C}^{\vee} \colon \nKK(-,(\mathcal{L}^p)^\rig) \otimes \C \to \mathbf{HA}(-, (\mathcal{L}^p)^\rig) \simeq \mathbf{HA}(-,\C)
\] where for the last equivalence, we have used Corollary \ref{cor:L^p-isoradial-HL}.

\begin{theorem}\label{thm:homology-UCT}
    The map \(\mathsf{Ch}_{\C}^\vee\) is an equivalence on the bootstrap class.
\end{theorem}

\begin{proof}
    It suffices to see that the map in question is an equivalence on \(\C\), which is essentially the same argument as in the proof of Theorem \ref{thm:UCT}.
\end{proof}

\begin{corollary}\label{cor:locally-finite-homology}
    For any smooth manifold \(M\), we have an isomorphism 
    \[
    \mathsf{Ch}_{\C,n}^\vee \colon \pi_n(\K^\diff(M)) \otimes \C \overset{\cong}\to\mathsf{HA}^n(C_c^\infty(M),\C) \cong \bigoplus_{j \in \Z}H_{n+2j}^{\mathrm{lf}}(M,\C),  
    \] for \(n = 0,1\), where the right hand side is locally finite homology with complex coefficients.
\end{corollary}

\begin{proof}
    The first identification is Theorem \ref{thm:homology-UCT}, and the second isomorphism follows from the identification of \(\mathbf{HA}(C_c^\infty(M)) \overset{\simeq}\to \mathbf{HA}(C_0(M))\) by the isoradiality of \(C_c^\infty(M) \to C_0(M)\) and then using \cite[Theorem 8.6, (b)]{puschnigg2003diffeotopy} to identify local cyclic homology with locally finite homology. 
\end{proof}

\appendix

\section{Some useful lemmas from category theory}

We recall some generalities from (higher)-category theory that are repeatedly used in this  article.

\begin{lemma}\label{lem: bunke seven}
    Let 
\[\begin{tikzcd}
	\mathcal{C'} & \mathcal{D'}
	\arrow["L", curve={height=-12pt}, from=1-1, to=1-2]
	\arrow["R", hook, from=1-2, to=1-1]
\end{tikzcd}\]
    be an adjunction between ($\infty$)-categories $\mathcal{C}'$ and $\mathcal{D}'$, where $\mathcal{C}'$ is bicomplete and the right adjoint $R$ is fully faithful. Then,  $\mathcal{D}'$ is also bicomplete. Moreover, $R$ preserves and reflects limits, and the colimit of a diagram $D: K \rightarrow \mathcal{D}'$ is given by 
     \[
    \underset{K}{\operatorname{colim}}D \simeq L(\underset{K}{\operatorname{colim}}R(D)).
    \]
\end{lemma}\label{lem: adjunction with ff right adjoint}
\begin{proof}
     See \cite[Subsection 7.1.4]{kerodon} and \cite[Proposition 7.21]{BUNKE3}.
\end{proof}

The following result is most likely known to the experts, but we give a proof since we could not find one in the literature. We also state and prove it in the $1$-category land because that is the generality in which we will need it.

\begin{proposition}\label{prop: abstract symm mon loc}
    Let $(\mathcal{C},\otimes)$ be a closed symmetric monoidal category, and $\mathcal{D}$ be a $\mathcal{C}$-enriched reflective subcategory of $\mathcal{C}$, i.e, there exists a $\mathcal{C}$-enriched adjunction 
\[\begin{tikzcd}
	{\mathcal{C}} & {\mathcal{D}}
	\arrow["L", curve={height=-12pt}, from=1-1, to=1-2]
	\arrow["\iota", hook, from=1-2, to=1-1]
\end{tikzcd}\]
with the right-adjoint $\iota$ being fully faithful. Then, 
\begin{enumerate}
    \item If $f : B \rightarrow C$ is a morphism in $\mathcal{C}$ such that $Lf$ is an isomorphism in $\mathcal{D}$, then the same is true for the morphisms \(1_A \otimes f: A \otimes B \rightarrow A \otimes C\) for all \(A \in \mathcal{C}\).
    \item $\mathcal{D}$ has a essentially unique symmetric monoidal refinement for which $L$ is a symmetric monoidal functor.
    \item The functor $\iota$ is naturally lax monoidal.
    \item $\mathcal{D}$ is also closed.
\end{enumerate} 
\end{proposition}
\begin{proof}
    In light of the enriched Yoneda lemma, to prove Assertion 1, we have to show that for all $D \in \mathcal{D}$ the map induced by $L(1_A \otimes f)$ on the internal homs
    \[
    \underline{\operatorname{Hom}}_{\mathcal{D}}(L(A \otimes C),D) \xrightarrow{L(1_A \otimes f)^*} \underline{\operatorname{Hom}}_{\mathcal{D}}(L(A \otimes B),D)
    \]
    is an isomorphism in $\mathcal{C}$. But note that we have a series of natural isomorphisms
    \begin{multline*}
    \underline{\operatorname{Hom}}_{\mathcal{D}}(L(A \otimes C),D) \cong \underline{\operatorname{Hom}}_{\mathcal{C}}(A \otimes C,\iota D) \\
    \cong \underline{\operatorname{Hom}}_{\mathcal{C}}(A,\underline{\operatorname{Hom}}_{\mathcal{C}}(C,\iota D)) \cong \underline{\operatorname{Hom}}_{\mathcal{C}}(A, \underline{\operatorname{Hom}}_{\mathcal{D}}(LC,D))
    \end{multline*} and similarly with $C$ replaced by $B$ (note that for the last isomorphism, we used the adjunctions being $\mathcal{C}$-enriched). Finally, our hypothesis yields that 
    \[
    \underline{\operatorname{Hom}}_{\mathcal{C}}(A, \underline{\operatorname{Hom}}_{\mathcal{D}}(LC,D)) \xrightarrow[\cong]{Lf^*} \underline{\operatorname{Hom}}_{\mathcal{C}}(A, \underline{\operatorname{Hom}}_{\mathcal{D}}(LB,D))
    \]
    is an isomorphism.  Assertion 2 is now a direct consequence of \cite[Definition 3.3.2]{hinich2013dwyer}; recall that every Bousfield localisation is a Dwyer-Kan localisation at the classes of morphisms it inverts. Assertion 3 follows from the general result that the right adjoint of a strongly monoidal left adjoint functor between monoidal categories is lax monoidal. Assertion 4 is a consequence of Day's reflection theorem, see \cite[Theorem 1.2]{DAY}.
\end{proof}

\printbibliography

\end{document}